\documentclass[11pt]{article}
\usepackage[colorlinks,linkcolor=blue,citecolor=orange]{hyperref}

\usepackage[letterpaper]{geometry}
\usepackage{amsmath,amsthm,amsfonts,amssymb,booktabs}
\usepackage{enumerate,color,xcolor}
\usepackage[dvipsnames]{xcolor}
\usepackage{graphicx}
\usepackage{subcaption}
\usepackage{url}
\usepackage{todo}
\usepackage{amssymb}
\usepackage{amsthm}
\usepackage{amsmath}
\usepackage{mathrsfs}
\usepackage{bm}

\numberwithin{equation}{section}
\theoremstyle{plain}

\newtheorem{theorem}{Theorem}[section]
\newtheorem{lemma}{Lemma}[section]
\newtheorem{assumption}{Assumption}[section]
\newtheorem{proposition}{Proposition}[section]
\newtheorem{definition}{Definition}[section]

\newtheorem{corollary}[theorem]{Corollary}

\newtheorem{remark}[theorem]{Remark}

\usepackage{comment}

\newcommand{\beq}{\begin{eqnarray}}
\newcommand{\eeq}{\end{eqnarray}}
\newcommand{\beqs}{\begin{eqnarray*}}
\newcommand{\eeqs}{\end{eqnarray*}}

\makeatletter
\def\munderbar#1{\underline{\sbox\tw@{$#1$}\dp\tw@\z@\box\tw@}}
\makeatother

\begin{document}

\begin{center}
  \Large \bf Accelerating Langevin Monte Carlo Sampling: A Large Deviations Analysis
\end{center}

\author{}
\begin{center}
  {Nian Yao\,\footnote{College of Mathematics Science, Shenzhen University, 518060 Shenzhen, China;
    yaonian@szu.edu.cn}},
       {Pervez Ali\,\footnote{Department of Mathematics, Florida State University, 1017 Academic Way, Tallahassee, FL-32306, United States of America; pa22g@fsu.edu}},
           {Xihua Tao\,\footnote{College of Mathematics Science, Shenzhen University, 518060 Shenzhen, China;
     2247640315@qq.com}},
  {Lingjiong Zhu\,\footnote{Department of Mathematics, Florida State University, 1017 Academic Way, Tallahassee, FL-32306, United States of America; zhu@math.fsu.edu}}
\end{center}

\begin{center}
 \today
\end{center}

%%%%%%%%%%%%%%%%%%%%%%%%%%%%%%%%%%%%%%

%%%%%%%%%%%%%%%%%%%%%%%%%%%%%%%%%%%%%%%%%%%%%%%%%%

\begin{abstract}
Langevin algorithms are popular Markov chain Monte Carlo methods that are often used to solve high-dimensional large-scale sampling problems in machine learning. The most classical Langevin Monte Carlo algorithm is based on the overdamped Langevin dynamics. There are many variants of Langevin dynamics that often show superior performance in practice. In this paper, 
we provide a unified approach to study the acceleration of
the variants of the overdamped Langevin dynamics through the lens of large deviations theory. Numerical experiments using both synthetic and real data are provided to illustrate the efficiency of these variants.
\end{abstract}

%%%%%%%%%%%%%%%%%%%%%%%%%%%%%%%%%%%%%%%%%%%%%
\section{Introduction}\label{sec:intro}

In this paper, we are interested in sampling a distribution $\mu$ supported on $\mathcal{X}$ with
the probability density function
\begin{equation}\label{eqn:pi}
\mu(\theta)\propto\exp(-U(\theta)),\qquad \theta\in\mathcal{X},
\end{equation}
where $\mathcal{X}$ is contained in a high-dimensional space. The sampling problem \eqref{eqn:pi}
has many applications in machine learning, such as Bayesian learning (inference)
where different choices of $U(\cdot)$ in \eqref{eqn:pi} correspond to different
problems such as Bayesian linear regression \cite{hoff2009first}, Bayesian logistic regression \cite{hoff2009first}, Bayesian deep learning \cite{wang2016towards,polson2017}
and Bayesian principal component analysis \cite{dubey2016variance}.

\textit{Langevin algorithms} are core Markov chain Monte Carlo (MCMC) methods in statistics that allow one to sample from a given density $\mu(\theta)$ of interest
defined in \eqref{eqn:pi}.
A common choice of the space $\mathcal{X}$ in \eqref{eqn:pi} is the Euclidean space, e.g. $\mathcal{X}=\mathbb{R}^{d}$.
Although Langevin algorithms have a long history and are well studied in computational physics and chemistry (see e.g. \cite{laio2002escaping,EB1980,BBK1984,PBS1988} and the books \cite{LRS2010,coffey2012langevin,leimkuhler2016molecular} and the references therein), their applications to machine learning lead to a number of computational challenges and require rethinking and redesigning \cite{ma2015complete,teh2016consistency}.
The classical Langevin Monte Carlo algorithm is based on the discretization of
{\it overdamped (or first-order) Langevin dynamics} \cite{Dalalyan,DM2017,DK2017,Raginsky}:
\begin{equation}\label{eq:overdamped-2}
d\theta_{t}=-\nabla U(\theta_{t})dt+\sqrt{2}dW_{t},
\end{equation}
where $U:\mathbb{R}^{d}\rightarrow\mathbb{R}$
and $W_{t}$ is a standard $d$-dimensional Brownian motion with $\theta_{0}\in\mathbb{R}^{d}$. Under some mild assumptions on $U(\cdot)$, the diffusion \eqref{eq:overdamped-2} admits a unique stationary distribution with the density $\mu(\theta) \propto e^{-U(\theta)}$,
also known as the \emph{Gibbs distribution} \cite{chiang1987diffusion,stroock-langevin-spectrum}. In computing practice, this diffusion is simulated by considering its discretization,
and one of the most commonly used discretization scheme is the Euler–Maruyama discretization of
\eqref{eq:overdamped-2}, often known as the unadjusted Langevin algorithm in the literature; see e.g. \cite{DM2017}:
\begin{equation}\label{discrete:overdamped}
\theta_{k+1}=\theta_{k}-\eta\nabla U(\theta_{k})+\sqrt{2\eta}\xi_{k+1},
\end{equation}
where $\xi_{k}$ are i.i.d. $\mathcal{N}(0,I_{d})$ Gaussian vectors.

The first non-asymptotic result of the discretized Langevin dynamics \eqref{discrete:overdamped}
is due to \cite{Dalalyan}, which was improved soon after
by \cite{DM2017} with a particular emphasis on the dependence on the dimension $d$.
Both works consider the total variation as the distance to measure the convergence.
Later, \cite{DM2016} studied the convergence in the 2-Wasserstein distance,
and \cite{DMP2016} studied variants of \eqref{discrete:overdamped} when $U$ is not smooth.
\cite{CB2018} studied the convergence in the Kullback-Leibler (KL) distance.
\cite{DK2017} studied the convergence when only stochastic gradients are available.

In the literature, there have also been active studies of various variants
of the overdamped Langevin dynamics and the discretization schemes.
One popular Langevin dynamics is the second-order,
also known as kinetic or {\it underdamped Langevin dynamics} 
\cite{Cheng,cheng-nonconvex,JianfengLu,dalalyan2018kinetic,GGZ2,Ma2019,GGZ}:
\begin{equation}\label{eqn:underdamped}
\begin{cases}
dr_{t}=-\gamma r_{t}dt-\nabla U(\theta_{t})dt+\sqrt{2\gamma}dW_{t},\\
d\theta_{t}=r_{t}dt,
\end{cases}
\end{equation}
where $W_{t}$ is a standard $d$-dimensional Brownian motion
with $r_{0},\theta_{0}\in\mathbb{R}^{d}$.
Under some mild assumptions on $U$, the diffusion \eqref{eqn:underdamped} admits a unique stationary distribution with the density $\mu(\theta,r) \propto e^{-U(\theta)-\frac{1}{2}|r|^{2}}$, whose $\theta$-marginal distribution coincides 
with the stationary distribution of \eqref{eq:overdamped-2}.

%%%%%%%%%%%%%%%%%%%%%%%%%%%%%%%%%%%%%%%
Another popular variant is the {\it non-reversible Langevin dynamics} \cite{HHS93,HHS05,DLP2016,DPZ17,FSS20,reyLDP,GGZ,HWGGZ20}:
\begin{equation}\label{non:reversible}
d\theta_{t}=-(I+J)\nabla U(\theta_{t})dt+\sqrt{2}dW_{t},
\end{equation}
where $W_{t}$ is a standard $d$-dimensional Brownian motion with $\theta_{0}\in\mathbb{R}^{d}$, and $J$ is an {\it anti-symmetric} matrix, i.e. $J^{\top}=-J$, 
and under mild conditions, $\mu(\theta)\propto e^{-U(\theta)}$
is the unique stationary distribution of \eqref{non:reversible}.
%%%%%%%%%%%%%%%%%%%%%%%%%%%%%%%%%

Other popular variants of Langevin dynamics include
high-order Langevin dynamics \cite{Mou2019,high-order-2025,high-order-Liu-2025}, Hessian-free high-resolution dynamics \cite{HFHR,WWZ2026}, mirror Langevin dynamics \cite{Hsieh2018,Chewi2020,ZPFP2020,TaoMirror2021},
replica exchange Langevin dynamics \cite{ChenICLR2019,Deng2020},
as well as the L\'{e}vy-driven Langevin dynamics
such as fractional Langevin Monte Carlo \cite{simcsekli2017fractional}
and fractional underdamped Langevin dynamics \cite{SZTG20}, 
the decentralized Langevin algorithms 
such as decentralized stochastic gradient Langevin dynamics \cite{DistMCMC19}, EXTRA stochastic gradient Langevin dynamics \cite{GIWZ2024}
and DIGing stochastic gradient Langevin dynamics \cite{DIGing}, 
and the constrained Langevin algorithms
such as projected Langevin Monte Carlo \cite{bubeck2015finite, bubeck2018sampling,Lamperski2021,zheng2022constrained},
skew-reflected non-reversible Langevin dynamics \cite{DFTWZ2025,WTWZ2025},
proximal Langevin Monte Carlo \cite{Brosse,SR2020}, 
decentralized stochastic gradient Langevin dynamics \cite{IslamZhu2026} 
and penalized Langevin algorithms \cite{GHZ2022}.

%%%%%%%%%%%%%%%%%%%%%%%%%%%%%%%%%%%%%%%%%%%%%%%
In this paper, we consider a {\it generalized Langevin dynamics}, that is, a Markov process $\{\bm{\mathbf{z}_t}\}_{t\ge0}$ evolving in $\mathcal{X}=\mathbb{R}^{n}$, 
and satisfy the following stochastic differential equation (SDE):
\begin{equation}\label{APP}
d\bm{\mathbf{z}_{t}}={\mathbf{f}}(\bm{\mathbf{z}_{t}})dt+\sqrt{2\mathcal{D}(\bm{\mathbf{z}_{t}})}d\bm{\mathbf{W}}_t,
\end{equation}
where $\bm{\mathbf{W}}_t$ is an $n$-dimensional Brownian motion with $\mathbf{z}_{0}\in\mathbb{R}^{n}$ and $\mathcal{D}(\bm{\mathbf{z}})$ is a positive semidefinite diffusion matrix. % that determines the intensity of the diffusion process driven by the Brownian motion. When the matrix $\mathcal{D}(\bm{\mathbf{z}})=0$, we know the dynamics of the equation (\ref{APP}) is deterministic and we do not consider this degenerate case. 
%%%%%%%%%%%%%%%%%%%%%%%%%
We define the Hamiltonian
\begin{equation}\label{eqn:Hamiltonian}
H({\mathbf{z}}) = H(\theta,r) = U(\theta)+g(\theta,r),
\end{equation}
where $\mathbf{z}$ could represent $\theta$ itself, or an augmented state space such that $\mathbf{z}=(\theta,r)$ and $\theta$ is the model parameter in $\mathbb{R}^{d}$ and $r$ is a set of auxiliary variables in $\mathbb{R}^{m}$ such that $d+m=n$. 
We are interested in the case when the choices of ${\mathbf{f}}({\mathbf{z}})$ and $\mathcal{D}({\mathbf{z}})$ yield the stationary distribution $\mu\propto \exp(-H({\mathbf{z}}))$. We write ${\mathbf{f}}({\mathbf{z}})$ as follows:
\begin{equation}
{\mathbf{f}}({\mathbf{z}})=-\left[\mathcal{D}({\mathbf{z}})+\mathcal{Q}({\mathbf{z}})\right]\nabla H({\mathbf{z}})+\Gamma({\mathbf{z}}),
\end{equation}
where $\Gamma(\mathbf{z})=(\Gamma_{1}(\mathbf{z}),\ldots,\Gamma_{n}(\mathbf{z}))$ with
\begin{equation}
\Gamma_i({\mathbf{z}}):=\sum_{j=1}^{n}\frac{\partial}{\partial\mathbf{z}_j}\left(\mathcal{D}_{ij}({\mathbf{z}})+\mathcal{Q}_{ij}({\mathbf{z}})\right),
\end{equation}
where $\mathcal{Q}(\bm{\mathbf{z}})$ is an anti-symmetric curl matrix representing the deterministic traversing effects seen
in Hamiltonian Monte Carlo (HMC) procedures. %Matrices $\mathcal{D}({\mathbf{z}})$ and $\mathcal{Q}({\mathbf{z}})$ can be adjusted to attain faster convergence to the posterior distribution.
%It is easy to see that the infinitesimal generator $\mathcal{L_\tau}$ of the dynamics (\ref{APP}) is given by
%\begin{equation}\label{FF}
%\begin{split}
%\mathcal{L_\tau}=&\sum_{i=1}^{n}\sum_{j=1}^{n}\Big[-\mathcal{D}_{ij}(\mathbf{z})\frac{\partial}{\partial \mathbf{z}_j}H(\mathbf{z})-\mathcal{Q}_{ij}(\mathbf{z})\frac{\partial}{\partial\mathbf{z}_j}H(\mathbf{z})+\frac{\partial}{\partial\mathbf{z}_j}\mathcal{D}_{ij}(\mathbf{z})+\frac{\partial}{\partial \mathbf{z}_j}\mathcal{Q}_{ij}(\mathbf{z})\Big]\frac{\partial}{\partial \mathbf{z}_i}\\
%&+\sum_{i=1}^{n}\sum_{j=1}^{n}\mathcal{D}_{ij}(\mathbf{z})\frac{\partial^2}{\partial\mathbf{z}_i\partial\mathbf{z}_j}.
%\end{split}
%\end{equation}
%By using the Fokker-Planck  equation, $\mu\propto e^{-H(\mathbf{z})}$ is a stationary distribution of the stochastic differential equation (SDE) (\ref{APP}) under the following condition satisfied by $\mathcal{Q}(\bm{\mathbf{z}})$:
%\begin{equation}\label{QW}
%\sum_{i=1}^{n}\sum_{j=1}^{n}\frac{\partial^2}{\partial\mathbf{z}_i\partial\mathbf{z}_j}\Big(\mathcal{Q}_{ij}(\bm{\mathbf{z}})e^{-H(\bm{\mathbf{z}})}\Big)=0.
%\end{equation}
It is known that
$\mu\propto e^{-H(\mathbf{z})}$ is a stationary distribution of the SDE (\ref{APP}) under the following condition for $\mathcal{Q}(\bm{\mathbf{z}})$ (see e.g. \cite{ma2015complete}):
\begin{equation}\label{QW}
\sum_{i=1}^{n}\sum_{j=1}^{n}\frac{\partial^2}{\partial\mathbf{z}_i\partial\mathbf{z}_j}\Big(\mathcal{Q}_{ij}(\bm{\mathbf{z}})e^{-H(\bm{\mathbf{z}})}\Big)=0.
\end{equation}
For any Langevin dynamics that converges to a given target distribution, it can be written in the provided framework (\ref{APP}) by choosing $\mathcal{D}(\bm{\mathbf{z}})$ and $\mathcal{Q}(\bm{\mathbf{z}})$ appropriately.

\begin{itemize}
\item[(1)]
In the general framework of (\ref{APP}), by choosing $\mathcal{D}=\mathbf{I}$
and $\mathcal{Q}=\mathbf{0}$, and
$\mathcal{X} =\mathbb{R}^d$, $\mathbf{z}_t=\theta_t\in\mathbb{R}^d$, and choosing
$H(\mathbf{z}_t)=U(\theta_t)$,
we recover the overdamped Langevin dynamics \eqref{eq:overdamped-2}.
\end{itemize}

\begin{itemize}
\item[(2)]
In the general framework of (\ref{APP}), supposing that
$\mathcal{X}=\mathbb{R}^d\times\mathbb{R}^d$ and setting $\mathbf{z}_t=(\theta_t,p_t)\in\mathbb{R}^d\times\mathbb{R}^d$, and choosing
$H(\mathbf{z}_t)=H(\theta_t,p_t)=U(\theta_t)+\frac{1}{2}| p_t |^2$, and
$$\mathcal{D}=\begin{pmatrix}
\mathbf{0}  & \mathbf{0} \\
\mathbf{0} & \gamma\mathbf{I}
\end{pmatrix} \quad{\rm and} \quad
\mathcal{Q}=\begin{pmatrix}
	\mathbf{0} & \mathbf{-I}\\
	\mathbf{I} & \mathbf{0}
\end{pmatrix}. $$
Then this reduces (\ref{APP}) to the underdamped Langevin dynamics \eqref{eqn:underdamped}.
\end{itemize}

\begin{itemize}
\item[(3)]
In the general framework of (\ref{APP}), by choosing $\mathcal{D}=\mathbf{I}$
and $\mathcal{Q}=\mathbf{J}$, a constant $d\times d$ anti-symmetric matrix, i.e. $\mathbf{J}^{\top}=-\mathbf{J}$,
and $\mathcal{X}=\mathbb{R}^d$, $\mathbf{z}_t=\theta_t\in\mathbb{R}^d$, and choosing
$H(\mathbf{z}_t)=U(\theta_t)$,
we recover the non-reversible Langevin dynamics \eqref{non:reversible}.
\end{itemize}

In addition, the general framework of (\ref{APP}) includes the following models. We write $\nabla$ for the gradient operator, $\nabla^2$ for the Hessian matrix, $\nabla\cdot$ for the divergence operator, $\Delta$ for the Laplacian operator, and $\nabla^3$ for the third-order tensor of partial derivatives. Please find the more detailed definitions in Appendix~\ref{Apendix definition}. 

\begin{itemize}
\item[(4)]
In the general framework of (\ref{APP}), by choosing $\mathcal{X}=\mathbb{R}^{d}$ and
\begin{equation}\label{mirror cofficient}
H(\mathbf{z}_t)=U(\theta_t), \quad
\mathcal{D}(\theta)=[\nabla^2\phi(\theta)]^{-1},
\quad
\mathcal{Q}(\theta)=\begin{pmatrix}
	0      & e^{U(\theta)} & \cdots & e^{U(\theta)}    \\
	-e^{U(\theta)} & 0  & \ddots &\vdots  \\
	\vdots & \vdots    & \ddots   &e^{U(\theta)}\\
	-e^{U(\theta)}      & -e^{U(\theta)} &\cdots & 0
\end{pmatrix},
\end{equation}
we reduce (\ref{APP}) to the well-known {\it mirror Langevin dynamics}  \cite{Hsieh2018,ZPFP2020,Chewi2020,TaoMirror2021,Ahn2021}, that is,
\begin{equation}\label{mirror Langevin dynamics}
d\theta_t=(\varPhi(\theta_t)-[\nabla^2\phi(\theta_t)]^{-1}\nabla U(\theta_t))dt+\sqrt{2[\nabla^2\phi(\theta_t)]^{-1}}dW_t,
\end{equation}
where $\varPhi(\theta):=-[\nabla^2\phi(\theta)]^{-1}\mathrm{Tr}\left(\nabla^3\phi(\theta)[\nabla^2\phi(\theta)]^{-1}\right)$ and $W_t$ is a standard Brownian motions in $\mathbb{R}^d$ with $\phi \in C^3(\mathcal{X})$. For mirror Langevin dynamics~\eqref{mirror Langevin dynamics}, $h := \nabla \phi : \mathcal{X} \mapsto \mathbb{R}^{d}$ is known as the {\it mirror map}, which transforms the coordinates from the primal space to the dual space and $\varPhi(\theta) =  [\nabla h(\theta)]^{-1}\mathrm{Tr}\left(\nabla^2 h(\theta)[\nabla h(\theta)]^{-1}\right)$ is a mirror curvature correction that ensures that the Langevin dynamics correctly follows the structure induced by the mirror map $h$.
\end{itemize}

\begin{itemize}
\item[(5)]  In the general framework of (\ref{APP}), by choosing
$\mathcal{X} =\mathbb{R}^d\times\mathbb{R}^d\times\mathbb{R}^d$ and setting $\mathbf{z}_t=(\theta_t,p_t,r_t)\in\mathbb{R}^d\times\mathbb{R}^d\times\mathbb{R}^d$, and choosing
\begin{equation}\label{H:function:high:order}
H(\mathbf{z}_t)=H(\theta_t,p_t,r_t)=U(\theta_t)+\frac{1}{2}|p_t|^2+\frac{1}{2}|r_t|^2,
\end{equation}
and
\begin{equation}\label{high-order cofficient}
\mathcal{D}=\begin{pmatrix}
\mathbf{0} & \mathbf{0} & \mathbf{0}\\
\mathbf{0} & \mathbf{0} & \mathbf{0}\\
\mathbf{0} & \mathbf{0} & \alpha\mathbf{I}
\end{pmatrix} \quad{\rm and} \quad
\mathcal{Q}=\begin{pmatrix}
	\mathbf{0} & \mathbf{-I} & \mathbf{0}\\
	\mathbf{I} & \mathbf{0} & -\gamma\mathbf{I}\\
	\mathbf{0} & \gamma\mathbf{I} & \mathbf{0}
\end{pmatrix}.
\end{equation}
Then this reduces (\ref{APP}) to the {\it high-order Langevin dynamics} introduced in \cite{Mou2019}; see also \cite{Monmarche2023,high-order-2025,high-order-Liu-2025}:
\begin{equation}\label{JI}
\begin{cases}
d\theta_t=p_t\ dt,\\
dp_t=-\nabla U(\theta_t)\ dt + \gamma r_t\ dt,\\
dr_t=-\gamma p_t\ dt - \alpha r_t\ dt + \sqrt{2\alpha}\ dW_t,
\end{cases}
\end{equation}
where $\alpha,\gamma>0$ are friction parameters and $W_t$ is a standard Brownian motions in $\mathbb{R}^d$.
\end{itemize}

\begin{itemize}
\item[(6)] In the general framework of (\ref{APP}), let
$\mathcal{X} =\mathbb{R}^d\times\mathbb{R}^d $ and  $\mathbf{z}_t=(\theta_t,r_t)\in\mathbb{R}^d\times\mathbb{R}^d$, and define the function
\begin{equation}\label{H:function:HFHR}
H(\mathbf{z}_t)=H(\theta_t,r_t)=U(\theta_t)+\frac{1}{2}|r_t |^2,
\end{equation}
along with the matrices
$$\mathcal{D}=\begin{pmatrix}
\beta\mathbf{I} & \mathbf{0}\\
\mathbf{0} & \alpha\mathbf{I}
\end{pmatrix}\quad {\rm and} \quad
\mathcal{Q}=\begin{pmatrix}
\mathbf{0} & \mathbf{-I} \\
\mathbf{I} & \mathbf{0} \\
\end{pmatrix}, $$\\		
then we obtain the {\it Hessian-free high-resolution dynamics} introduced in \cite{HFHR}:
\begin{equation}\label{Hessian-free high-resolution}
\begin{cases}
d\theta_t=r_t\ dt - \beta\nabla U(\theta_t)\ dt + \sqrt{2\beta}\ d\bar{W}_t,\\
dr_t=-\alpha r_t\ dt - \nabla U(\theta_t)\ dt + \sqrt{2\alpha}\ d\bar{B}_t,
\end{cases}
\end{equation}
where $\beta>0,\alpha>0$ are friction parameters, and $\bar{W}_t,\bar{B}_t$ are two independent standard Brownian motions in $\mathbb{R}^d$.
\end{itemize}

%%%%%%%%%%%%%%%%%%%%%%%%%%%%%%%%%%%%%%%%%%%%%%%%%%%%%%%%%%%%%%%%%%%%%%%%%%%%%%%%%%%%%%

For any fixed time $ {t>0} $, the empirical measure of the generalized Langevin dynamics $ \{\mathbf{z}_{t}\}_{t\ge0} $ defined in (\ref{APP}) is given as
\begin{equation*}
\pi_{t}= \frac{1}{t}\int_{0}^{t}\delta_{\mathbf{z}_{s}}ds,
\end{equation*}
where $ \delta_{x} $ is the Dirac measure at $x\in \mathcal{X}$, such that $ \{\pi_{t}\}_{t\ge0} $ is a sequence of random measures, which are random elements of $\mathcal{P}(\mathcal{X})$, that is, the space of probability measures on $\mathcal{X}=\mathbb{R}^n$.
We will show that under mild conditions, \eqref{APP} is ergodic
and by the ergodic theorem,
\begin{equation}\label{ergodic:mean}
\pi_{t}\rightarrow\mu\propto e^{-H(\mathbf{z})},\qquad\text{almost surely as $t\rightarrow\infty$}.
\end{equation}
To understand how fast $\pi_{t}$ converges to the stationary distribution $\mu$,
we adopt a {\it large deviations} approach.
We present the large deviation principle about the invariant measure of the empirical measure related to the process, and then use the large deviation principle to judge the acceleration effect of the algorithm.
The large deviation principle was first formulated
in the pioneering work by Varadhan \cite{Varadhan1966}.
In contrast to the law of large numbers, which studies
a typical event, the large deviations study the small probability
of rare events \cite{DZ1998,Varadhan1,Varadhan2}.
In a series of seminal papers, Donsker and Varadhan
studied the large deviation principle for the empirical measure $\pi_{t}$
for $t\rightarrow\infty$, where the underlying is a Markov process \cite{DV1,DV2,DV3,DV4}.
Informally, $\mathbb{P}(\pi_{t}\in\cdot)$ satisfies a large deviation principle
with rate function $I(\cdot)$ if for any $\nu\in\mathcal{P}(\mathcal{X})$,
\begin{equation*}
\mathbb{P}(\pi_{t}\approx\nu)=e^{-tI(\nu)+o(t)},
\end{equation*}
as $t\rightarrow\infty$, where $I(\nu)\geq 0$ and $I(\nu)=0$ if and only if $\nu=\mu$.
In other words, when $\nu\neq\mu$, $\mathbb{P}(\pi_{t}\approx\nu)$ is exponentially
small as time $t\rightarrow\infty$, and the rate function $I(\nu)$ characterizes
how small this probability is. This suggests the larger the value of the rate function,
the smaller the probability the empirical measure deviates away from the Gibbs distribution,
and hence one expects faster convergence.

From technical perspective, although Donsker-Varadhan large deviations theory for
large-time asymptotics works for the general Markov process, it often requires restrictive assumptions,
including for example compact domain (see e.g. \cite{DV1}).
The Langevin dynamics we are interested in lives in the unbounded Euclidean space,
and our analysis relies on a more recent large deviations result
where the space can be unbounded and the underlying topology is the weighted topology \cite{LDP-GG}.

The idea of applying large deviations analysis to study Langevin algorithms and the variants
is not new. For example, the large deviations for overdamped and underdamped Langevin dynamics
are studied in \cite{LDP-GG}. 
%It is shown for example in \cite{LDP-GG} that $\mathbb{P}(\pi_{t}\in\cdot)$
%follows a large deviation principle with the rate function
%\begin{equation}\label{I:o}
%I_o(\nu)=\frac{1}{4}\int_{\chi}|\nabla\upsilon|^2\ d\nu,
%\end{equation}
%where $|\cdot|$ denotes the Euclidean norm. 
The large deviations for some variants
of Langevin algorithms have been studied in \cite{reyLDP, ChenICLR2019}.
However, these works are mostly restricted to a particular variant of the Langevin algorithm,
e.g. the non-reversible Langevin dynamics in \cite{reyLDP}, skew-reflected non-reversible Langevin dynamics in \cite{WTWZ2025} and replica exchange Langevin dynamics in \cite{ChenICLR2019}.
We take a more unified approach by studying a generalized Langevin dynamics, described in \cite{ma2015complete},
that includes overdamped, underdamped, non-reversible Langevin dynamics
as special cases,
as well as other variants from the recent literature literature, for example, the mirror Langevin dynamics \cite{Hsieh2018},
the high-order Langevin dynamics \cite{Mou2019}, 
and the Hessian-free high-resolution dynamics \cite{HFHR} whose large deviations have never been studied in the literature to the best of our knowledge.

The contributions of our paper can be summarized as follows:
\begin{itemize}
\item 
In Section~\ref{sec:main}, we first establish large deviations for the empirical measures for the generalized Langevin dynamics \eqref{APP} (Theorem~\ref{thm:main}) assuming Hypoellipticity (Assumption~\ref{Hypoellipticity}), Controllability (Assumption~\ref{irreducibility}) and Lyapunov condition (Assumption~\ref{Lyapunov condition}). 
Then, we apply the general framework of the large deviations for the generalized Langevin dynamics \eqref{APP} 
to study the large deviations for variants of Langevin dynamics, including
mirror Langevin dynamics \eqref{mirror Langevin dynamics} (Theorem~\ref{Mirror}), 
high-order Langevin dynamics \eqref{JI} (Theorem~\ref{thm:high:order}),
and Hessian-free high-resolution dynamics \eqref{Hessian-free high-resolution} (Theorem~\ref{thm: Hessian-Free}),
whose large deviations have never been established in the literature
to the best of our knowledge.
The technical novelty lies upon a careful analysis to
check Hypoellipticity (Assumption~\ref{Hypoellipticity}), Controllability (Assumption~\ref{irreducibility}) and Lyapunov condition (Assumption~\ref{Lyapunov condition}) by constructing
novel Lyapunov functions, under mild conditions on the target function $U(\cdot)$, for each variant of Langevin dynamics.
To the best of our knowledge, this is the first general and unified framework
for the large deviations of Langevin dynamics and a study of
most familiar variants of Langevin dynamics of interest in the literature.
\item
In Section~\ref{sec:comparison}, we use the large deviation rate functions
obtained in Section~\ref{sec:main} as a measure to analyze the speed
of convergence to the invariant measure. 
The larger the rate function, the more concentrated the empirical measure is around the invariant measure, thus indicating faster convergence. 
Since variants of Langevin dynamics, such as underdamped Langevin dynamics \eqref{eqn:underdamped}, high-order Langevin dynamics \eqref{JI},
and Hessian-free high-resolution dynamics \eqref{Hessian-free high-resolution}, 
may live in a higher-dimensional space such as $\mathbb{R}^{2d}$ or $\mathbb{R}^{3d}$
than the overdamped Langevin dynamics \eqref{eq:overdamped-2} in $\mathbb{R}^{d}$, a direct comparison
of their rate functions impossible.
We introduce a novel method to expand the space of overdamped Langevin dynamics to match the dimensions with those of the variants of Langevin dynamics under study. 
We then show comparisons of their rate functions for some parameter regimes
either on the whole space of probability measures or a subspace of probability measures (Proposition~\ref{Hessian-free high-resolution VS overdamped}, Proposition~\ref{underdamped VS overdamped}, Proposition~\ref{High-order VS overdamped}). In the former case, by applying contraction principle from large deviations, we show the acceleration for the underlying process $\theta_{t}$ in $\mathbb{R}^{d}$ by establishing the comparison for the LDP rate functions for its empirical measure (Corollary~\ref{cor:HFHR:vs:overdamped}).
\item 
In Section~\ref{sec:numerical}, we conduct numerical experiments 
for various Langevin algorithms based on the Euler–Maruyama
discretization of variants of Langevin dynamics. We compare
the numerical results with the unadjusted Langevin algorithm, i.e.
the Euler–Maruyama discretization of overdamped Langevin dynamics.
These numerical experiments show superior performance
or comparable performance using variants of Langevin algorithms.
\end{itemize}

The rest of the paper is organized as follows.
In Section~\ref{sec:main}, we introduce the main results
of our paper, including a large deviation principle
for the generalized Langevin dynamics
and, under various mild assumptions on the target distribution, 
the large deviations for mirror Langevin dynamics, 
high-order Langevin dynamics and Hessian-free high-resolution dynamics.
In Section~\ref{sec:comparison}, we study the acceleration of 
variants of Langevin dynamics compared to the overdamped Langevin dynamics
by comparing their rate functions from the large deviations theory.
Numerical experiments are provided in Section~\ref{sec:numerical}.
We conclude in Section~\ref{sec:conclude}.
The notations used in this paper are summarized in Appendix~\ref{Apendix definition},
all the technical proofs are provided in Appendix~\ref{Technical Proofs}, 
and additional technical lemmas are presented in Appendix~\ref{sec:technical:lemmas}.

 %%%%%%%%%%%%%%%%%%%%%%%%%%%
\section{Main Results}\label{sec:main}

\subsection{Preliminary}

%In this paper, we will show that under mild conditions, the generalized Langevin dynamics defined by SDE (\ref{APP}) is ergodic. 
%Our main interest is to study the large deviations, 
%that is the exponentially small probability that the empirical
%measure deviations away from the ergodic mean $\mu$ given in \eqref{ergodic:mean}.

We use $C_c^\infty(\mathcal{X})$ (resp. $C_b(\mathcal{X})$) to denote
the space of smooth functions with compact support
(resp. continuous and bounded functions), as well as the space of smooth functions
growing at most polynomially and whose derivatives also grow at most polynomially: 
\begin{equation}\label{S:space}
\mathscr{S}=\left\{\varphi\in C^\infty(\mathcal{X})\Big|\forall\alpha\in \mathbb{N}^d,\,\, \exists N>0,\,\,\, \text{such that}~ \sup_{x\in\mathcal{X}} \frac{|\partial^\alpha \varphi(x)|}{(1+|x|^2)^N}<+\infty \right\},
\end{equation}
where $\partial^\alpha=\partial_{x_1}^{\alpha_1}\cdots\partial_{x_d}^{\alpha_d}$ and $\alpha=(\alpha_1,\alpha_2,\cdots,\alpha_d)$.

Let us recall that given a sequence of empirical measures $\{\pi_{t}\}_{t\ge0}$, we say that $\pi_{t}$ satisfies a {\it large deviation principle} (LDP) on $\mathcal{P}(\mathcal{X})$ equipped with the $\tau^\kappa$ topology (here $\kappa:\mathcal X\to[1,\infty)$ is a measurable weight controlling tails, see Remark~\ref{rem:kappa-def}) with the rate function $I: \mathcal{P}(\mathcal{X}) \to \mathbb{R}$
if $I$ is non-negative, lower semicontinuous and for any $\tau^\kappa$-measurable set $\varTheta\subset\mathcal{P}(\mathcal{X})$,
\begin{align*}
-\mathop{\inf}\limits_{\nu\in\mathring{\varTheta}} I(\nu)\le \mathop{\liminf}\limits_{t\rightarrow\infty}\frac{1}{t}\log \mathbb{P}(\pi_t\in\varTheta)
\le \mathop{\limsup}\limits_{t\rightarrow\infty}\frac{1}{t}\log \mathbb{P}(\pi_t\in\varTheta)\le -\mathop{\inf}\limits_{\nu\in \bar{\varTheta}} I(\nu),
\end{align*}
where $\mathring{\varTheta}$ denotes the interior of $\varTheta$ and $\bar \varTheta$ stands for its closure; see e.g. \cite{DZ1998,Varadhan1,Varadhan2} for a survey on the theory of large deviations. 

\begin{remark}\label{rem:kappa-def}
Throughout, $\kappa:\mathcal X\to[1,\infty)$ denotes a measurable \emph{weight function} used to control tails and to define the $\kappa$-weak topology $\tau^\kappa$ on $\mathcal P(\mathcal X)$. We set
\[
\mathcal P_\kappa(\mathcal X):=\Bigl\{\nu\in\mathcal P(\mathcal X): \int \kappa\,d\nu<\infty\Bigr\},
\]
and define $\tau^\kappa$ as the coarsest topology making all maps
\[
\nu\ \longmapsto\ \int f\,d\nu\qquad\text{continuous for all }f\in C_b(\mathcal X)\text{ with }|f|\le C\,\kappa \text{ for some }C<\infty.
\]
Equivalently, $\tau^\kappa$ is the topology of convergence against $\kappa$-dominated continuous test functions. If $\kappa$ is bounded on $\mathcal X$, then $\tau^\kappa$ coincides with the usual weak topology. Typical choices include $\kappa(x)=1+|x|^p$ or $\kappa(x)=\exp(\alpha|x|)$, depending on the growth of the dynamics.
\end{remark}

By following \cite{LDP-GG}, we introduce the following three basic assumptions: hypoellipticity of the generator, controllability (irreducibility of the dynamics), and a Lyapunov condition.
First, let us introduce the hypoellipticity assumption.

\begin{assumption}[Hypoellipticity]\label{Hypoellipticity} functions $\mathbf{f}$ and $\sqrt{\mathcal{D}}$ in (\ref{APP}) belong to $\mathscr{S}^n$ and $\mathscr{S}^{n\times n}$,
respectively, and the generator $\mathcal{L}_\tau$ defined in (\ref{FF}) satisfies the Hypoelliptic H\"ormander condition. More precisely, $\mathcal{L}_\tau$ can be written as
\begin{equation}\label{HYP0}
\mathcal{L}_\tau=\sum_{j=1}^{r}X^*_jX_j+X_0,
\end{equation}
where $X_0,\ldots,X_r$ denote first-order homogeneous differential operators in an open set $\Omega \subset \mathbb{R}^n$ with $C^\infty$ coefficients such that $\{X_j\}_{j=1}^r$, $\{[X_i,X_j]\}_{i,j=0}^r$, $\{[[X_i,X_j], X_k]\}_{i,j,k=0}^r,\ldots$ span the space
$\Omega$ at any given point $x\in \mathbb{R}^n$ for a finite number of commutators $n_x \in \mathbb{N}$, where the definition of $[X,Y]$ can be found in Appendix~\ref{Apendix definition}.
\end{assumption}

Next, we introduce the following controllability condition
about the irreducibility of the dynamics.

\begin{assumption}[Controllability]\label{irreducibility} For any $x, y\in \mathcal{X}$ and $T>0$, there exists a control $u \in C^0([0,T],\mathbb{R}^n)$
such that the path $\phi\in C^0([0,T],\mathcal{X})$ satisfying
\begin{equation}\label{Controllability}
\begin{cases}
\phi(0)=x,\\
\dot{\phi}(t)=\mathbf{f}(\phi(t))+\sqrt{2\mathcal{D}(\phi(t))}u(t),
\end{cases}
\end{equation}
and $\phi(T)=y$ is well-defined, where $\mathbf{f}$ and $\mathcal{D}$ are defined in \eqref{APP}.
\end{assumption}

Finally, by following the Lyapunov condition given in Proposition~2.9 in \cite{LDP-GG}, we introduce the following assumption.

\begin{assumption}[Lyapunov condition]\label{Lyapunov condition}There exists a function $W:\mathcal{X}\to[1,+\infty) $ of class $C^2(\mathcal{X})$ with compact level sets and is such that $|\sqrt{\mathcal{D}}\nabla W|$ has compact level sets. For any $\theta \in (0,1)$,
\begin{align}\label{Lyapunov1}
-\mathcal{L_\tau}W-\theta|\sqrt{\mathcal{D}}\nabla W|^2\sim |\sqrt{\mathcal{D}}\nabla W|^2,
\end{align}
where for any two functions $g,f:\mathcal{X}\rightarrow\mathbb{R}$, $g$ is said to be equivalent to $f$ (denoted by $g \sim f$) if there exist constants $c, c' > 0$ and $R, R' \in \mathbb{R}$ such that
$c'g(x) - R' \leq f(x) \leq c g(x) + R$ for any $x\in\mathcal{X}$.
\end{assumption}

\begin{assumption}[Witten-Lyapunov condition]\label{Witten-Lyapunov condition}
There exists a function $W:\mathcal X\to[1,+\infty)$ of class $C^{2}(\mathcal X)$,
with compact level sets and such that
\begin{equation}\label{eq:Psi-def}
\Psi_\tau:= -\,\frac{\mathcal{L_\tau} W}{W} 
\end{equation} has compact level sets. Moreover, there exists a $C^{2}(\mathcal X)$ function
$\mathscr W:\mathcal X\to[1,+\infty)$ such that for some constants
$C_{1}>0$, $C_{2}\in\mathbb R$,
\begin{equation}\label{eq:213}
\mathscr W^{2}\,\le\, C_{1}\,W, \qquad
\Psi \,\sim\, -\,\frac{\mathcal L\mathscr W}{\mathscr W}, \qquad
-\,2\,\frac{\mathcal L\mathscr W}{\mathscr W}\,\le\,
-\,\frac{\mathcal L W}{W}\,+\,C_{2}.
\end{equation}
\end{assumption}
\begin{assumption}[Standing assumption on the weight $\kappa$]\label{ass:kappa}
Throughout, we consider an arbitrary function $\kappa:\mathcal X\to[1,+\infty)$
belonging to $\mathscr S$ such that
\begin{itemize}
  \item $\kappa \ll \Psi_\tau$ \quad (with $\Psi_\tau:=-\mathcal L_\tau W/W$ from Assumption~\ref{Lyapunov condition});
  \item either (i) $\kappa$ is bounded, or (ii) $\kappa$ has compact level sets and there exists $c\in\mathbb R$ such that
  \begin{equation}\label{eq:kappaW}
    \mathcal L_\tau(\kappa W)\ \le\ c\,\kappa W .
  \end{equation}
\end{itemize}
\end{assumption}

By Lemma~\ref{Lyapunovkappa} in Appendix~\ref{sec:technical:lemmas}, the Witten–Lyapunov and Lyapunov drift conditions hold with $W_\eta(\mathbf{z})=e^{\eta W(\mathbf{z})}$, and $\Psi_\tau\sim|\sqrt{\mathcal{D}}\nabla W|^{2}$. All results that require the Lyapunov condition can thus be applied under Assumption~\ref{Lyapunov condition}.

%%%%%%%%%%%%%%%%%%%%%%%%%%%

Under Assumptions~\ref{Hypoellipticity}-\ref{Lyapunov condition}, 
it is known that $\pi_{t}$ converges to the invariant distribution $\mu$ exponentially 
fast as $t\rightarrow\infty$. We have the following
result from \cite{LDP-GG}.

\begin{proposition}[Proposition 2.10. in \cite{LDP-GG}]
\label{prop:Wasserstein}
Suppose Assumptions ~\ref{Hypoellipticity}-\ref{Lyapunov condition} and \ref{ass:kappa} hold 
so that the Markov semigroup $(P_t)_{t\ge0}$ is well-defined and admits a unique invariant probability measure $\mu$. Let $W:\mathcal X\to[1,\infty)$ be
the Lyapunov function in Assumption~\ref{Lyapunov condition}, and denote by
\[
d_W(\nu,\eta):=\sup_{\|\varphi\|_{B_W^\infty}\le1}\Bigl|\int_{\mathcal X}\varphi\,d\nu-\int_{\mathcal X}\varphi\,d\eta\Bigr|,
\qquad 
\|\varphi\|_{B_W^\infty}:=\sup_{x\in\mathcal X}\frac{|\varphi(x)|}{W(x)} .
\]
Then there exist constants $C,c>0$ such that, for any initial measure
$\nu\in\mathcal P_W(\mathcal X)$,
\begin{equation}\label{eq:contraction-dW}
d_W(\nu P_t,\mu)\ \le\ C\,e^{-ct}\, d_W(\nu,\mu),\qquad \forall\,t\ge0.
\end{equation}
In particular, $\mu\in\mathcal P_W(\mathcal X)$ and the process is exponentially ergodic in the weighted total variation distance $d_W$, see exact definition in Appendix~\ref{Apendix definition}. 
\end{proposition}

However, Proposition~\ref{prop:Wasserstein} does not have an explicit expression
for the convergence speed $c$ except that such a positive $c$ exists. 
As a result, one cannot rely on Proposition~\ref{prop:Wasserstein} to compare
the performance of the variants of Langevin dynamics. 
This motivates us to adopt a different approach in our paper, by obtaining large deviations
for the variants of Langevin dynamics, and characterizing their rate functions
as a measure to compare their convergence to the invariant distributions.

%%%%%%%%%%%%%%%%%%%%%%%%%%%
We introduce the carr\'e du champ operator \cite{BGL2013} associated
with $\mathcal{L}:=b~\cdot\nabla+\mathcal{S}:\nabla^2$ which is defined as follows. For two regular functions $\varphi$ and $\psi$:
\begin{equation}\label{Gamma function}
\mathscr{C}(\varphi,\psi)=\frac{1}{2}\big(\mathcal{L}(\varphi\psi)-\varphi\mathcal{L}\psi-\psi\mathcal{L}\varphi\big)=\nabla\varphi\cdot\mathcal{S}\nabla\psi,
\end{equation}
where $\nabla^2$ stands for the Hessian matrix, and for two matrices $A, B$ belonging to $\mathbb{R}^{d\times d}$, we write $A : B = {\rm Tr}(A^{\top}B)$.

For any $\varphi\in C_c^\infty(\mathcal{X})$, we also introduce the seminorms
$$
|\varphi|_{\mathscr{H}^1(\nu)}^2=\int_{\mathcal{X}}\mathscr{C}(\varphi,\varphi)\ d\nu,
$$
and
$$
|\varphi|_{\mathscr{H}^{-1}(\nu)}^2=\sup_{\psi\in C_c^\infty}\left\{2\int_{\mathcal{X}}\varphi\psi\ d\nu-|\psi|_{\mathscr{H}^1(\nu)}^2\right\}.
$$
Let $\tilde{\nabla}$ denote the adjoint of gradient operator $\nabla$ in $L^2(\nu)$. 

In order to provide the exact expression of the rate function of generalized Langevin dynamics \eqref{APP}, we consider the generator $\mathcal{L}_\tau$ of the dynamics \eqref{APP}, and decompose it into symmetric and anti-symmetric parts with respect to $\mu$. First, for any closed operator $T$, we denote $T^*$ as its adjoint on $L^2(\mu)$, where $\mu$ is the invariant probability measure of the dynamics. 
Let $\mu(\mathrm dx)=Z^{-1}e^{-U(x)}\mathrm dx$ be the invariant measure. 
On $L^{2}(\mu)$, the adjoint of the gradient is
\[
\nabla^{*}=-\operatorname{div}+\nabla U\cdot,
\qquad 
\partial_i^{*}=-\partial_i+(\partial_i U).
\]
Given a symmetric diffusion matrix $ 
\mathcal D(x)=\left(\mathcal D(x)\right)^{\top}\succeq 0$ and an
antisymmetric matrix $\mathcal Q(x)=-\mathcal Q(x)^{\top}$, we write the infinitesimal generator in the
canonical form:
\begin{equation}\label{eq:canonic}
\mathcal L_\tau=\nabla^{*}(\mathcal D+\mathcal Q)\,\nabla
\sum_{i=1}^{n}\sum_{j=1}^{n}\partial_i^{*}\,(     \mathcal D_{ij}+\mathcal Q_{ij})\,\partial_j.
\end{equation}
Then, we can decompose the generator $\mathcal{L}_\tau$ of the dynamics $(\ref{APP})$ into its symmetric and anti-symmetric parts as follows:
\begin{equation}\label{GG}
\mathcal{L}_\tau=\mathcal{L}_S+\mathcal{L}_A,
\quad
\mathcal{L}_S=\frac{\mathcal{L}_\tau+\mathcal{L}_\tau^*}{2},
\quad
\mathcal{L}_A=\frac{\mathcal{L}_\tau-\mathcal{L}_\tau^*}{2}.
\end{equation}
With $(\operatorname{div}\mathcal D)_j:=\sum_{i=1}^{n}\partial_i \mathcal D_{ij}$ and
$(\operatorname{div}Q)_j:=\sum_{i=1}^{n}\partial_i Q_{ij}$, for any smooth $f$,
\begin{align}
\mathcal L_S f
:=\nabla^{*}\mathcal D\nabla f
&=\sum_{i=1}^{n}\sum_{j=1}^{n}\partial_i^{*}\!\big(\mathcal D_{ij}\partial_j f\big)
 = \underbrace{\sum_{i=1}^{n}\sum_{j=1}^{n} \mathcal D_{ij}\,\partial_{ij}f}_{\text{second order}}
   + \underbrace{\big(\operatorname{div}\mathcal D-\mathcal D\nabla U\big)\!\cdot\nabla f}_{\text{first order}},
\label{eq:symmetric-part}\\[2mm]
\mathcal L_A f
:=\nabla^{*}\mathcal Q\nabla f
&=\sum_{i=1}^{n}\sum_{j=1}^{n}\partial_i^{*}\!\big(\mathcal Q_{ij}\partial_j f\big)
 = \big( \mathcal Q \nabla U-\operatorname{div}  \mathcal Q\big)\!\cdot\nabla f ,
\label{eq:antisymm-part}
\end{align}
where the second–order term vanishes in $\mathcal L_A$ thanks to $\mathcal Q^\top=-\mathcal Q$
(and the symmetry of the Hessian), i.e. $\sum_{i=1}^{n}\sum_{j=1}^{n} \mathcal Q_{ij}\partial_{ij}f=0$.
\begin{remark}[Interpretation of $\nabla^{*}  
\mathcal Q\nabla$]
The notation $\nabla^{*}\mathcal Q\nabla$ is a convenient quadratic-form notation;
for an antisymmetric $\mathcal Q$, it represents the \emph{pure first–order} operator
\eqref{eq:antisymm-part}.
\end{remark}

The rate function $I(\nu)$ (often referred to as the Donsker-Varadhan functional in the literature \cite{DV1,DV2,DV3,DV4}) takes the form
\begin{equation}\label{TR}
\ I(\nu)=\sup\left\{-\int_{\mathcal{X}}\frac{\mathcal{L}u}{u}d\nu,\ u\in\text{dom}(\mathcal{L})^+\right\},
\qquad\text{for any $\nu\in\mathcal{P}({\mathcal{X}})$},
\end{equation}
with
$$
\text{dom}(\mathcal{L})^+=\left\{u\in\mathcal{D}(\mathcal{L})\cap\mathcal{C}^0({\mathcal{X}})\ \Big|\ u>0,\ -\frac{\mathcal{L}u}{u}\in B_\kappa^\infty({\mathcal{X}})\right\},
$$
where $B_\kappa^\infty({\mathcal{X}})$ is defined in Appendix~\ref{Apendix definition}.

According to the symmetric and anti-symmetric parts of the generator, we can use the decomposition of the operator to obtain the symmetric and anti-symmetric
parts of the rate function respectively, as described in the following technical lemma.

\begin{lemma}[Theorem~3.3 in \cite{LDP-GG}]\label{RF}
Suppose that Assumptions~\ref{Hypoellipticity}-\ref{Lyapunov condition} and \ref{ass:kappa} hold. Consider a measure $\nu\in\mathcal{P}_\kappa({\mathcal{X}})$ such that $d\nu=e^\upsilon\ d\mu$ with $\upsilon\in\mathscr{H}^1(\nu)$ and  $\mathcal{L}_A\upsilon\in\mathscr{H}^{-1}(\nu)$. Then, the rate function $I$ defined in (\ref{TR}) reads:
$$I(\nu)=I_S(\nu)+I_A(\nu),$$
where
$$I_S(\nu)=\frac{1}{4}|\upsilon|_{\mathscr{H}^{1}(\nu)}^2,$$
and
$$I_A(\nu)=\frac{1}{4}\big|\mathcal{L}_A(\upsilon)\big|_{\mathscr{H}^{-1}(\nu)}^2.$$
\end{lemma}

Since $d\nu=e^{\upsilon}d\mu$ is the Radon–Nikodym derivative of $\nu$ with respect to
$\mu$, we can equivalently write
\[
I(\nu)
= \frac14\Bigl\lvert\log\frac{d\nu}{d\mu}\Bigr\rvert_{\mathscr{H}^1(\nu)}^2
 +\frac14\Bigl\lvert
   \mathcal L_A\!\left(\log\frac{d\nu}{d\mu}\right)
  \Bigr\rvert_{\mathscr{H}^{-1}(\nu)}^2 .
\]

%%%%%%%%%%%%%%%%%%%%%%%%%%%
\subsection{Main Results}		
%%%%%%%%%%%%%%%%%%%%%%%%%%%
Now we introduce our main result, that is, 
the large deviation principle for the generalized Langevin dynamics \eqref{APP} under Hypoellipticity (Assumption~\ref{Hypoellipticity}), Controllability (Assumption~\ref{irreducibility}) and Lyapunov condition (Assumption~\ref{Lyapunov condition}), 
and we obtain an explicit characterization of its rate function.

\begin{theorem}\label{thm:main}  
Suppose that Assumptions ~\ref{Hypoellipticity}-\ref{Lyapunov condition} and \ref{ass:kappa} hold, and let the function $\kappa$ be as in Lemma~\ref{Lyapunovkappa} and $\nu$ be any probability
measure in $\mathcal P_\kappa(\mathcal{X})$ of the form $d\nu=e^\upsilon\ d\mu$. Then the empirical measure $\{\pi_t\}_{t\ge 0}$ of the generalized Langevin dynamics defined by (\ref{APP}) satisfies a large deviation principle
in the $\tau^\kappa$-topology with the rate function given by
\begin{equation}\label{LDP rate function}
I_\tau(\nu)=\frac{1}{4}\int_{{\mathcal{X}}}\nabla\upsilon\cdot \mathcal{D}\nabla\upsilon\ d\nu\ +\frac{1}{4}\int_{{\mathcal{X}}}\nabla\psi_\upsilon\cdot \mathcal{D}\nabla\psi_\upsilon\ d\nu,
\end{equation}
where $\psi_\upsilon$ is the unique solution in $\mathscr{H}^{1}(\nu)$ to the Poisson equation
\begin{equation}\label{eqn:tilde:nabla}
\tilde{\nabla}\left(\mathcal{D}\nabla\psi_\upsilon\right)=\mathcal{L}_A\upsilon,
\end{equation}
where $\mathcal{L}_A=\sum_{i=1}^{n}\sum_{j=1}^{n}\left[\frac{\partial}{\partial\mathbf{z}_j}\mathcal{Q}_{ij}(\mathbf{z})-\mathcal{Q}_{ij}(\mathbf{z})\frac{\partial}
{\partial\mathbf{z}_j}H(\mathbf{z})\right]
\frac{\partial}{\partial\mathbf{z}_i}$. Moreover, $I_\tau(\nu)=+\infty$ for $\nu\notin\mathcal P_\kappa(\mathcal{X})$. That is to say for $\tau^\kappa$-measurable set $\varTheta\subset\mathcal P_\kappa(X)$, it holds
$$
-\mathop{\inf}\limits_{\nu\in\mathring{\varTheta}} I_{\tau}(\nu)\le\mathop{\liminf}\limits_{t\rightarrow+\infty}\frac{1}{t}\log \mathbb{P}(\pi_t\in\varTheta)\le\mathop{\limsup}\limits_{t\rightarrow+\infty}\frac{1}{t}\log \mathbb{P}(\pi_t\in\varTheta)\le -\mathop{\inf}\limits_{\nu\in \bar{\varTheta}} I_{\tau}(\nu).
$$ 
\end{theorem}%\proofat{app:pf-thm:main}

\begin{proof}
The proof will be provided in Appendix~\ref{app:pf-thm:main}.
\end{proof}

\begin{remark}
Since $\kappa$ is bounded on $\mathcal{X}$,  $\mathcal P_\kappa(\mathcal{X})=\mathcal P(\mathcal{X})$ and the statement of
Theorem~\ref{thm:main} (as well as Lemma~2.6) automatically holds on $\mathcal P(\mathcal{X})$.
\end{remark}

Next, we are interested in applying the general framework of large deviations for
the generalized Langevin dynamics \eqref{APP} to study the large deviations
for the variants of Langevin dynamics.
For each variant of Langevin dynamics, 
we carefully verify Hypoellipticity (Assumption~\ref{Hypoellipticity}), Controllability (Assumption~\ref{irreducibility}) and Lyapunov condition (Assumption~\ref{Lyapunov condition}).
According to the form of generalized Langevin dynamics \eqref{APP}, we can
divide the variants of Langevin dynamics of interest 
into two cases to study:
\begin{itemize}
\item $\mathbf{z}$ represents the model parameter $\theta$ itself, that is, $g(\theta,r)=0$ and $H(\mathbf{z})=U(\theta)$. This case includes the {\it overdamped Langevin dynamics} \eqref{eq:overdamped-2}, the {\it non-reversible Langevin dynamics} \eqref{non:reversible},  and the {\it mirror Langevin dynamics} \eqref{mirror Langevin dynamics}.
\end{itemize}

\begin{itemize}
\item $\mathbf{z}$ represents the extended state space, that is, the case of $g(\theta,r)\neq 0$.
This case includes the {\it underdamped Langevin dynamics} \eqref{eqn:underdamped}, the {\it high-order Langevin dynamics} \eqref{JI} and the {\it Hessian-free high-resolution dynamics} \eqref{Hessian-free high-resolution}.
\end{itemize}

%For the model framework in this case, we don't get the general large deviation principle. But we apply it to two examples,
%and construct the Langevin dynamics process in this case by selecting the appropriate positive semidefinite diffusion matrix $\mathcal{D}$
%and anti-symmetric curl matrix $\mathcal{Q}$.

Before we proceed, let us recall that the LDPs for overdamped Langevin dynamics \eqref{eq:overdamped-2}, underdamped Langevin dynamics \eqref{eqn:underdamped} and non-reversible Langevin dynamics \eqref{non:reversible} have already been obtained
in the literature; see e.g. \cite{LDP-GG}, restated as Lemma~\ref{LDP-A}, Lemma~\ref{lem:underdamped} and Lemma~\ref{lemma:non:reversible} in Appendix~\ref{sec:technical:lemmas}.

%%%%%%%%%%%%%%%%%%%%%%%%%%%%%%%%%%%%%%%%%%%%%%%%%%%%%%%%%

We will derive the LDPs for mirror Langevin dynamics \eqref{mirror Langevin dynamics}, 
high-order Langevin dynamics \eqref{JI} and Hessian-free high-resolution dynamics \eqref{Hessian-free high-resolution}, which to the best of our knowledge have never been studied in the previous literature.

%%%%%%%%%%%%%%%%%%%%%%%%%%%%%%%%%%%%%%%%%%%%%%%%%%%%%%%%%
\subsubsection{Mirror Langevin dynamics}\label{sec:mirror}

Let $\mu(\mathrm d\theta)\propto e^{-U(\theta)}\mathrm d\theta$ and work on $L^{2}(\mu)$, so that
$\nabla^{*}=-\mathrm{div}+\nabla U\cdot$. Let $\phi:\mathcal X\to\mathbb R$ be a twice–differentiable strictly convex function. With
\[
\mathcal{Q}\equiv 0,\qquad
\mathcal{D}(\theta)=\big[\nabla^{2}\phi(\theta)\big]^{-1}\succeq 0,
\]
the infinitesimal generator of the mirror Langevin dynamics can be written in the canonical form
\begin{equation}\label{eq:ML-generator}
\mathcal L_M=\nabla^{*}\mathcal D\,\nabla,
\qquad\text{i.e.}\qquad
\mathcal L_M f
=\operatorname{tr}\!\big(\mathcal D\nabla^{2} f\big)
+\big(\mathrm{div}\,\mathcal D-\mathcal D\nabla U\big)\!\cdot\nabla f .
\end{equation}
Then its symmetric and anti-symmetric parts are
\begin{align*}
\mathcal L_{MS}&=\nabla^{*}\big[\nabla^{2}\phi(\theta)\big]^{-1}\,\nabla\\
&=\operatorname{tr}\!\big(\big[\nabla^{2}\phi(\theta)\big]^{-1}(\theta)\,\nabla^{2} f(\theta)\big)
+ \big(\operatorname{div}\big[\nabla^{2}\phi(\theta)\big]^{-1}(\theta)-\big[\nabla^{2}\phi(\theta)\big]^{-1}(\theta)\nabla U(\theta)\big)\!\cdot\nabla f(\theta),\\
\mathcal L_{M,\mathrm A} &=\nabla^{*}\mathcal Q\,\nabla \equiv 0
\qquad (\text{since } \mathcal Q\equiv 0 \text{ in mirror Langevin dynamics}).
\end{align*}
\begin{remark}[On antisymmetric part]
In the canonical decomposition $\mathcal L=\nabla^{*}(\mathcal D+\mathcal Q)\nabla$,
the antisymmetric contribution is first order:
$\nabla^{*}\mathcal Q\nabla=(\mathcal Q\nabla U-\mathrm{div}\,\mathcal Q)\cdot\nabla$.
For mirror Langevin dynamics, $\mathcal Q\equiv0$, and thus only the symmetric part
$\nabla^{*}\mathcal D\nabla$ remains; the lengthy coordinate computations in
Appendix~\ref{app:pf-thm:main} become unnecessary.
\end{remark}
Let $W_M(\theta)=e^{\eta U(\theta)}$ with $\eta\in(0,1)$. Using
$\nabla W_M=\eta W_M\nabla U$ and
$\nabla^{2}W_M=\eta W_M\nabla^{2}U+\eta^{2}W_M\,\nabla U\\\otimes\nabla U$,
substituting $f=W_M$ into \eqref{eq:ML-generator} yields the Witten–Lyapunov drift
\begin{equation}\label{eq:PsiM}
\Psi_M(\theta):=-\frac{\mathcal L_M W_M}{W_M}
=\eta\Big((1-\eta)\,\nabla U\!\cdot\!\big[\nabla^2\phi(\theta)\big]^{-1}\nabla U
- \mathrm{tr}\!\Big(\big[\nabla^2\phi(\theta)\big]^{-1}\nabla^2U\Big)
+ \nabla\!\cdot\!\Big(\big[\nabla^2\phi(\theta)\big]^{-1}\Big)\!\cdot\nabla U\Big),
\end{equation}
which is the quantity displayed in \eqref{mirrorLyapinov2}. Therefore condition \eqref{Lyapunov2} is exactly the
statement that \(\Psi_M(\theta)\to+\infty\) as \(|\theta|\to\infty\).
Consequently, \(W_M\) has compact level sets and satisfies the Lyapunov
drift condition (Assumption~\ref{Lyapunov condition}) for the infinitesimal generator \(\mathcal L_M\). 

\begin{theorem}\label{Mirror} 
Assume that the potential $U\in \mathscr{S}$ has compact level sets, $e^{-U}\in L^1(\mathcal{X})$, and for any $\eta\in(0, 1)$, it holds that
\begin{equation}\label{Lyapunov2}
(1-\eta)|\sqrt{[\nabla^2\phi(\theta)]^{-1}}\cdot\nabla U|^2-[\nabla^2\phi(\theta)]^{-1}:\nabla^2U-\nabla\cdot [\nabla^2\phi(\theta)]^{-1}\nabla U \xrightarrow[|\theta|\rightarrow +\infty]\quad+\infty.
\end{equation}
Then the mirror Langevin dynamics \eqref{mirror Langevin dynamics} admits the function
\begin{align*}
W_M(\theta) = e^{\eta U(\theta)}
\end{align*}
for any $\eta\in(0, 1)$ as a Lyapunov function in the sense of Assumption~\ref{Lyapunov condition}. Moreover,
\begin{align}\label{mirrorLyapinov2}
\Psi_M
&:=-\frac{\mathcal{L}_M W_M}{W_M}
\nonumber
\\
&=\eta\left((1-\eta)|\sqrt{[\nabla^2\phi(\theta)]^{-1}}\nabla U|^2-[\nabla^2\phi(\theta)]^{-1}:\nabla^2 U-\nabla \cdot[\nabla^2\phi(\theta)]^{-1} \cdot \nabla U \right)
\end{align}
has compact level sets and, for any $\kappa:\mathcal{X} \rightarrow [1,\infty)$ belonging to $\mathscr{S}$, bounded or with
compact level sets and 
$$\frac{\Psi_M(\mathbf{z})}{\kappa(\mathbf{z})}\xrightarrow[|\mathbf{z}|\rightarrow +\infty]\quad +\infty,$$
the empirical measure $\{\pi_t\}_{t\ge0}$  satisfies a large deviation principle
in the $\tau^\kappa$-topology and the corresponding rate function is defined by
\begin{equation}\label{mirror rate function}
I_M(\nu)=\frac{1}{4}\int_{\mathcal{X}}\nabla\upsilon\cdot [\nabla^2\phi(\theta)]^{-1}\nabla\upsilon\ d\nu,
\end{equation}
where $d\nu=e^\upsilon d\mu$. That is to say for $\tau^\kappa$-measurable set $\varTheta\subset\mathcal{P}(\mathcal{X})$, it holds
$$
-\mathop{\inf}\limits_{\nu\in\mathring{\varTheta}} I_{M}(\nu)\le\mathop{\liminf}\limits_{t\rightarrow+\infty}\frac{1}{t}\log \mathbb{P}(\pi_t\in\varTheta)\le\mathop{\limsup}\limits_{t\rightarrow+\infty}\frac{1}{t}\log \mathbb{P}(\pi_t\in\varTheta)\le -\mathop{\inf}\limits_{\nu\in \bar{\varTheta}} I_{M}(\nu).
$$
\end{theorem}%\proofat{app:pf-Mirror}

\begin{proof}
The proof will be provided in Appendix~\ref{app:pf-Mirror}.
\end{proof}

%%%%%%%%%%%%%%%%%%%%%%%%%%%%%%%%%%%
\subsubsection{High-order Langevin dynamics}\label{High-order}

For the high-order Langevin dynamics \eqref{JI}, according to the values of $\mathcal{D}$ and $\mathcal{Q}$ set in (\ref{high-order cofficient}), 
through formula (\ref{BH}), we can compute that the infinitesimal generator of \eqref{JI} is given by
\begin{equation}\label{High-order generator}
\mathcal{L}_{H}=\mathcal{L}_{\mathrm{HS}} + \mathcal{L}_{\mathrm{HA}}=\alpha\Delta_r- (\gamma p+\alpha r)\cdot\nabla_r  + p\cdot\nabla_\theta -( \nabla U - \gamma r)\cdot\nabla_p,
\end{equation}
where
\begin{align*}
&\mathcal{L}_{\mathrm{HS}}:=-\alpha r\cdot\nabla_r + \alpha\Delta_r,\\
&\mathcal{L}_{\mathrm{HA}}:=p\cdot\nabla_\theta -\nabla U\cdot\nabla_p +\gamma  r\cdot\nabla_p - \gamma  p\cdot\nabla_r.
\end{align*}

In order to obtain the large deviation principle for the high-order Langevin dynamics (\ref{JI}), we make the following classical assumptions for the growth of potential function $U$.

\begin{assumption}\label{HHJ}
The potential $U\in\mathscr{S}$ has compact level sets, satisfies $e^{-U}\in L^1(\mathcal{X})$, there exist $k > 1$ and $M_U, c_U, m_U > 0$ such that for all $\theta \in \mathbb{R}^d$ with $|\theta| \geq c_U$,
\begin{itemize}
\item[(a)] $m_U |\theta|^k \leq U(\theta) \leq M_U |\theta|^k$ \quad \text{and} \quad $m_U |\theta|^k \leq \theta \cdot \nabla U(\theta)$.

\item[(b)]$|\nabla U(\theta)| \leq M_U |\theta|^{k-1}$.
\end{itemize}
\end{assumption}

In the next lemma, we show that Hypoellipticity (Assumption~\ref{Hypoellipticity}) and Controllability (Assumption~\ref{irreducibility}) are satisfied
for the high-order Langevin dynamics \eqref{JI}.

\begin{lemma}\label{Hypoellipticity-Controllability} (Hypoellipticity and Controllability). The generator $\mathcal{L}_{H}$ defined in (\ref{High-order generator})
satisfies the Hypoelliptic H\" ormander and Controllability conditions.
\end{lemma}%\proofat{app:pf-Hypoellipticity-Controllability}

\begin{proof}
The proof will be provided in Appendix~\ref{app:pf-Hypoellipticity-Controllability}.
\end{proof}

\newcommand{\GammaH}[1]{\alpha\,|\nabla_r #1|^{2}}
In the next lemma, we show that Lyapunov condition (Assumption~\ref{Lyapunov condition})
is satisfied
for the high-order Langevin dynamics \eqref{JI}.

\begin{proposition}[Witten--Lyapunov $\Rightarrow$ Lyapunov drift for high-order Langevin]\label{Witten-Lyapunov condition-high:order}.
Suppose that $(\bm{\mathbf{z}}_t)_{t\ge0}=(\theta_t,p_t,r_t)_{t\ge0}$ is as in \eqref{JI} and $U$ satisfies Assumption~\ref{HHJ}. Let $a>0$ be sufficiently small. Then for any $\alpha, \gamma>0$ and $k\in(1,2]$,
there exists $\delta\in (\frac{2-k}{k},1]$ such that
\begin{equation}\label{ABA}
W_\delta(\theta,p,r)=e^{\varphi_{\mathrm{HL}}^\delta(\theta,p,r)}
\end{equation}
is a Lyapunov function, where 
$\varphi_{\mathrm{HL}}(\theta,p,r)=\varphi_0(\theta,p,r) - \inf_{(\tilde{\theta},\tilde{p},\tilde{r})\in\mathbb{R}^{3d}} \varphi_0(\tilde{\theta},\tilde{p},\tilde{r}) + 1$ 
and $\varphi_0(\theta,p,r)\\:
=hH(\theta,p,r)+aL(\theta)\cdot p+ap\cdot r$ for some $h>0$ with $H(\theta,p,r)$ given in \eqref{H:function:high:order} and $L(\theta)$ in \eqref{1}.
\medskip

\noindent\emph{(WL inequality).}
There exist $c_0>0$, $M<\infty$ and a compact set $K\Subset\mathcal X$ such that
\begin{equation}\label{eq:WL-high}
\mathcal L_H\big(\log W_\delta\big)\;+\;\GammaH{\log W_\delta}
\ \le\
-\,c_0\,\GammaH{\log W_\delta}\;+\; M\,\mathbf 1_K .
\end{equation}

\noindent\emph{(Consequently: Lyapunov drift = Assumption \ref{Lyapunov condition}).}
There exists $\theta\in(0,1)$ and constants $c_1,c_2>0$, $R,R'<\infty$ such that
\begin{equation}\label{eq:Lyap-high}
c_1\,\alpha\,|\nabla_r W_\delta|^{2}-R'
\ \le\
-\mathcal L_H W_\delta - \theta\,\alpha\,|\nabla_r W_\delta|^{2}
\ \le\
c_2\,\alpha\,|\nabla_r W_\delta|^{2}+R .
\end{equation}
Equivalently,
\[
-\mathcal L_H W_\delta - \theta\,|\sqrt{\mathcal D}\nabla W_\delta|^{2}
\ \sim\
|\sqrt{\mathcal D}\nabla W_\delta|^{2}\qquad(\mathcal D=\mathrm{diag}(0,0,\alpha I_d)),
\]
so Assumption~\ref{Lyapunov condition} holds with $W=W_\delta$.
In particular, there exist $h,A,B,C,D>0$ such that
\begin{equation}\label{LPC}
-\frac{\mathcal{L}_{H}W_\delta}{W_\delta}\ge A|\theta|^{2(k-1)} + B|p|^2 + C|r|^2 - D.
\end{equation}
\end{proposition}%\proofat{app:pf-Witten-Lyapunov condition-high:order}

\begin{proof}
The proof will be provided in Appendix~\ref{app:pf-Witten-Lyapunov condition-high:order}.
\end{proof}

Now, we are ready to state the following large deviations result 
for the high-order Langevin dynamics \eqref{JI}.

\begin{theorem}\label{thm:high:order}
Assume that $(\bm{\mathbf{z}}_t)_{t\ge0}=(\theta_t,p_t,r_t)_{t\ge0}$ is as 
in \eqref{JI} and $U$ satisfies Assumption~\ref{HHJ}, and consider a smooth function $\kappa$
with $\kappa(\theta,p,r)=1 + |\theta|^\lambda + |p|^\sigma + |r|^\omega$ with $\lambda\in[0,2), \sigma\in[0,2), \omega\in[0,2)$. Then $(\bm{\mathbf{z}}_t)_{t\ge0}$
is ergodic with respect to the measure $\mu$, with Lyapunov function defined in (\ref{ABA}). Moreover, the empirical measure
$$
\pi_{t}:= \frac{1}{t}\int_{0}^{t}\delta_{(\theta_s,p_s,r_s)}ds
$$
satisfies an LDP in the $\tau^\kappa$-topology. Finally, for any $\nu\in\mathcal{P}_\kappa(\mathcal{X})$ such that $d\nu=e^\upsilon d\mu$, the rate function reads
\begin{equation}\label{eqn: rate function High order}
I_{H}(\nu)=\frac{\alpha}{4}\int_{\mathcal{X}}|\nabla_r\upsilon|^2\ d\nu + \frac{1}{4\alpha}\int_{\mathcal{X}}|\nabla_r\psi|^2\ d\nu,
\end{equation}
where $\psi$ is the unique solution in $\mathscr{H}^{1}(\nu)$ to the Poisson equation
\begin{equation}\label{Poisson-equation-High-order}
- \Delta_r \psi + (r - \nabla_r \upsilon) \cdot \nabla_r \psi=\mathcal{L}_{\mathrm{HA}}\upsilon.
\end{equation}
\end{theorem}%\proofat{app:pf-thm:high:order}

\begin{proof}
The proof will be provided in Appendix~\ref{app:pf-thm:high:order}.
\end{proof}

%%%%%%%%%%%%%%%%%%%%%%%%%%%%%%%%%%%%

%%%%%%%%%%%%%%%%%%%%%%%%%%%%%%%%%%%%%%%%
\subsubsection{Hessian-free high-resolution dynamics}

According to the choices of $\mathcal{D}$ and $\mathcal{Q}$ in (\ref{Hessian-free high-resolution}),
through formula (\ref{BH}), we can compute that the infinitesimal generator of \eqref{Hessian-free high-resolution} is given by
\begin{equation*}
\begin{aligned}
&\mathcal{L}_{\mathrm{RS}}=-\beta\nabla U(\theta)\cdot\nabla_\theta - \alpha r\cdot\nabla_r + \beta\Delta_\theta + \alpha\Delta_r,\\
&\mathcal{L}_{\mathrm{RA}}=r\cdot\nabla_\theta - \nabla U(\theta)\cdot\nabla_r,
\end{aligned}
\end{equation*}
which implies that
\begin{equation}\label{BCB}
\mathcal{L}_R=\mathcal{L}_{\mathrm{RS}} + \mathcal{L}_{\mathrm{RA}}=(-\beta\nabla U(\theta) + r)\cdot\nabla_\theta +(- \alpha r- \nabla U(\theta))\cdot\nabla_r +\beta\Delta_\theta + \alpha\Delta_r.
\end{equation}

In order to obtain the large deviation principle for the Hessian-free high-resolution dynamics (\ref{Hessian-free high-resolution}), we make the following classical assumptions for the growth of potential function $U$.
\begin{assumption}\label{HFHR}
The potential $U\in\mathscr{S}$ has compact level sets, satisfies $e^{-U}\in L^1(\mathcal{X})$ and  
\begin{itemize}
  \item[(a)] There exist $c_1>0$, $C_1\in\mathbb{R}$ such that $\theta\cdot\nabla U(\theta)\ge c_1 |\theta|^2-C_1$ for all $\theta\in\mathbb{R}^{d}$.
  \item[(b)] There exists $m_{U}>0$ such that $|\nabla U(\theta_1) - \nabla U(\theta_2)| \leq m_U |\theta_1-\theta_2|$ for all $\theta_{1},\theta_{2}\in\mathbb{R}^{d}$.
\end{itemize}
\end{assumption}

The Hessian-free high-resolution dynamics (\ref{Hessian-free high-resolution}) automatically satisfies Hypoellipticity (Assumptions~\ref{Hypoellipticity}) and Controllability (Assumption~\ref{irreducibility}) since its infinitesimal generator is an ellpitic operator. In the next lemma, we will show that the Lyapunov condition (Assumption~\ref{Lyapunov condition})
is satisfied
for the Hessian-free high-resolution dynamics (\ref{Hessian-free high-resolution}).

\begin{lemma}\label{Witten-Lyapunov condition-Hessian-Free}(Witten-Lyapunov condition).
Suppose that $(\bm{\mathbf{z}}_t)_{t\ge0}=(\theta_t,r_t)_{t\ge0}$ is as in \eqref{Hessian-free high-resolution} and assume that $U$ satisfies Assumption~\ref{HFHR}. Then for any $\alpha,\beta>0$ and $a\in(0,1)$,
\begin{equation}\label{ABA1}
W_a(\theta,r)=e^{aH(\theta,r)}
\end{equation}
is a Lyapunov function, where $H(\theta,r)=U(\theta)+\frac12 |r|^2$ given in \eqref{H:function:HFHR}. More specifically, there exist
$A,B,C>0$ such that
\begin{equation}\label{LPC1}
-\frac{\mathcal{L}_RW_a}{W_a}\ge A|\theta|^2 + B|r|^2 - C.
\end{equation}
\end{lemma}%\proofat{app:pf-Witten-Lyapunov condition-Hessian-Free}

\begin{proof}
The proof will be provided in Appendix~\ref{app:pf-Witten-Lyapunov condition-Hessian-Free}.
\end{proof}

Now, we are ready to state the following large deviations result 
for the Hessian-free high-resolution dynamics (\ref{Hessian-free high-resolution}).

\begin{theorem}\label{thm: Hessian-Free}
Assume that $(\bm{\mathbf{z}}_t)_{t\ge0}=(\theta_t,r_t)_{t\ge0}$ is as in \eqref{Hessian-free high-resolution} and $U$ satisfies Assumption ~\ref{HFHR}, and consider a smooth function $\kappa$
with $\kappa(\theta,r)=1 + |\theta|^\lambda+ |r|^\omega$ and $\lambda\in[0,2), \omega\in[0,2)$. Then $(\bm{\mathbf{z}}_t)_{t\ge0}$
is ergodic with respect to the measure $\mu$, with Lyapunov function defined in (\ref{ABA1}). Moreover, the empirical measure
$$
\pi_{t}:= \frac{1}{t}\int_{0}^{t}\delta_{(\theta_s,r_s)}ds
$$
satisfies a LDP in the $\tau^\kappa$-topology. Finally, for any $\nu\in\mathcal{P}_\kappa(\mathcal{X})$ such that $d\nu=e^{\upsilon} d\mu$, the rate function reads
\begin{equation}\label{rate function Hessian Free}
I_R(\nu)=\frac{\beta}{4}\left(\int_{\mathcal{X}}|\nabla_\theta\upsilon|^2\ d\nu + \int_{\mathcal{X}}|\nabla_\theta\psi_\upsilon|^2\ d\nu\right)
+ \frac{\alpha}{4}\left(\int_{\mathcal{X}}|\nabla_r\upsilon|^2\ d\nu + \int_{\mathcal{X}}|\nabla_r\psi_\upsilon|^2\ d\nu\right),
\end{equation}
where $\psi_\upsilon$ is the unique solution in $\mathscr{H}^{1}(\nu)$ to the Poisson equation
\begin{equation}\label{Poisson-equation-Hessian-Free}
- \beta \, \Delta_{\boldsymbol{\theta}} \psi_{\upsilon}
- \alpha \, \Delta_{{r}} \psi_{\upsilon}
+ \beta \, (\nabla_{\boldsymbol{\theta}} U - \nabla_{\boldsymbol{\theta}} \upsilon) \cdot \nabla_{\boldsymbol{\theta}} \psi_{\upsilon}
+ \alpha \, ({r} - \nabla_{{r}} \upsilon) \cdot \nabla_{{r}} \psi_{\upsilon}= \mathcal{L}_{\mathrm{RA}} \upsilon.
\end{equation}
\end{theorem}%\proofat{app:pf-thm: Hessian-Free}

\begin{proof}
The proof will be provided in Appendix~\ref{app:pf-thm: Hessian-Free}.
\end{proof}

%%%%%%%%%%%%%%%%%%%%%%%%%%%%%%%%
\section{Comparisons}\label{sec:comparison}

In this section, by using the large deviation rate function of empirical measure as a tool to analyze the convergence rate to the invariant measure, and taking the large deviation rate function as a measure of the rate of convergence to equilibrium, we study whether variants of Langevin dynamics have acceleration effect over the overdamped Langevin dynamics. The intuition is that
the larger the rate function, the more concentrated the empirical measure is around the invariant measure, thus indicating faster convergence. 
First, we consider the acceleration of the generalized Langevin dynamics \eqref{APP}
We have the following result.

\begin{corollary}\label{thm:accelerate}
Suppose that Assumptions~\ref{Hypoellipticity}-\ref{Lyapunov condition} hold. 
Let $\nu$ be any measure in $\mathcal{P}_k (\mathcal{X}=\mathbb{R}^d)$ of the form
$d\nu=e^\upsilon\ d\mu$ and $\nu\neq\mu$. If the matrix $\mathcal{D}-\mathbf{I}$ is a positive semidefinite diffusion matrix, then we have
\begin{equation}
I_\tau(\nu)\ge I_o(\nu),
\end{equation}
where $I_\tau(\nu)$ is the rate function of the generalized Langevin dynamics \eqref{APP} and $I_o(\nu)$ is the rate function \eqref{eq:Io-variational} of the overdamped Langevin dynamics (\ref{eq:overdamped-2}). 
\end{corollary}%\proofat{app:pf-thm:accelerate}

\begin{proof}
The proof will be provided in Appendix~\ref{app:pf-thm:accelerate}.
\end{proof}

Corollary~\ref{thm:accelerate} shows the acceleration effect of the generalized Langevin dynamics compared to the overdamped Langevin dynamics.
%%%%%%%%%%%%%%%%%%%%%%%%%%%%%%%%%%%%%%%%
In the rest of this section, 
we show the acceleration of variants of Langevin dynamics
over the classical overdamped Langevin dynamics by comparing
their rate functions from the large deviations.
It is known that the rate function for non-reversible Langevin dynamics \eqref{non:reversible}
is greater than that for overdamped Langevin dynamics \eqref{eq:overdamped-2}, 
and hence the acceleration \cite{LDP-GG}.
In the rest of this section, we will compare the rate functions for mirror Langevin dynamics \eqref{mirror Langevin dynamics}, Hessian-free high-resolution dynamics \eqref{Hessian-free high-resolution}, underdamped Langevin dynamics \eqref{eqn:underdamped}, and high-order Langevin dynamics \eqref{JI}
with that of overdamped Langevin dynamics \eqref{eq:overdamped-2}.

Note that for underdamped Langevin dynamics, high-order Langevin dynamics,
and Hessian-free high-resolution dynamics, these processes live in higher dimensions
than overdamped Langevin dynamics, which makes a direct comparison
of their rate functions impossible.
%%%%%%%%%%%%%%%%%%%%%%%%%%%%%%%%%%%%
We introduce a novel method to expand the space of overdamped Langevin dynamics to match the dimension with that of the variant of Langevin dynamics under study. In particular
expanding to the space of overdamped Langevin dynamics to $\mathbb{R}^{2d}$ when we compare it with underdamped Langevin dynamics and Hessian-free high-resolution dynamics and to $\mathbb{R}^{3d}$ when we compare it with high-order Langevin dynamics.
If the comparison can be made on the whole expanded space, 
we then apply the contraction principle from the large deviations theory to
compare the rate functions for the large deviations of $\frac{1}{t}\int_{0}^{t}\delta_{\theta_{s}}ds$. 

First, we compare the rate function of the mirror Langevin dynamics \eqref{mirror Langevin dynamics}
with that of the overdamped Langevin dynamics \eqref{eq:overdamped-2}.
The following result can be obtained as a corollary of Corollary~\ref{thm:accelerate}.

\begin{corollary}\label{Mirror VS overdamped} 
Assume the same condition as Theorem~\ref{Mirror}. Let $\nu$ be any measure in $\mathcal{P}_k (\mathcal{X}=\mathbb{R}^d)$ of the form $d\nu=e^\upsilon\ d\mu$ and $\nu\neq\mu$. Then for the mirror Langevin dynamics defined by (\ref{mirror Langevin dynamics}), if the matrix $[\nabla^2\phi(\mathbf{z})]^{-1}-I$ is a positive semidefinite diffusion matrix, then we have
\begin{equation}
I_M(\nu)\ge I_o(\nu),
\end{equation}
where $I_o(\nu)$ is the rate function of the overdamped Langevin dynamics (\ref{eq:overdamped-2}) and $I_M(\nu)$ is the rate function of the mirror Langevin dynamics (\ref{mirror rate function}).
\end{corollary}
%\proofat{app:pf-Mirror VS overdamped}

\begin{proof}
The proof will be provided in Appendix~\ref{app:pf-Mirror VS overdamped}.
\end{proof}

Corollary~\ref{Mirror VS overdamped} shows the acceleration effect of the mirror Langevin dynamics \eqref{mirror Langevin dynamics} compared to the overdamped Langevin dynamics \eqref{eq:overdamped-2}.

Next, we compare the rate functions of Hessian-free high-resolution dynamics \eqref{Hessian-free high-resolution} and underdamped Langevin dynamics \eqref{eqn:underdamped} with overdamped Langevin dynamics \eqref{eq:overdamped-2}. 
Note that overdamped Langevin dynamics lives in $\mathbb{R}^{d}$, 
whereas Hessian-free high-resolution dynamics and underdamped Langevin dynamics live in $\mathbb{R}^{2d}$.
Thus, we lift the overdamped Langevin dynamics to $\mathbb{R}^{2d}$, before we do the comparison.
We consider the \textit{expanded second-order overdamped Langevin dynamics}, for $(\theta_t, r_t)\in\mathbb{R}^d\times\mathbb{R}^d$,
\begin{equation}\label{eqn:expanded-2th-overdamped}
\begin{cases}
d\theta_t=-\nabla U(\theta_t)dt+\sqrt{2}dW_t,\\
dr_t=-r_t dt+\sqrt{2}dB_t,
\end{cases}
\end{equation}
where $W_t$ and $B_t$ are independent $d$-dimensional Brownian motions. We can easily check that its invariant measure is $d\mu=e^{-U(\theta)-\frac{1}{2}| r|^2}d\theta dr$, that is the same as the invariant measure of Hessian-free high-resolution dynamics. 

Lemma~\ref{LDP-A} is a direct application of Lemma~\ref{Lyapunovkappa}, Proposition~\ref{prop:Wasserstein} and Theorem \ref{thm:main}. By LDP for the overdamped Langevin dynamics (Lemma~\ref{LDP-A}), the LDP rate function for the expanded second-order overdamped Langevin dynamics \eqref{eqn:expanded-2th-overdamped} is 
\begin{equation}\label{rate function expanded 2th overdamped}
I_{e2o}(\nu)=\frac{1}{4}\int_{\mathcal{X}}|\nabla\upsilon|^2\ d\nu, 
\end{equation}
where \(\mathcal X=\mathbb{R}^{2d}\) and \(\nu\in\mathcal P_{\kappa}\mathcal X)\) satisfies \(\nu\ll\mu\) with density \(e^{\upsilon}\), i.e. \(d\nu=e^{\upsilon}\,d\mu\).

\begin{proposition}\label{Hessian-free high-resolution VS overdamped} 
Assume the same condition as Theorem~\ref{thm: Hessian-Free}. Let $\nu$ be any measure in $\mathcal{P}_{\kappa} (\mathcal{X}=\mathbb{R}^{2d})$ of the form $d\nu=e^\upsilon\ d\mu$. Then for the Hessian-free high-resolution dynamics defined by (\ref{Hessian-free high-resolution}), if $\min(\alpha,\beta)\geq 1$, then we have
\begin{equation}
I_R(\nu)\ge I_{e2o}(\nu),
\end{equation}
where $I_{e2o}(\nu)$ , given in \eqref{rate function expanded 2th overdamped}, is the rate function of the expanded second-order overdamped Langevin dynamics (\ref{eqn:expanded-2th-overdamped}) and $I_R(\nu)$ is the rate function of the Hessian-free high-resolution dynamics (\ref{Hessian-free high-resolution}). 
\end{proposition}

%\proofat{app:pf-Hessian-free high-resolution VS overdamped}

\begin{proof}
The proof will be provided in Appendix~\ref{app:pf-Hessian-free high-resolution VS overdamped}.
\end{proof}

Proposition~\ref{Hessian-free high-resolution VS overdamped} show the acceleration effect of the Hessian-free high-resolution dynamics \eqref{Hessian-free high-resolution} compared to the expanded second-order overdamped Langevin dynamics \eqref{eqn:expanded-2th-overdamped}. 
Note that, 
for the Hessian-free high-resolution dynamics \eqref{Hessian-free high-resolution}, the $\theta$-marginal distribution of its stationary distribution is the target distribution $\mu(\theta)\propto e^{-U(\theta)}$ which coincides with
the $\theta$-marginal distribution of the stationary distribution
of the expanded second-order overdamped Langevin dynamics \eqref{eqn:expanded-2th-overdamped}, which is indeed the stationary distribution
of the overdamped Langevin dynamics \eqref{eq:overdamped-2}.
By contraction principle (see e.g. Theorem~4.2.1. in \cite{DZ1998}), for the Hessian-free high-resolution dynamics \eqref{Hessian-free high-resolution}, $\frac{1}{t}\int_{0}^{t}\delta_{\theta_{s}}ds$ satisfies a LDP in the $\tau^{\kappa}$-topology with the rate function
\begin{equation}\label{I:R:theta}
I_{R,\theta}(\nu_{\theta}):=\inf_{\nu\in\mathcal{P}_{\kappa}(\mathbb{R}^{2d}):M_{\theta}(\nu)=\nu_{\theta}}I_{R}(\nu),
\end{equation}
where $M_{\theta}(\nu)$ is the $\theta$-marginal distribution of $\nu$.
Similarly, by contraction principle, the rate function \eqref{eq:Io-variational} for overdamped Langevin dynamics~\eqref{eq:overdamped-2} can be re-written as
\begin{equation}\label{I:o:contraction}
I_{o}(\nu_{\theta}):=\inf_{\nu\in\mathcal{P}_{\kappa}(\mathbb{R}^{2d}):M_{\theta}(\nu)=\nu_{\theta}}I_{e2o}(\nu).
\end{equation}
Therefore, we obtain the following corollary from Proposition~\ref{Hessian-free high-resolution VS overdamped}.

\begin{corollary}\label{cor:HFHR:vs:overdamped}
Under the same assumptions as in Proposition~\ref{Hessian-free high-resolution VS overdamped}, for any $\nu_{\theta}\in\mathcal{P}_{\kappa}(\mathbb{R}^{d})$, 
\begin{equation}
I_{R,\theta}(\nu_{\theta})\geq I_{o}(\nu_{\theta}),
\end{equation}
where $I_{R,\theta}$ and $I_{o}$ are defined in \eqref{I:R:theta} and \eqref{eq:Io-variational}.
\end{corollary}

%\proofat{app:pf-cor:HFHR:vs:overdamped}

\begin{proof}
The proof will be provided in Appendix~\ref{app:pf-cor:HFHR:vs:overdamped}.    
\end{proof}

Corollary~\ref{cor:HFHR:vs:overdamped} shows the acceleration
effect of the Hessian-free high-resolution dynamics \eqref{Hessian-free high-resolution} compared to the expanded second-order overdamped Langevin dynamics \eqref{eqn:expanded-2th-overdamped} in terms of the convergence of the $\theta$-marginal dynamics to the $\theta$-marginal of the target distribution
\[
d\mu = e^{-\,U(\theta)\;-\;\tfrac12\,|r|^{2}}\, d\theta\, dr .
\]

Next, we compare the rate functions of underdamped Langevin dynamics \eqref{eqn:underdamped} and high-order Langevin dynamics \eqref{JI} with overdamped Langevin dynamics \eqref{eq:overdamped-2}.
We will make comparisons on a subspace of the probability measures
that is defined as follows. It remains an open problem how to compare
the rate functions on the whole space of the probabilities
for underdamped Langevin dynamics \eqref{eqn:underdamped} and high-order Langevin dynamics \eqref{JI}, which will be left as a future research direction.

\begin{definition} 
The class $\mathcal{P}^H(\mathcal{X})$ consists of $\nu\in\mathcal{P}(\mathcal{X}=\mathbb{R}^{d}\times\mathbb{R}^{d}\times\cdots\times\mathbb{R}^{d})$ of the form $d\nu=e^\upsilon\ d\mu$ and $\nu\neq\mu$, where $\upsilon$ is a function of the variable associated with the last $\mathbb{R}^{d}$ component of $\mathcal{X}=\mathbb{R}^{d}\times\mathbb{R}^{d}\times\cdots\times\mathbb{R}^{d}$.
\end{definition}

\begin{assumption}
\label{ass:potential}
The potential \(U \in \mathscr S\) has compact level sets, satisfies
\(e^{-U}\in L^1(\mathcal X)\), and there exist \(c_U>0\), \(C_U\in\mathbb R\)
such that
\[
\theta\cdot \nabla U(\theta) \geq c_U |\theta|^2 - C_U .
\]
\end{assumption}

We have the following comparison result for the rate function of underdamped Langevin dynamics \eqref{eqn:underdamped} with overdamped Langevin dynamics \eqref{eq:overdamped-2}.

\begin{proposition}\label{underdamped VS overdamped} 
Assume Assumption \ref{ass:potential}. Let $\nu(\theta,r)$ be any measure in $\mathcal{P}_k^H(\mathcal{X}=\mathbb{R}^{2d})$ of the form $d\nu=e^\upsilon\ d\mu$ and $\nu\neq\mu$ for $\upsilon=\upsilon(r)$. Then for the underdamped Langevin dynamics defined by (\ref{eqn:underdamped}), if $\gamma\ge 1$, we have
\begin{equation}
I_u(\nu)\ge I_{e2o}(\nu),
\end{equation}
where $I_{e2o}(\nu)$ is the rate function of the 
expanded second-order overdamped Langevin dynamics (\ref{eqn:expanded-2th-overdamped}) and $I_u(\nu)$ is the rate function of the underdamped Langevin dynamics (\ref{eqn:underdamped}).
\end{proposition}%\proofat{app:pf-underdamped VS overdamped}

\begin{proof}
The proof will be provided in Appendix~\ref{app:pf-underdamped VS overdamped}.
\end{proof}

Proposition~\ref{underdamped VS overdamped} shows the acceleration
effect of the underdamped Langevin dynamics \eqref{eqn:underdamped} compared to the expanded second-order overdamped Langevin dynamics \eqref{eqn:expanded-2th-overdamped}
on a class of probability measures $\mathcal{P}^{H}(\mathcal{X}=\mathbb{R}^{2d})$.

Finally, we compare the rate functions of high-order Langevin dynamics \eqref{JI} with overdamped Langevin dynamics \eqref{eq:overdamped-2}.
Note that overdamped Langevin dynamics lives in $\mathbb{R}^{d}$, 
whereas high-order Langevin dynamics live in $\mathbb{R}^{3d}$.
Thus, we lift the overdamped Langevin dynamics to $\mathbb{R}^{3d}$, before we do the comparison.
We consider the {\it expanded third-order overdamped Langevin dynamics}, for $(\theta_t, p_t, r_t)\in\mathbb{R}^d\times\mathbb{R}^d\times\mathbb{R}^d$,
\begin{equation}\label{eqn:expanded-3th-overdamped}
\begin{cases}
d\theta_t=-\nabla U(\theta_t)dt+\sqrt{2}dW_t,\\
dp_t=-p_t dt+\sqrt{2}dB^1_t,\\
dr_t=-r_t dt+\sqrt{2}dB^2_t,
\end{cases}
\end{equation}
where $W_t$, $B^1_t$ and $B^2_t$ are independent $d$-dimensional Brownian motions. We can easily check that its invariant measure is $d\mu=e^{-U(\theta)-\frac{1}{2}|p|^2-\frac{1}{2}|r|^2}d\theta dpdr$, that is the same as the invariant measure of high-order Langevin dynamics. By LDP for the overdamped Langevin dynamics (Lemma~\ref{LDP-A}), the LDP rate function for the expanded third-order overdamped Langevin dynamics \eqref{eqn:expanded-3th-overdamped} is 
\begin{equation}\label{rate function expanded 3th overdamped}
I_{e3o}(\nu)=\frac{1}{4}\int_{\mathcal{X}}|\nabla\upsilon|^2\ d\nu, 
\end{equation}
where $\nu\in\mathcal{P}(\mathcal{X}=\mathbb{R}^{3d})$ of the form $d\nu=e^\upsilon\ d\mu$ and $\nu\neq\mu$.

\begin{proposition}\label{High-order VS overdamped} 
Assume that $U$ satisfies Assumption~\ref{HHJ}. For any $\nu(\theta,p,r)\in\mathcal{P}^H(\mathcal{X}=\mathbb{R}^{3d})$, $\nu$ is of the form $d\nu=e^\upsilon\ d\mu$ and $\nu\neq\mu$ for $\upsilon=\upsilon(r)$. Then for the high-order Langevin dynamics defined by (\ref{JI}), if $\alpha\ge 1$, then we have
\begin{equation}
I_{H}(\nu)\ge I_{e3o}(\nu),
\end{equation}
where $I_{e3o}(\nu)$ is the rate function of the 
expanded third-order overdamped Langevin dynamics (\ref{eqn:expanded-3th-overdamped}) and $I_{H}(\nu)$ is the rate function of the high-order Langevin dynamics (\ref{eqn: rate function High order}). 
\end{proposition}%\proofat{app:pf-High-order VS overdamped}

\begin{proof}
The proof will be provided in Appendix~\ref{app:pf-High-order VS overdamped}.
\end{proof}

Proposition~\ref{High-order VS overdamped} shows the acceleration
effect of the underdamped Langevin dynamics \eqref{JI} compared to the expanded third-order overdamped Langevin dynamics \eqref{eqn:expanded-3th-overdamped}
on a class of probability measures $\mathcal{P}^{H}(\mathcal{X}=\mathbb{R}^{3d})$.

\section{Numerical Experiments}\label{sec:numerical}

In this section, we conduct numerical experiments 
for the variants of Langevin algorithms based on the Euler–Maruyama
discretization of underdamped Langevin dynamics \eqref{eqn:underdamped}, non-reversible Langevin dynamics \eqref{non:reversible}, mirror Langevin dynamics \eqref{mirror Langevin dynamics},
high-order Langevin dynamics \eqref{JI}
and Hessian-free high-resolution dynamics \eqref{Hessian-free high-resolution}, 
and compare our results with the unadjusted Langevin algorithm, 
i.e. the Euler–Maruyama discretization of overdamped Langevin dynamics \eqref{eq:overdamped-2}.
In the literature, finer discretizations for Langevin algorithms have been used; see e.g. the discretization scheme in \cite{Cheng,dalalyan2018kinetic}
and the splitting schemes in \cite{leimkuhler2016molecular} for underdamped Langevin dynamics. 
However, to illustrate our theory, which is based on the continuous-time diffusion, instead of showing acceleration due to finer discretizations, we use Euler-Maruyama discretization for all the variants of 
Langevin dynamics for the sake of fair comparison.
We focus on applying variants of Langevin algorithms
to Bayesian linear regression using synthetic data and Bayesian logistic regression using both synthetic and real data. 
%%%%%%%%%%%%%%%%%%%%%%%%%%%%%%%%%%%%%%%

%%%%%%%%%%%%%%%%%%%%%%%%%%%%%%%%%
\subsection{Bayesian linear regression}

\subsubsection{Bayesian linear regression framework}

Consider a dataset $\mathcal{A}=\{(X_j,y_j)\}_{j=1}^n$, where $X_j\in\mathbb{R}^d$ are independent feature vectors and $y_j\in\mathbb{R}$ are real-valued responses. We assume the standard Gaussian noise model

\begin{equation}
    \mathbb{P}(y_j \mid X_j, x)
    \;=\;
    \frac{1}{\sqrt{2\pi\sigma^2}}
    \exp\!\left(-\frac{1}{2\sigma^2}\big(y_j-x^\top X_j\big)^2\right),
    \qquad j=1,\dots,n.
    \label{eq:linear-likelihood}
\end{equation} where $x\in\mathbb{R}^d$ denotes the regression coefficients and $\sigma^2$ is the observation-noise variance. Under the Bayesian paradigm, we impose a Gaussian prior $p(x)=\mathcal{N}(0,\lambda I)$ to regularize the parameter space. The posterior distribution $\pi(x)\propto e^{-U(x)}$ combines the likelihood and prior through the potential function
\begin{align}
    U(x)
    &=
    -\sum_{j=1}^n \log \mathbb{P}(y_j \mid X_j, x) - \log p(x)
    \nonumber\\
    &=
    \frac{1}{2\sigma^2}\sum_{j=1}^n \big(y_j-x^\top X_j\big)^2
    +\frac{\lambda}{2}\|x\|^2
    \;+\; C,
    \label{eq:blr-potential}
\end{align}
where $C$ is a constant independent of $x$.
The function $U(x)$ in \eqref{eq:blr-potential} is convex and continuously differentiable. Moreover, due to the $\ell_2$-regularization term inherited from the Gaussian prior, $U(x)$ is $\lambda$-strongly convex (and hence admits a unique minimizer), which ensures identifiability of the regression coefficients and is favorable for efficient sampling. For completeness, its gradient is
\begin{equation}
    \nabla U(x)
    \;=\;
    -\frac{1}{\sigma^2}\sum_{j=1}^n \big(y_j-x^\top X_j\big)X_j
    +\lambda x.
    \label{eq:blr-gradient}
\end{equation}

\subsubsection{Bayesian linear regression with synthetic data}
\begin{figure}[h!]
     \begin{subfigure}[b]{0.45\textwidth}
         \centering
         \includegraphics[width=0.9\linewidth]{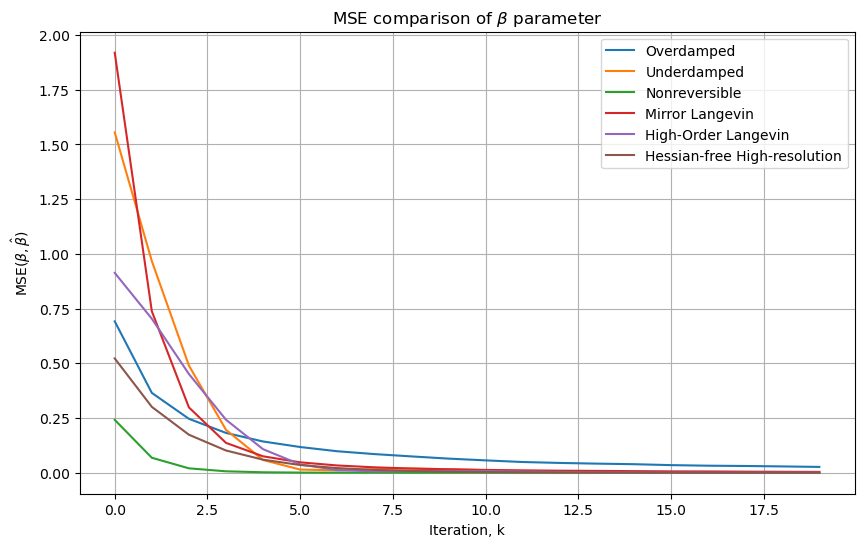}
         \caption{MSE of $\beta$ parameters}
         \label{mse_beta}
     \end{subfigure}
     \begin{subfigure}[b]{0.45\textwidth}
         \centering
         \includegraphics[width=0.9\linewidth]{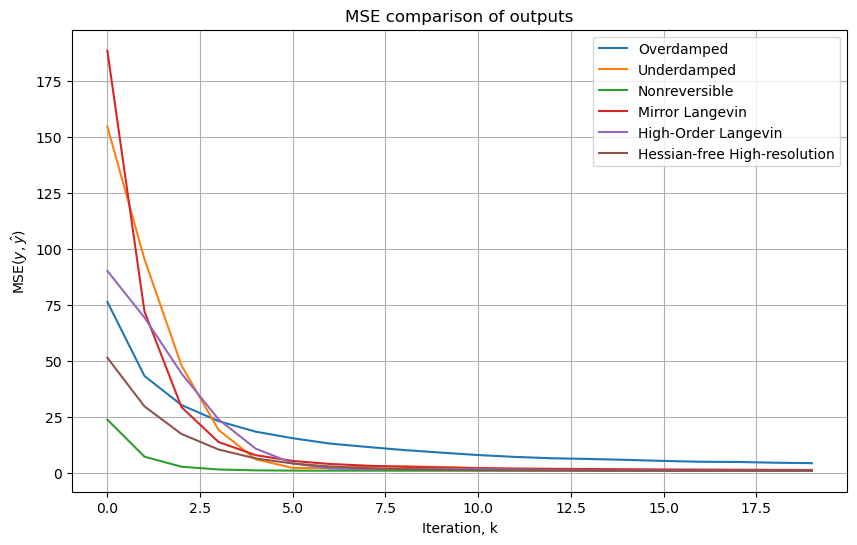}
         \caption{MSE of outputs}
         \label{mse_y}
     \end{subfigure}
     
     \hfill
        \caption{The plots show the accuracy and mean square error (MSE) over the synthetic data using Bayesian linear regression with dimension 5,000 × 100, in which all variants of the Langevin algorithms outperform overdamped Langevin algorithm in with an appropriate choice of hyper-parameters.}
        \label{graph5}
\end{figure}

We generate a synthetic dataset for Bayesian linear regression in dimension $d=100$, including an explicit intercept term. A total of $N=5000$ samples are produced.
Feature vectors are drawn independently as
\[
\tilde{X}_i \sim \mathcal{N}(0, I_{d-1}), \qquad i=1,\dots,N-1.
\]
and an intercept is appended to form
\[
X_i = (1,\tilde{X}_i^\top)^\top \in \mathbb{R}^{d}.
\]
A ground-truth parameter vector $\beta \sim \mathcal{N}(0, I_d)$ is sampled, and responses are generated according to the noisy linear model
\[
y_i = X_i^\top \beta + \xi_i,
\qquad \xi_i \sim \mathcal{N}(0,\sigma^
{2}),
\]
with noise level $\sigma = 1$.
We evaluate the algorithms using mean squared error, the parameter mean squared error $\|\beta_k-\beta\|^2$ and the prediction error $\|y - \hat{y}_{k}\|^2$. The evolution of these quantities is shown in Figure~\ref{mse_beta} and Figure~\ref{mse_y}. Across all experiments, the proposed Langevin variants  outperform the overdamped Langevin dynamics, achieving faster convergence and  lower errors in both parameter recovery and output prediction.
For non-reversible Langevin, we generated the anti-symmetric matrix $J$ by first generating a $d\times d$ random sparse matrix $A$ with standard normal random entries and then we set $J = A - A^{\top}$ such that $J$ is an anti-symmetric matrix i.e $J^{\top} = -J$; for mirror Langevin, we set the mirror map as $h(x) = \frac{1}{4} \sum_{i = 1}^{d} x_{i}^{4}$ for any $x=(x_{1},\ldots,x_{d})\in\mathbb{R}^{d}$. We choose the step size $\eta = 0.0003$ for overdamped,  $\eta = 0.0001$ for non-reversible and Hessian-free high-resolution,  $\eta = 0.005$ for underdamped and high order Langevin, and $\eta = 0.00003$ for mirror Langevin.
We chose $\gamma = 75$ for underdamped Langevin, $\gamma = 150$ and $\alpha = 220$ for high-order Langevin, and $\gamma = 0.5$ and $\alpha = 10$ for Hessian-free high-resolution. We can see from this Figure~\ref{mse_beta} and Figure~\ref{mse_y} that Overdamped Langevin is outperformed by all the other variants.

%%%%%%%%%%%%%%%%%%%%%%%%%%%%%%%%%%%%%%%%%%%%%%%%%%%%%

\subsection{Bayesian logistic regression}

\subsubsection{Bayesian logistic regression framework}

Consider a dataset $\mathcal{A} = \{(X_j, y_j)\}_{j=1}^n$, where $X_j \in \mathbb{R}^d$ denotes independent feature vectors and $y_j \in \{0,1\}$ represents binary labels. The likelihood function follows the logistic model:
\begin{equation}
    \mathbb{P}(y_j = 1 | X_j, x) = \sigma(x^\top X_j) = \left(1 + e^{-x^\top X_j}\right)^{-1},
    \label{eq:sigmoid-model}
\end{equation}
where $\sigma(\cdot)$ denotes the sigmoid function and $x \in \mathbb{R}^{d}$ denotes the regression coefficients.

Under the Bayesian paradigm, we impose a Gaussian prior $p(x) = \mathcal{N}(0, \lambda I)$ with $\lambda = 10$ to regularize the parameter space. The posterior distribution $\pi(x) \propto e^{-U(x)}$ combines the likelihood and prior through the potential function:
\begin{equation}
    U(x) = -\sum_{j=1}^n \log \mathbb{P}(y_j | X_j, x) - \log p(x)= \sum_{j=1}^n \log\left(1 + e^{-x^\top X_j}\right) + \frac{\lambda}{2}\|x\|^2.
    \label{eq:potential-function}
\end{equation}
Note that the function $U(x)$ in \eqref{eq:potential-function} is strongly convex and smooth due to the $\ell_2$-regularization term inherited from the Gaussian prior. This guarantees unique identifiability of the parameters and facilitates efficient sampling.

%%%%%%%%%%%%%%%%%%%%%%%%%%%%%%%%%%%%%%%%%%%%%%%%%%%%%%%%%%%%%%%%
\subsubsection{Bayesian logistic regression with synthetic data}

In this section, we validate our methodology through synthetic data experiments.
Let feature vectors $X_j \sim \mathcal{N}(0, 10I_{d})$,  true parameters $x = [x_1, x_2,\dots,  x_d]^\top$ with prior $x \sim \mathcal{N}(0, 10 I_{d})$ and binary labels for each $(X_j, p_j)$ where $p_j \sim \mathcal{U}[0,1]$, 
and $y_j = 1$ if $p_j \leq \sigma(x_{\text{true}}^\top X_j)$
and $y_j=0$ otherwise.
This data generation process ensures the labels $y_j$ adhere to the logistic model while maintaining controlled experimental conditions. The uniform threshold introduces stochasticity consistent with the Bernoulli likelihood structure.

%%%%%%%%%%%%%%%%%%%%%%%%%%%%%%%%%%%%%%%%%%%%%%%%%%%%%%%%%%%%%%%%%%%%%%%%%%%%%%%%%%%%%%%%%
For synthetic data, we generated $5000$ data points with $31$ random features and generated binary labels. 
For synthetic data, for non-reversible Langevin, we generated the anti-symmetric matrix $J$ by first generating a $d\times d$ random matrix $A$ with standard normal random entries and then we set $J = A - A^{\top}$ such that $J$ is an anti-symmetric matrix i.e $J^{\top} = -J$; for mirror Langevin, we set the mirror map as $h(x) = \frac{1}{4} \sum_{i = 1}^{d} x_{i}^{4}$ for any $x=(x_{1},\ldots,x_{d})\in\mathbb{R}^{d}$. We chose the stepsize $\eta = 0.0003$ for overdamped, non-reversible, Hessian-free high-resolution and mirror Langevin and the stepsize $\eta = 0.003$ for underdamped and high-order Langevin. 
We chose $\gamma = 4$ for underdamped Langevin, $\gamma = 20$ and $\alpha = 15$ for high-order Langevin, and $\gamma = 1$ and $\alpha = 30$ for Hessian-free high-resolution. 

%%%%%%%%%%%%%%%%%%%%%%%%%%%%%%%%%%%%%%%%%%%%%%%%%%%%%%%%%%%%%%%%%%%%%%%%%%%%%%%%%

\begin{figure}[!h]
     \begin{subfigure}[b]{0.3\textwidth}
         \centering
         \includegraphics[width=0.9\linewidth]{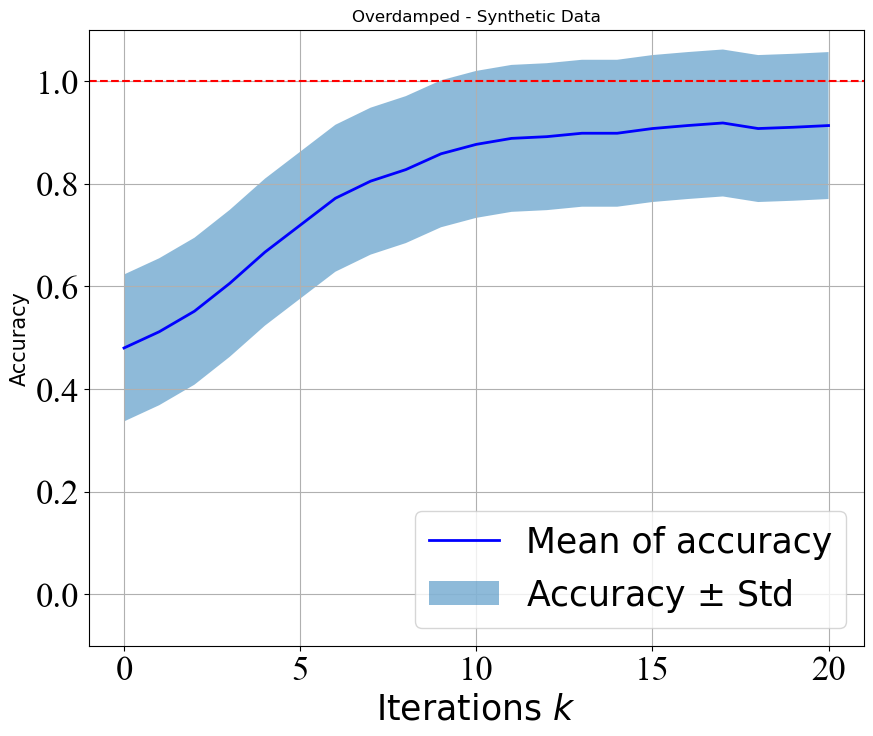}
         \caption{Overdamped Langevin}
         \label{fig:OD_s}
     \end{subfigure}
     \hfill
     \begin{subfigure}[b]{0.3\textwidth}
         \centering
         \includegraphics[width=0.9\linewidth]{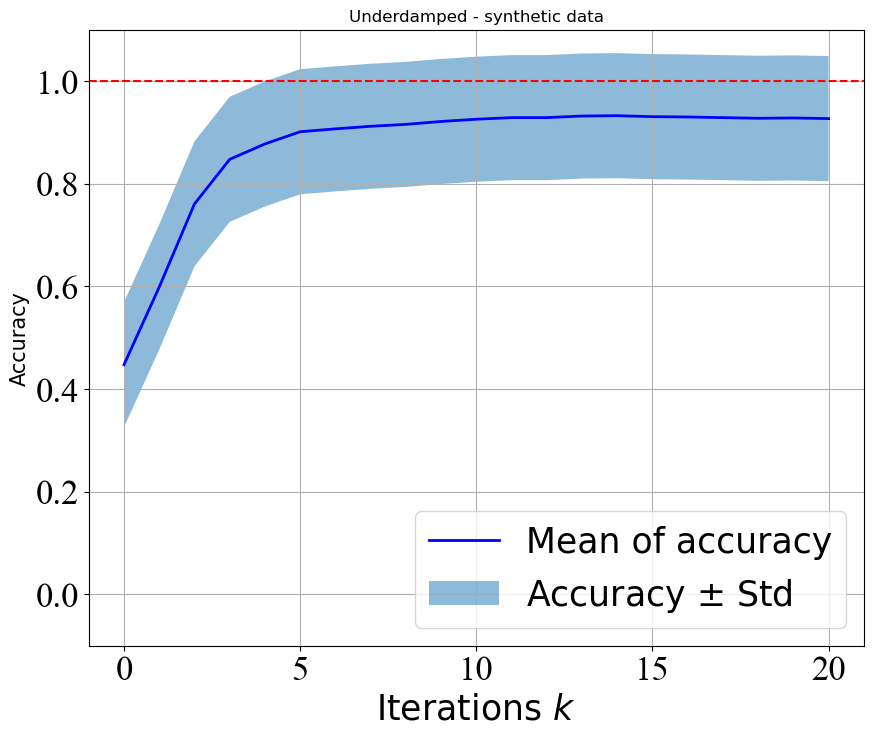}
         \caption{Underdamped Langevin}
         \label{fig:UD_s}
     \end{subfigure}
     \hfill
     \begin{subfigure}[b]{0.3\textwidth}
         \centering
         \includegraphics[width=0.9\linewidth]{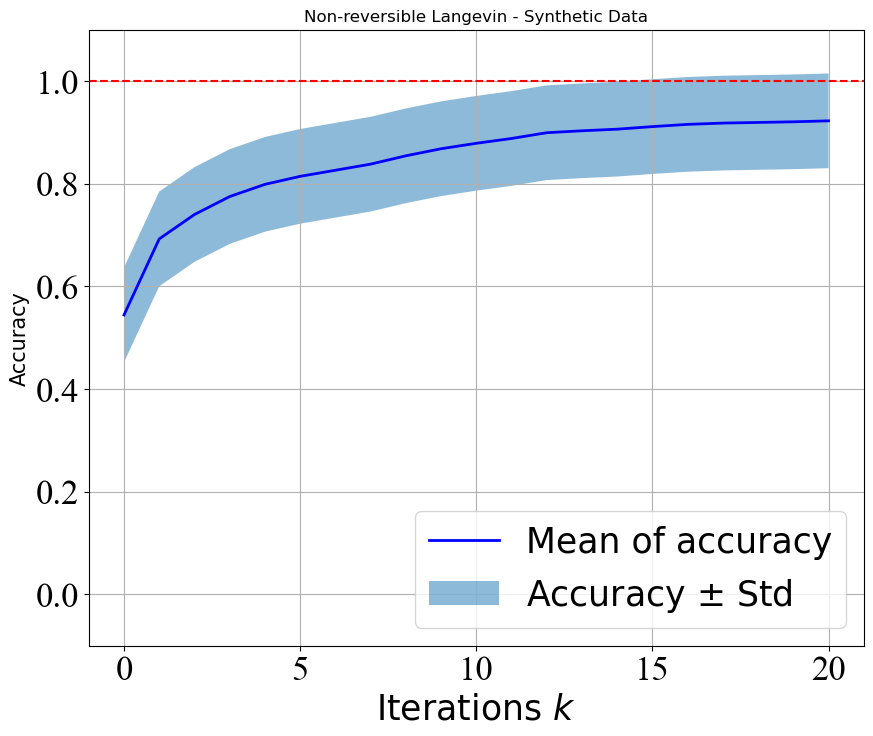}
         \caption{Non-reversible Langevin}
         \label{fig:NR_s}
     \end{subfigure}
\vspace{.6ex}
     \begin{subfigure}[b]{0.3\textwidth}
         \centering
         \includegraphics[width=0.9\linewidth]{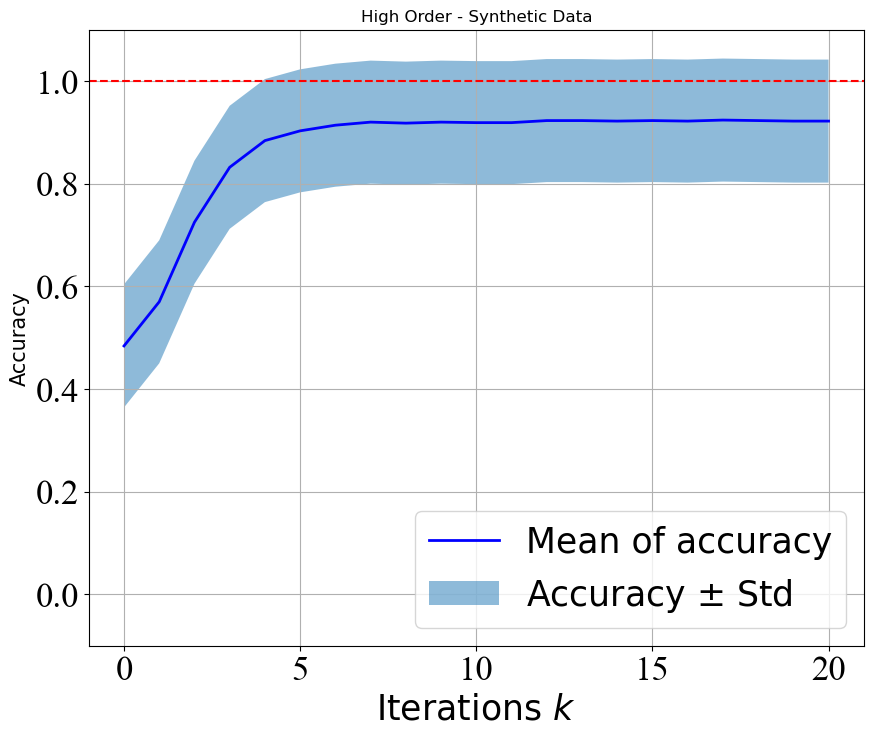}
         \caption{High-order Langevin}
         \label{fig:HO_s}
     \end{subfigure}
     \hfill
     \begin{subfigure}[b]{0.33\textwidth}
         \centering
         \includegraphics[width=0.9\linewidth]{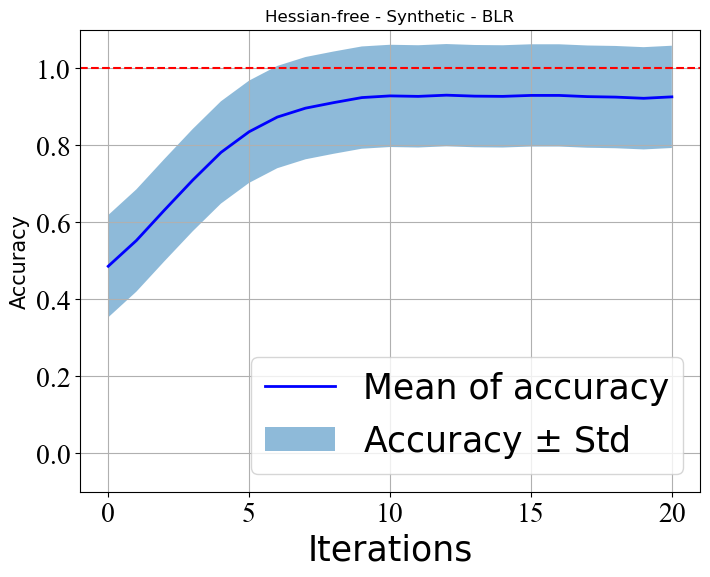}
         \caption{Hessian-free high-resolution}
         \label{fig:HFHR_s}
     \end{subfigure}
      \hfill
     \begin{subfigure}[b]{0.3\textwidth}
         \centering
         \includegraphics[width=0.9\linewidth]{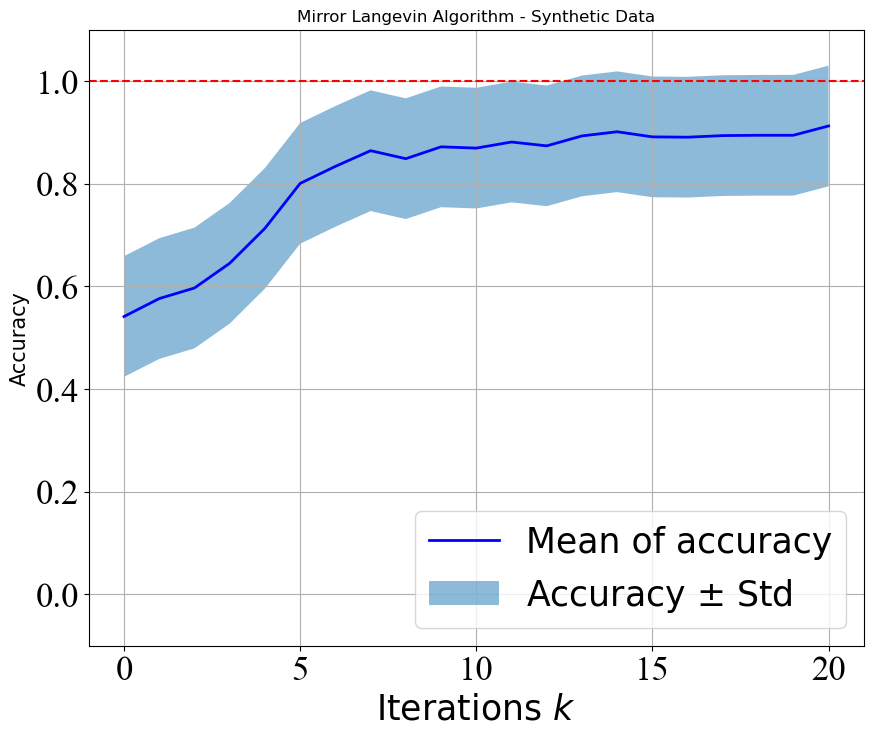}
         \caption{Mirror Langevin}
         \label{fig:ML_s}
     \end{subfigure}

     \hfill
     
        \caption{The plots show the accuracy over the synthetic data with dimension 5,000 $\times$ 31, in which all variants of the Langevin algorithms outperform overdamped Langevin algorithm in Figure~\ref{fig:OD_s} with an appropriate choice of hyperparameters.}
        \label{fig:graphs1}
\end{figure}

\begin{figure}[!h]
     \begin{subfigure}[b]{0.3\textwidth}
         \centering
         \includegraphics[width=0.9\linewidth]{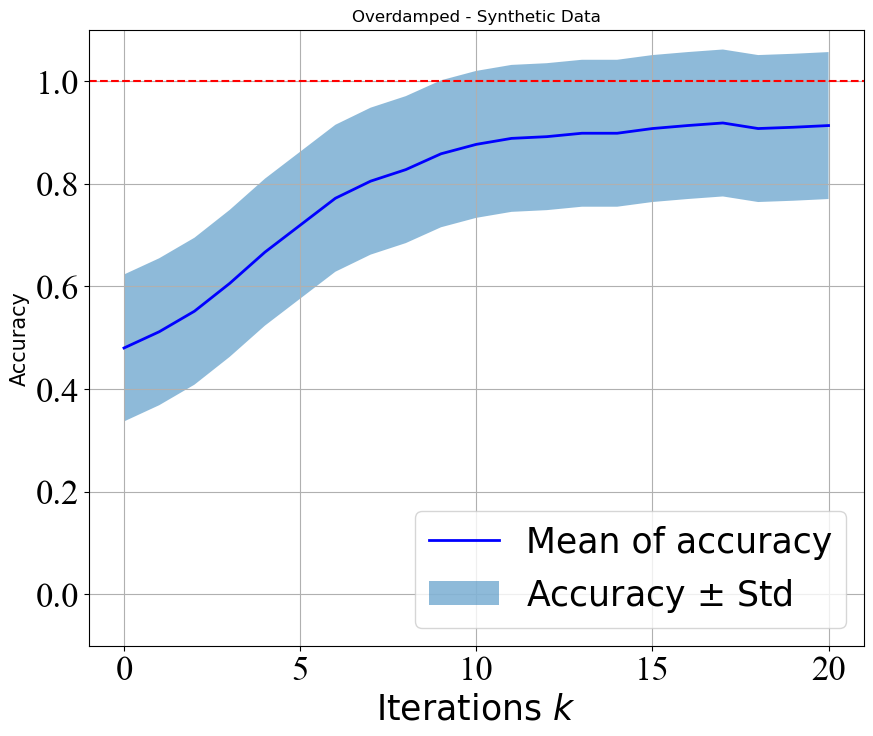}
         \caption{Overdamped Langevin}
         \label{fig:OD2_s}
     \end{subfigure}
     \hfill
     \begin{subfigure}[b]{0.3\textwidth}
         \centering
         \includegraphics[width=0.9\linewidth]{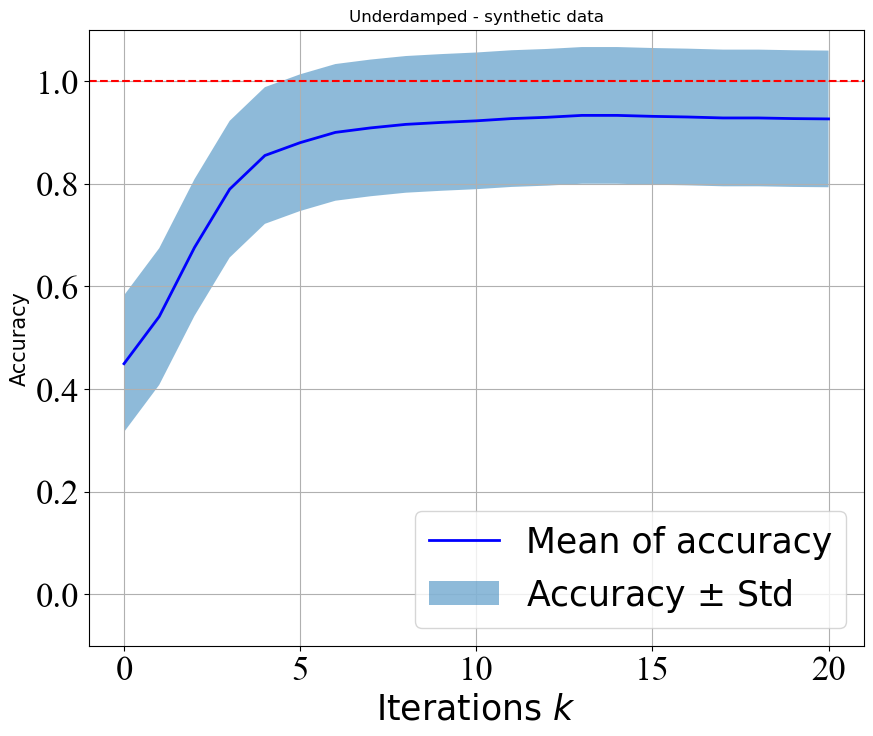}
         \caption{Underdamped Langevin}
         \label{fig:UD2_s}
     \end{subfigure}
     \hfill
     \begin{subfigure}[b]{0.3\textwidth}
         \centering
         \includegraphics[width=0.9\linewidth]{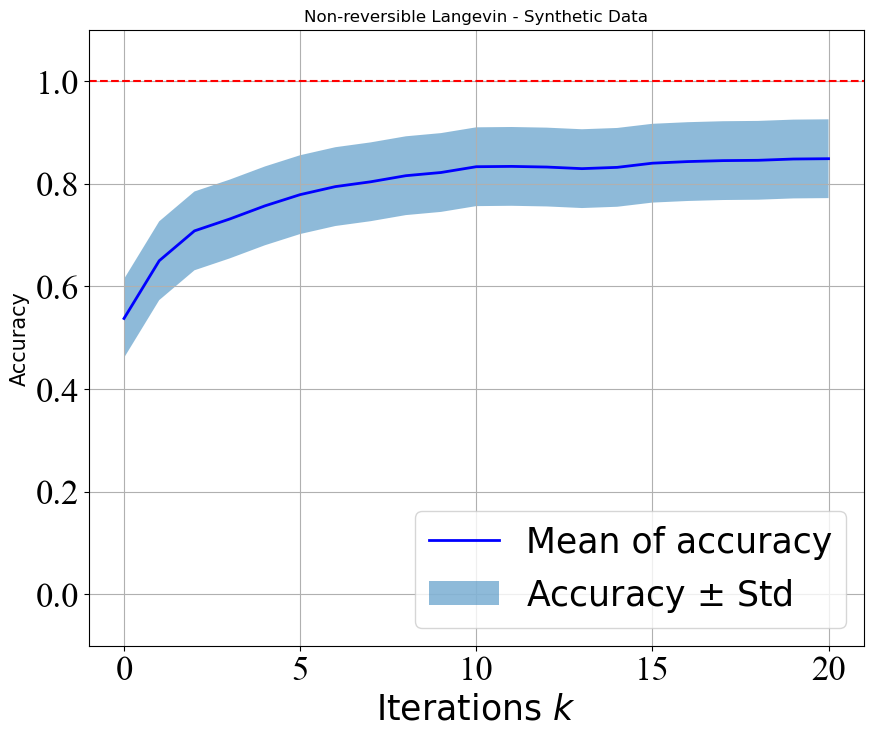}
         \caption{Non-reversible Langevin}
         \label{fig:NR2_s}
     \end{subfigure}
\vspace{.6ex}
     \begin{subfigure}[b]{0.3\textwidth}
         \centering
         \includegraphics[width=0.9\linewidth]{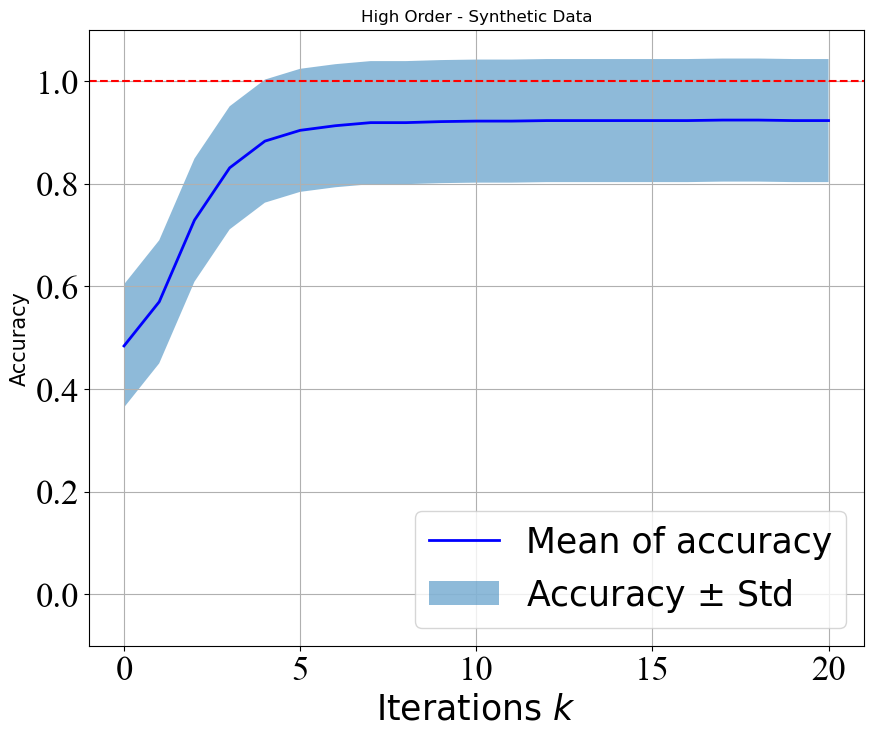}
         \caption{High-order Langevin}
         \label{fig:HO2_s}
     \end{subfigure}
     \hfill
     \begin{subfigure}[b]{0.33\textwidth}
         \centering
         \includegraphics[width=0.9\linewidth]{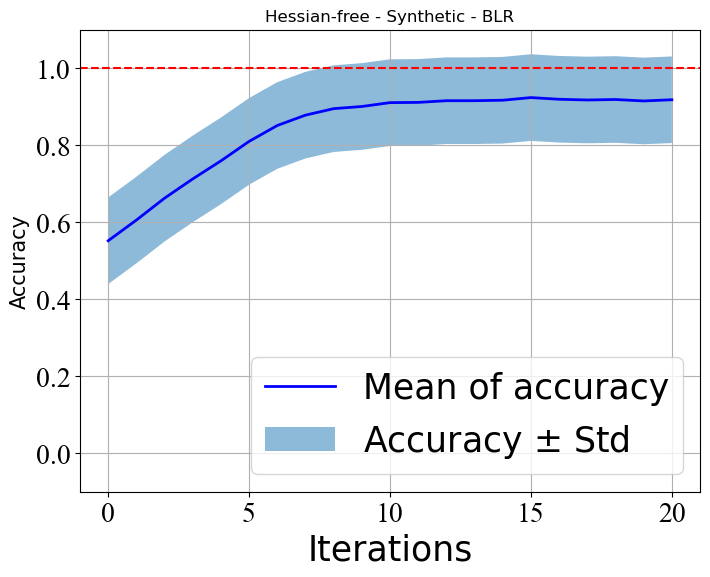}
         \caption{Hessian-free high-resolution}
         \label{fig:HFHR2_s}
     \end{subfigure}
      \hfill
     \begin{subfigure}[b]{0.3\textwidth}
         \centering
         \includegraphics[width=0.9\linewidth]{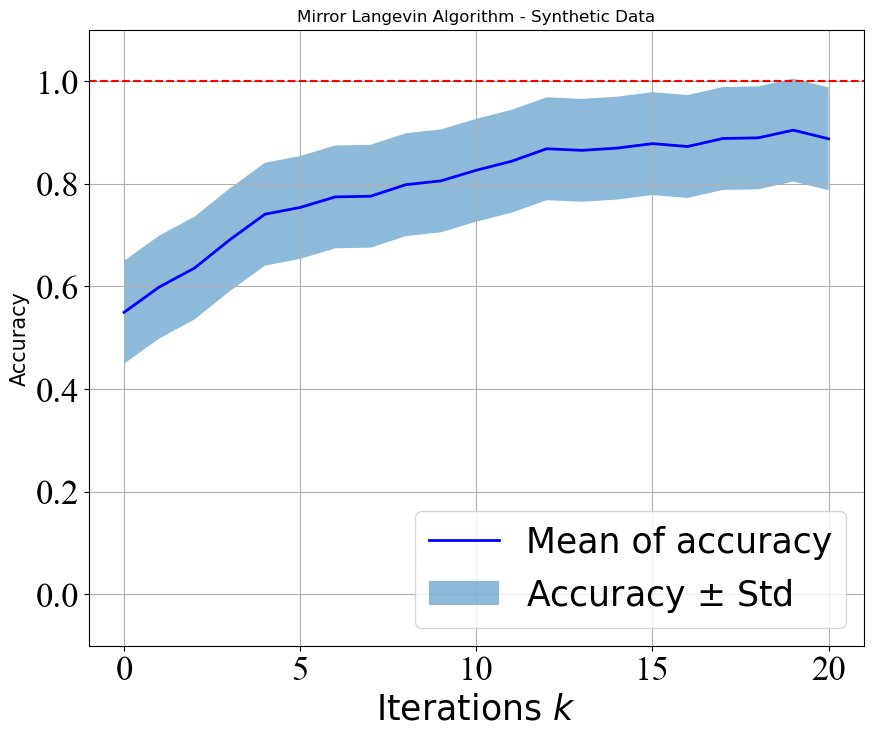}
         \caption{Mirror Langevin}
         \label{fig:ML2_s}
     \end{subfigure}
     \hfill
     
        \caption{With a slight change of hyperparameters, we can see from this figure that underdamped Langevin (Figure~\ref{fig:UD2_s}), high-order Langevin (Figure~\ref{fig:HO2_s}) and Hessian-free high-resolution (Figure~\ref{fig:HFHR2_s}) can outperform overdamped Langevin (Figure~\ref{fig:OD2_s}); however, mirror Langevin (Figure~\ref{fig:ML2_s}) and non-reversible (Figure~\ref{fig:NR2_s}) cannot, even though their performance is comparable with overdamped Langevin (Figure~\ref{fig:OD2_s}).}
        \label{graphs2}
\end{figure}

The accuracy is reported only for test data sets. Our numerical results using synthetic data, summarized in Figure~\ref{fig:graphs1}, show that for proper choices of hyperparameters (and the anti-symmetric matrix for non-reversible Langevin and the mirror map for mirror Langevin), the variants of Langevin algorithms, i.e., underdamped Langevin (Figure~\ref{fig:UD_s}), non-reversible Langevin (Figure~\ref{fig:NR_s}), high-order Langevin (Figure~\ref{fig:HO_s}), Hessian-free high-resolution (Figure~\ref{fig:HFHR_s}) and mirror Langevin (Figure~\ref{fig:ML_s}) can all have faster convergence than overdamped Langevin (Figure~\ref{fig:OD_s}).
%%%%%%%%%%%%%%%%%%%%%%%%%%%%%%%%%%%%%%%%%%%%%%%%%%%%%%%%%%%%%%%%%%%%%%%%%%%%
On the other hand, the performance of these variants of Langevin algorithms
is sensitive to the choices of hyperparameters (and the anti-symmetric matrix for non-reversible Langevin and the mirror map for mirror Langevin). 
In another set of experiments reported in Figure~\ref{graphs2}, 
with a slight change in the choices of hyperparameters, underdamped Langevin (Figure~\ref{fig:UD2_s}), high-order Langevin (Figure~\ref{fig:HO2_s}) and Hessian-free high-resolution (Figure~\ref{fig:HFHR2_s}) can outperform overdamped Langevin (Figure~\ref{fig:OD2_s}), whereas mirror Langevin (Figure~\ref{fig:ML2_s}) and non-reversible (Figure~\ref{fig:NR2_s}) cannot, even though their performance is comparable with overdamped Langevin (Figure~\ref{fig:OD2_s}). One possible explanation is that
the performance of non-reversible Langevin dynamics \eqref{non:reversible} depends on the choice of the anti-symmetric matrix, 
and the performance of mirror Langevin dynamics \eqref{mirror Langevin dynamics} depends on the choice of the mirror map,
such that when the hypermarameters are changed, one has to fine tune
the choices of the anti-symmetric matrix and the mirror map simultaneously
in order to maintain the good performance.

%%%%%%%%%%%%%%%%%%%%%%%%%%%%%%%%%%%%%%%%%%%%%%%%%%%%%%%%%
\subsection{Bayesian logistic regression with real data}

In this section, we validate our methodology through real data experiments.
We consider the UCI ML Breast Cancer Wisconsin (Diagnostic) data set \cite{misc_breast_cancer_wisconsin_(diagnostic)_17}. The data set contains $569$ instances with $31$ features, where each sample describes characteristics of the cell nuclei present in a digitized image of a fine needle aspirate (FNA) of a breast mass.

\begin{figure}[!h]
     \begin{subfigure}[b]{0.3\textwidth}
         \centering
         \includegraphics[width=0.9\linewidth]{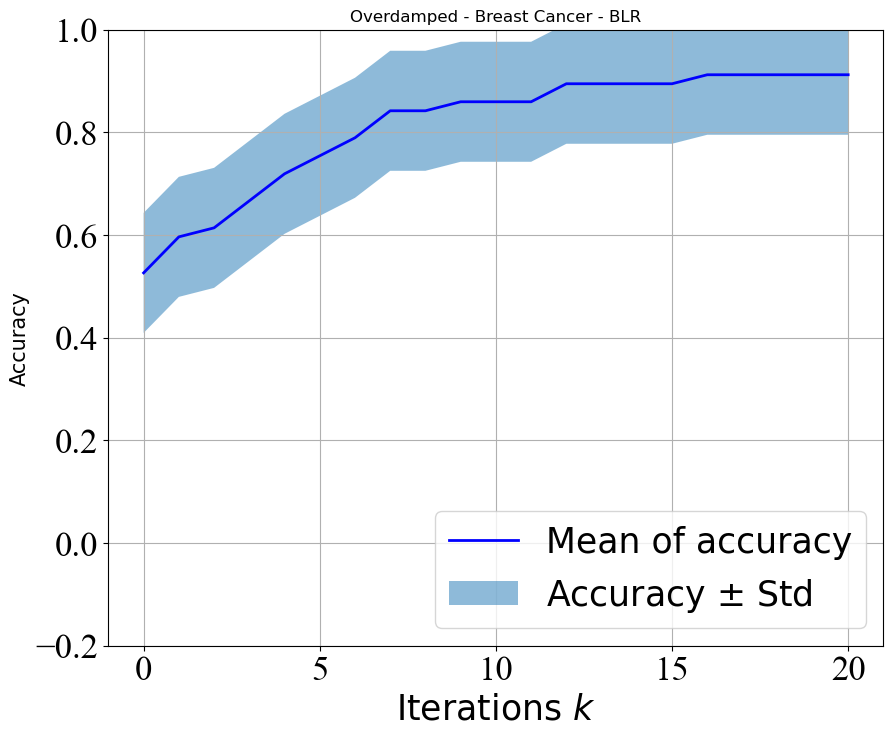}
         \caption{Overdamped Langevin}
         \label{fig:OD}
     \end{subfigure}
     \hfill
     \begin{subfigure}[b]{0.3\textwidth}
         \centering
         \includegraphics[width=0.9\linewidth]{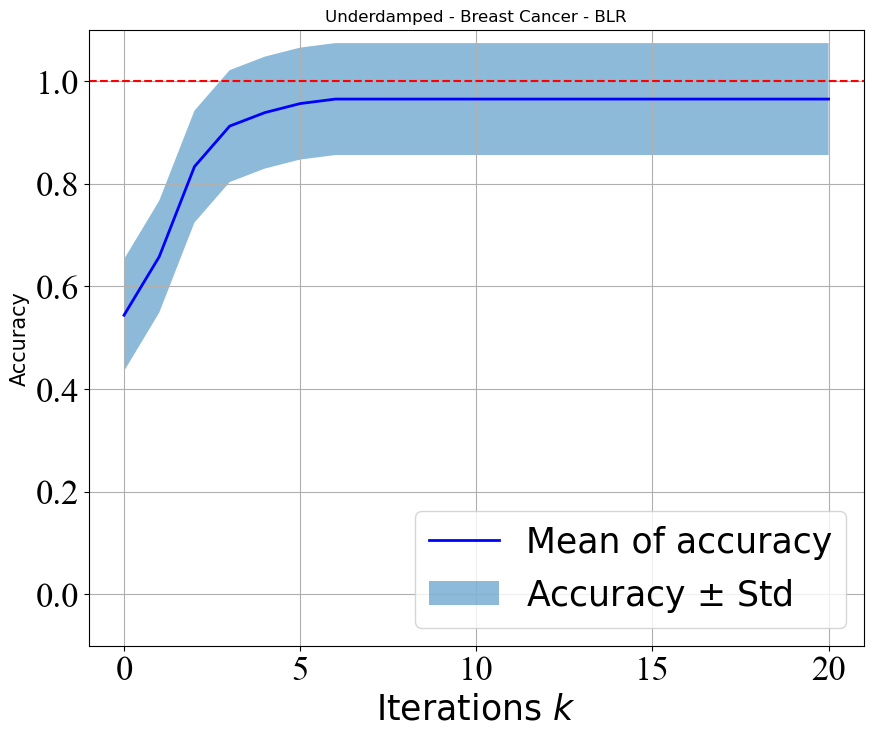}
         \caption{Underdamped Langevin}
         \label{fig:UD}
     \end{subfigure}
     \hfill
     \begin{subfigure}[b]{0.3\textwidth}
         \centering
         \includegraphics[width=0.9\linewidth]{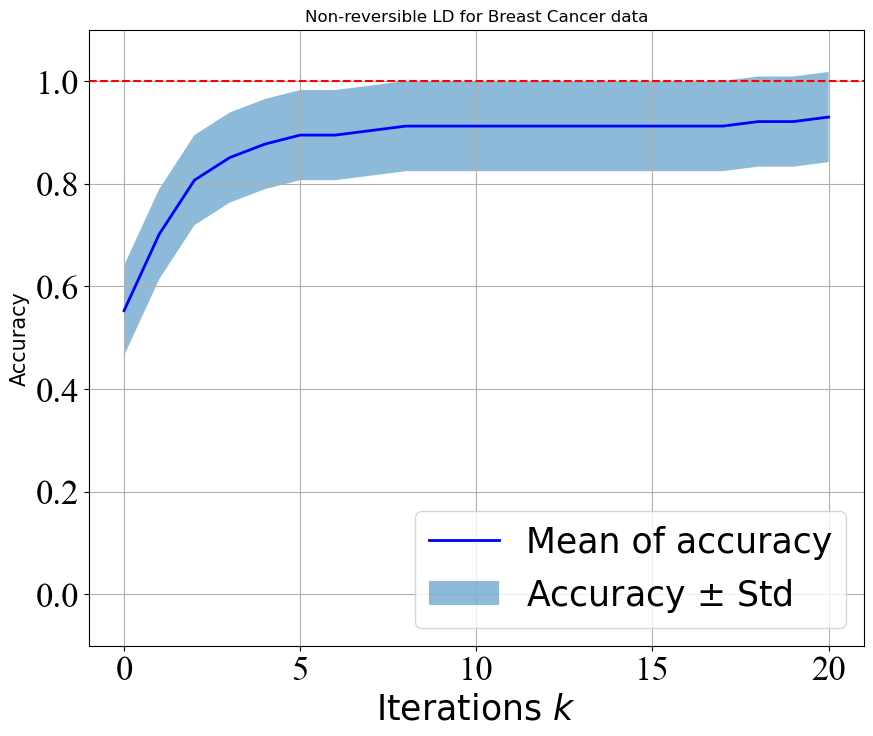}
         \caption{Non-reversible Langevin}
         \label{fig:NR}
     \end{subfigure}
\vspace{.6ex}
     \begin{subfigure}[b]{0.3\textwidth}
         \centering
         \includegraphics[width=0.9\linewidth]{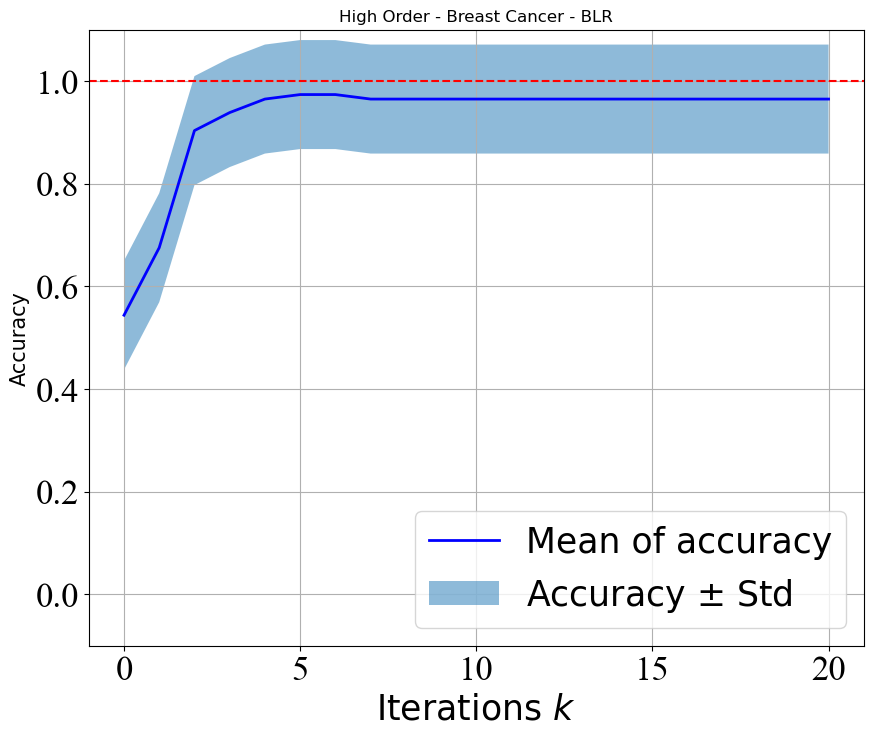}
         \caption{High-order Langevin}
         \label{fig:HO}
     \end{subfigure}
     \hfill
     \begin{subfigure}[b]{0.32\textwidth}
         \centering
         \includegraphics[width=0.9\linewidth]{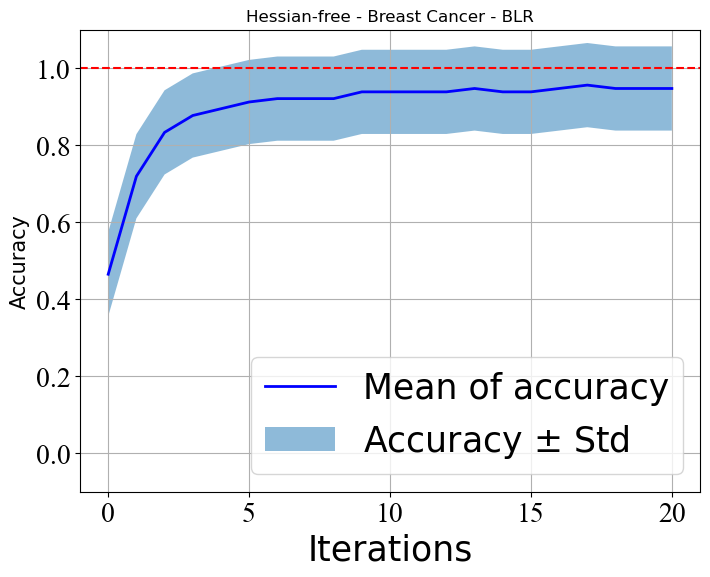}
         \caption{Hessian-free high-resolution}
         \label{fig:HFHR}
     \end{subfigure}
     \hfill
     \begin{subfigure}[b]{0.3\textwidth}
         \centering
         \includegraphics[width=0.9\linewidth]{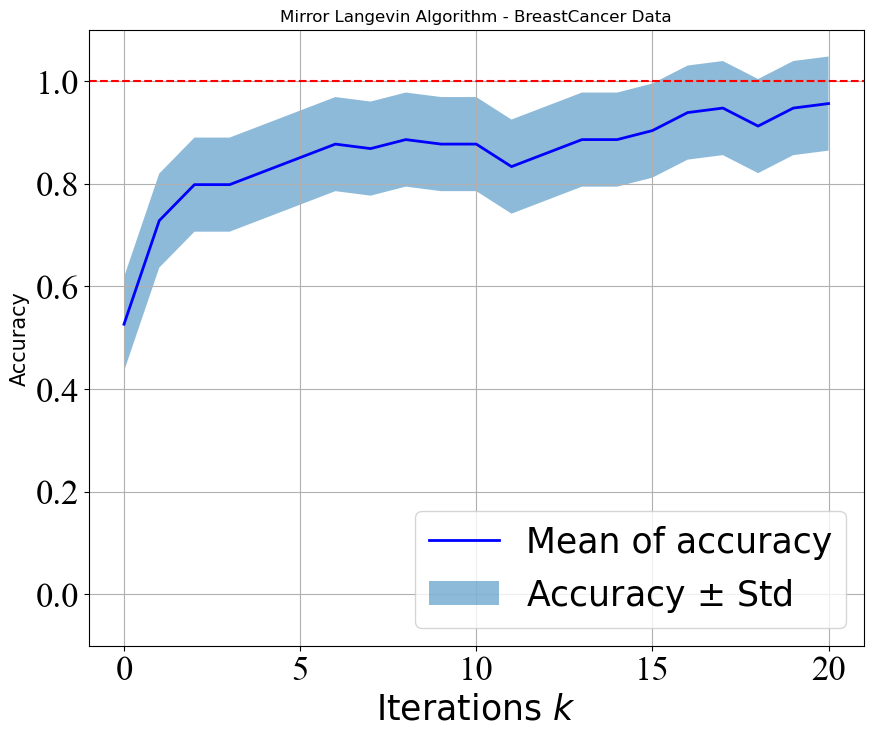}
         \caption{Mirror Langevin}
         \label{fig:ML}
     \end{subfigure}
     \hfill
        \caption{The plots show the accuracy over the real data with dimension 569 $\times$ 31, in which all variants of the Langevin algorithms outperform overdamped Langevin algorithm in Figure~\ref{fig:OD} with an appropriate choice of hyperparameters.}
        \label{fig:graph3}
\end{figure}

\begin{figure}[!h]
     \begin{subfigure}[b]{0.3\textwidth}
         \centering
         \includegraphics[width=0.9\linewidth]{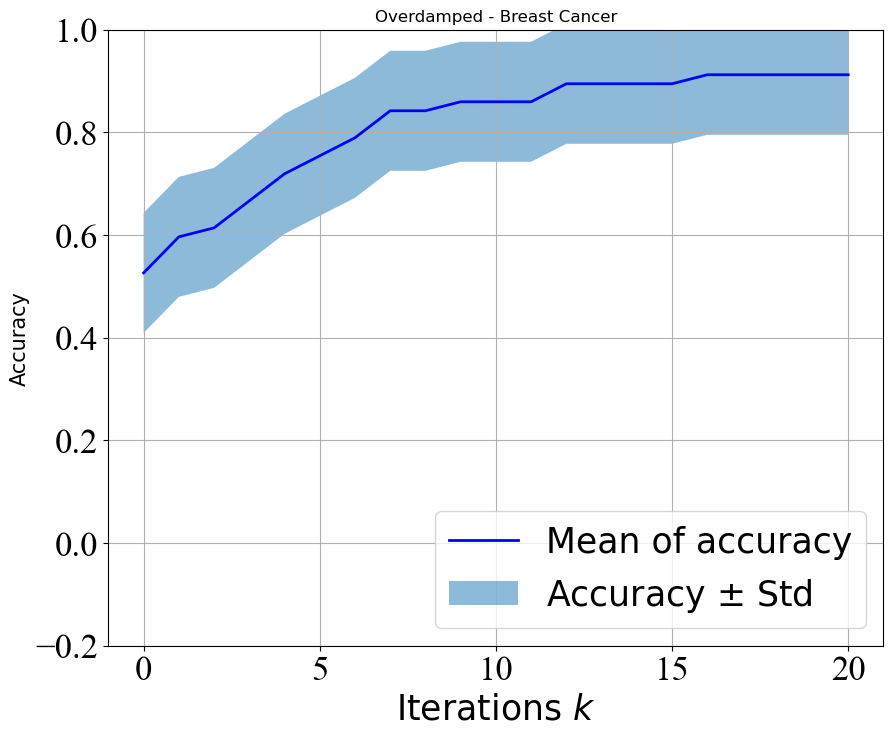}
         \caption{Overdamped Langevin}
         \label{fig:OD2}
     \end{subfigure}
     \hfill
     \begin{subfigure}[b]{0.3\textwidth}
         \centering
         \includegraphics[width=0.9\linewidth]{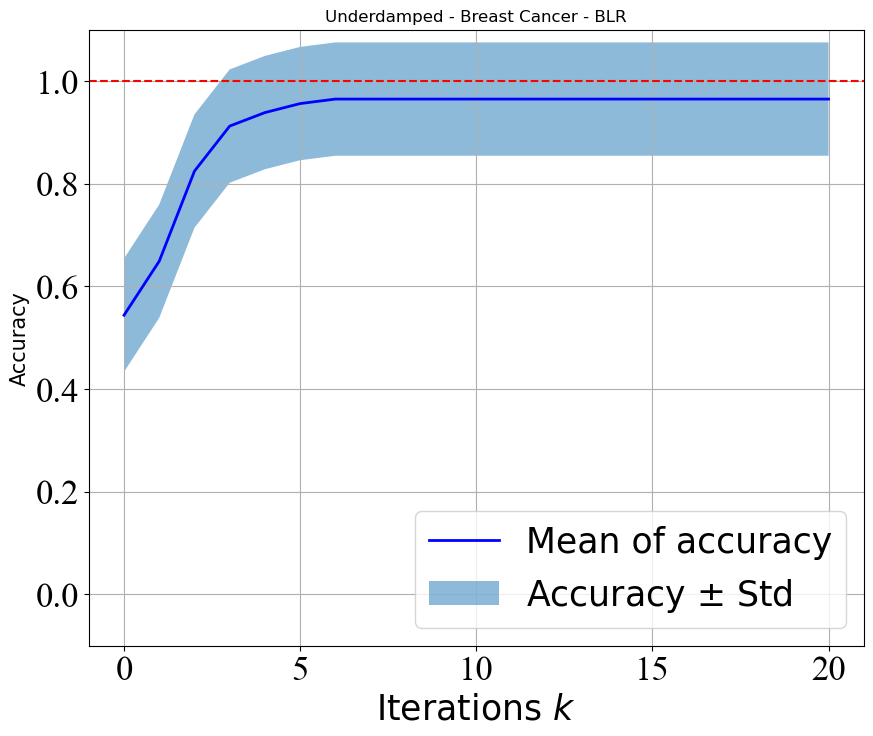}
         \caption{Underdamped Langevin}
         \label{fig:UD2}
     \end{subfigure}
     \hfill
     \begin{subfigure}[b]{0.3\textwidth}
         \centering
         \includegraphics[width=0.9\linewidth]{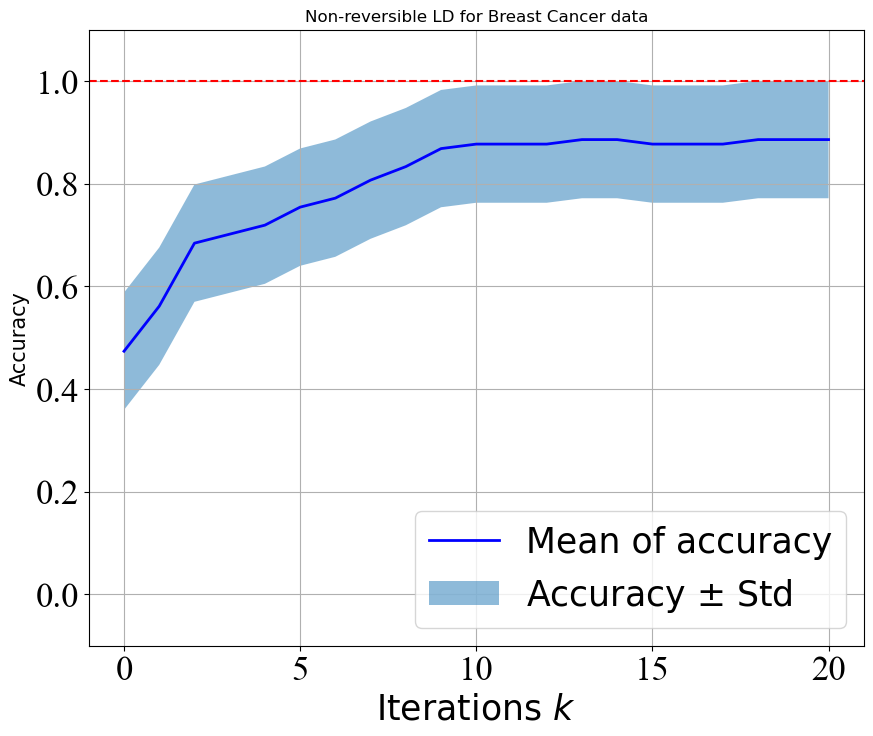}
         \caption{Non-reversible Langevin}
         \label{fig:NR2}
     \end{subfigure}
\vspace{.6ex}
     \begin{subfigure}[b]{0.3\textwidth}
         \centering
         \includegraphics[width=0.9\linewidth]{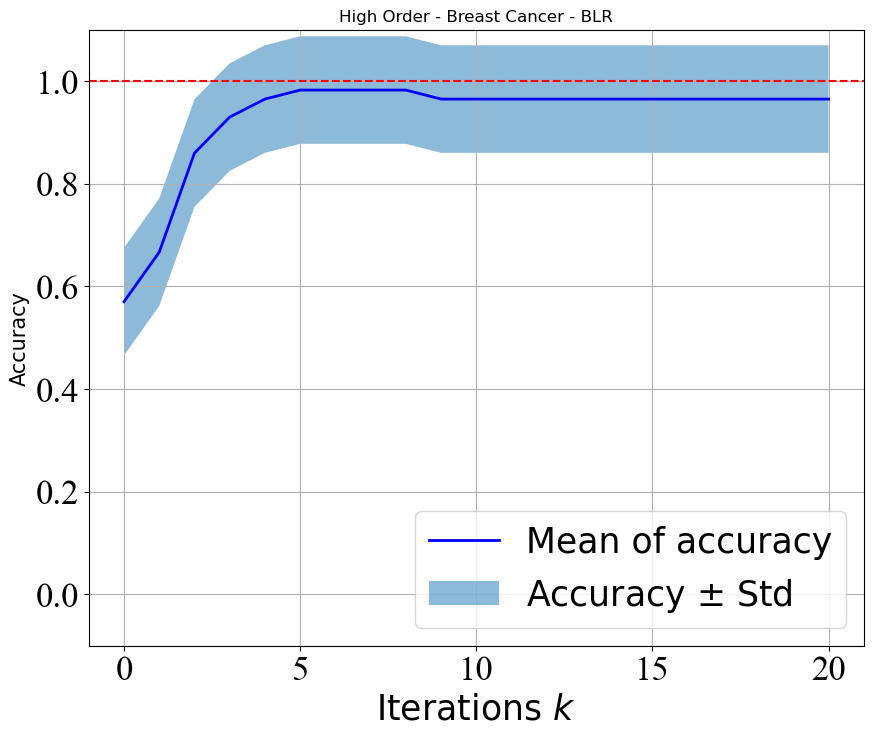}
         \caption{High-order Langevin}
         \label{fig:HO2}
     \end{subfigure}
     \hfill
     \begin{subfigure}[b]{0.32\textwidth}
         \centering
         \includegraphics[width=0.9\linewidth]{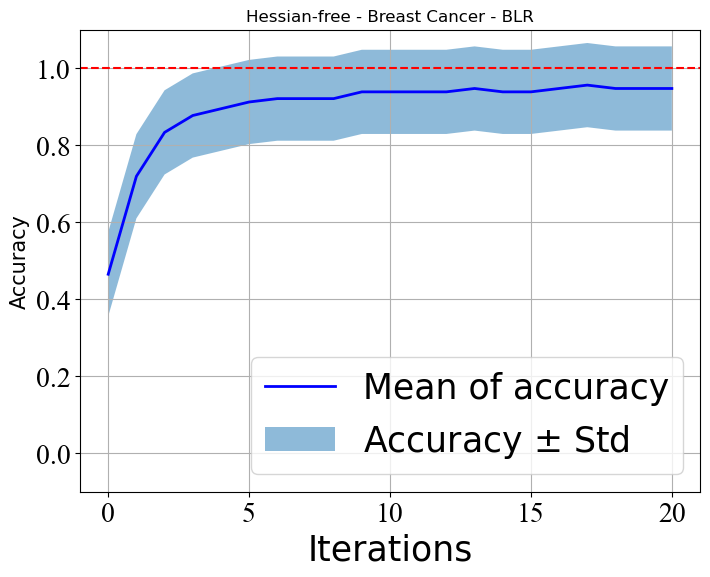}
         \caption{Hessian-free high-resolution}
         \label{fig:HFHR2}
     \end{subfigure}
     \hfill
     \begin{subfigure}[b]{0.3\textwidth}
         \centering
         \includegraphics[width=0.9\linewidth]{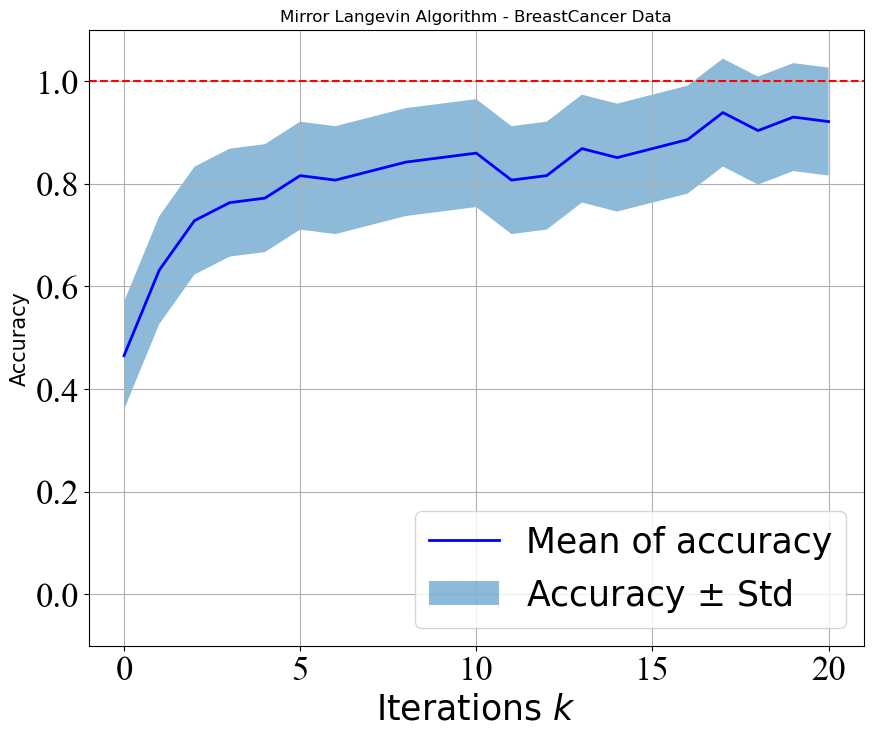}
         \caption{Mirror Langevin}
         \label{fig:ML2}
     \end{subfigure}
     \hfill
        \caption{With a slight change of hyperparameters, we can see from this figure that underdamped Langevin (Figure~\ref{fig:UD2}), high-order Langevin (Figure~\ref{fig:HO2}) and Hessian-free high-resolution (Figure~\ref{fig:HFHR2_s}) can outperform overdamped Langevin (Figure~\ref{fig:OD2}); however, mirror Langevin (Figure~\ref{fig:ML2}) and non-reversible (Figure~\ref{fig:NR2}) cannot, even though their performance is comparable with overdamped Langevin (Figure~\ref{fig:OD2}).}
        \label{graph4}
\end{figure}

%%%%%%%%%%%%%%%%%%%%%%%%%%%%%%%%%%%%%%%%%%%%%%%%%%%%%%%%%%%%%%%%%%%%%%%%%%%

For real data, we generated the anti-symmetric matrix for non-reversible Langevin and the mirror map $h(x)$ for mirror Langevin similarly as for synthetic data, and kept the same stepsizes as well. 
We chose $\gamma = 35$ for underdamped Langevin, $\gamma = 35$ and $\alpha = 35$ for high-order Langevin, and $\gamma = 1$ and $\alpha = 30$ for Hessian-free high-resolution. 
%%%%%%%%%%%%%%%%%%%%%%%%%%%%%%%%%%%%%%%%%%%%%%%%%%%%%%%%%%%%%%%%%%%%%%%%%%%%%%%%%%%%%%%%%%%%%%

The accuracy is reported only for test data sets. Our numerical results using real data show that for some particular selection of hyperparameters, the variants of Langevin algorithms can all outperform the overdamped Langevin algorithm (Figure~\ref{fig:graph3}). 
These experiments demonstrate the practical applicability of our methods and validate their performance on real-world classification tasks.
%%%%%%%%%%%%%%%%%%%%%%%%%%%%%%%%
On the other hand, under a different choice of hyperparameters as in Figure~\ref{graph4},
underdamped, high-order Langevin and Hessian-free high-resolution algorithms have faster convergence than overdamped Langevin; see Figures~\ref{fig:UD2}, \ref{fig:HFHR2}, \ref{fig:HO2} and \ref{fig:OD2}; 
whereas non-reversible Langevin, and mirror Langevin exhibit similar performance
compared with overdamped Langevin; see Figures~\ref{fig:NR2},  and \ref{fig:OD2}.
Non-reversible and mirror Langevin seem to be more sensitive to the choice of hyperparameters,
which might be due to the fact that their performance also depends on the choice of the anti-symmetric matrix and the mirror map.

%%%%%%%%%%%%%%%%%%%%%%%%%%%%%%%%%
\section{Concluding Remarks}\label{sec:conclude}

In this paper, we studied variants of Langevin dynamics through the lens of large deviations theory. 
We showed the acceleration of convergence of the variants of overdamped Langevin dynamics to the Gibbs distribution, including 
underdamped Langevin dynamics, non-reversible Langevin dynamics, 
mirror Langevin dynamics, high-order Langevin dynamics, Hessian-free
high-resolution dynamics by comparing the rate functions from the large deviations theory. 
We provided numerical experiments using both synthetic and real data, based
on the Euler–Maruyama discretizations of these variants of Langevin dynamics
and demonstrated their efficiency.

\section*{Acknowledgments}
The first author would like to thank Prof. Liming Wu, Prof. Fuqing Gao and Prof. Qi Lv for the helpful discussions.
The second and fourth authors would like to thank Prof. Mert G\"{u}rb\"{u}zbalaban
and Dr. Yuanhan Hu.
Nian Yao was supported in part by the Natural Science Foundation of China under Grant 12071361, the Natural Science Foundation of Guangdong Province under Grant 2020A1515010822
and Shenzhen Natural Science Fund (the Stable Support Plan Program 20220810152104001). Lingjiong Zhu is partially supported by the NSF grants DMS-2053454, DMS-2208303 and a Simons Foundation
Collaboration Grant.

\bibliographystyle{alpha} % Style BST file (imsart-number.bst or imsart-nameyear.bst)
\bibliography{langevin}       % Bibliography file (usually '*.bib')

\newpage

%\begin{center}
%\Large \bf Accelerating Langevin Monte Carlo Sampling: A Large Deviations Analysis \vspace{3pt}\\ {\normalsize Supplementary Document}
%\end{center}

%The supplementary document is organized as follows.
%\begin{itemize}
%\item
%In Appendix~\ref{Apendix definition}, we summarize the notations in the paper.
%\item 
%In Appendix~\ref{Technical Proofs}, we prove the technical proofs of all the results in the paper.
%\item
%{\color{blue}In Appendix~\ref{sec:technical:lemmas}, we present additional technical lemmas.}
%\end{itemize}

\appendix

\section{Notations}\label{Apendix definition}

~~~~~~
For a function $f:\mathbb{R}^{n}\rightarrow\mathbb{R}$, 
we introduce the following notations:
\begin{itemize}
\item
The gradient of $f$ is defined as
$\nabla f := \left( \frac{\partial f}{\partial x_1}, \frac{\partial f}{\partial x_2}, \dots, \frac{\partial f}{\partial x_n} \right)
$.
\item
$\tilde{\nabla}$ denote the adjoint of gradient operator $\nabla$ in $L^2(\nu)$.

\item
The Hessian matrix of $f$ is defined as
$\nabla^2 f:=\left( \frac{\partial^2 f}{\partial x_i \partial x_j } \right)_{i,j=1}^n$. 
%\[
%\nabla^2 f = \begin{pmatrix}
%\frac{\partial^2 f}{\partial x_1^2} & \frac{\partial^2 f}{\partial x_1 \partial x_2} & \cdots & \frac{\partial^2 f}{\partial x_1 \partial x_n} \\
%\frac{\partial^2 f}{\partial x_2 \partial x_1} & \frac{\partial^2 f}{\partial x_2^2} & \cdots & \frac{\partial^2 f}{\partial x_2 \partial x_n} \\
%\vdots & \vdots & \ddots & \vdots \\
%\frac{\partial^2 f}{\partial x_n \partial x_1} & \frac{\partial^2 f}{\partial x_n \partial x_2} & \cdots & \frac{\partial^2 f}{\partial x_n^2}
%\end{pmatrix}
%\]
\item
The third-order tensor of partial derivatives
of $f$ is defined as
$\nabla^3 f := \left( \frac{\partial^3 f}{\partial x_i \partial x_j \partial x_k} \right)_{i,j,k=1}^n$.
\item
The Laplacian of $f$ is given by
$\Delta f := \nabla \cdot \nabla f = \sum_{i=1}^n \frac{\partial^2 f}{\partial x_i^2}$.
\end{itemize}

For a vector field $\mathbf{v}:\mathbb{R}^{n}\rightarrow\mathbb{R}^{n}$, 
its divergence is defined as
$\nabla \cdot \mathbf{v} :=\sum_{i=1}^{n} \frac{\partial v_i}{\partial x_i}$.

For an operator $\mathcal{L}:=b~\cdot\nabla+\mathcal{S}:\nabla^2$,  
the carr\'{e} du champ operator $\mathcal{C}$ associated with $\mathcal{L}$ is defined as follows; see e.g. \cite{BGL2013}. For two regular functions $\varphi$, $\psi$:
\[
\mathcal{C}(\varphi, \psi) := \frac{1}{2} \left( \mathcal{L}(\varphi \psi) - \varphi \mathcal{L}\psi - \psi \mathcal{L}\varphi \right) = \nabla \varphi \cdot \mathcal{S} \nabla \psi.
\]

Next, we define various spaces for real-valued functions on $\mathcal{X}$
that are used throughout the paper as follows.

\begin{itemize}

\item $C^{\infty}(\mathcal{X})$ is the space of smooth functions.
\item $C_c^{\infty}(\mathcal{X})$ is the space of smooth functions with compact support.
\item $C^k(\mathcal{X})$ is the space of functions are smooth to the $k$-th derivative.
\item $C_b(\mathcal{X})$ is the space of continuous and bounded functions.
\item $\mathscr{S}$ is the space of smooth functions growing at most polynomially and whose derivatives also grow at most polynomially:
\[
\mathscr{S} = \left\{ \varphi \in C^{\infty}(\mathcal{X}) \, \middle| \, \forall \alpha \in \mathbb{N}^d, \, \exists N > 0 \text{ such that } \sup_{x \in \mathcal{X}} \frac{|\partial^{\alpha} \varphi(x)|}{(1 + |x|^2)^N} < +\infty \right\},
\]
where $\partial^{\alpha} = \partial_{x_1}^{\alpha_1} \dots \partial_{x_d}^{\alpha_d}$ with $\alpha = (\alpha_1, \dots, \alpha_d)$.

\item $B_{\infty}(\mathcal{X})$ is the space of bounded measurable functions, that is endowed with the norm
\[
\|\varphi\|_{B_{\infty}} := \sup_{x \in \mathcal{X}} |\varphi(\mathcal{X})|.
\]

\item $B_{W}^{\infty}(\mathcal{X})$ is defined as
\[
B_{W}^{\infty}(\mathcal{X}) = \left\{ \varphi : \mathcal{X} \to \mathbb{R} \text{ measurable} \, \middle| \, \|\varphi\|_{B_{W}^{\infty}} := \sup_{x \in \mathcal{X}} \frac{|\varphi(x)|}{W(x)} < +\infty \right\}, 
\]
where $W : \mathcal{X} \to [1, +\infty)$ is any given measurable function.

\item $\mathcal{P}_W(\mathcal{X})$ is the associated space of probability measures for duality on measure spaces:
\[
\mathcal{P}_W(\mathcal{X}) = \left\{ \nu \in \mathcal{P}(\mathcal{X}) \, \middle| \, \nu(W) < +\infty \right\}.
\]
%The associated weighted total variation distance is:
%\[
%\forall \nu, \eta \in \mathcal{P}_W(X), \ d_W(\nu, \eta) = \sup_{\|\varphi\|_{B_{W}^{\infty}} \leq 1} \left( \int_X \varphi \, d\nu - \int_X \varphi \, d\eta \right) = \int_X W(x) |\nu - \eta|(dx),
%\]

\item $C^0([0,T],\mathbb{R}^d)$ refers to the set of continuous functions from the interval $[0,T]$ to $\mathbb{R}^d$. %The $C^0$ space is typically the space of continuous functions without requiring differentiability. 

\item  For any probability measure $\mu \in \mathcal{P}(\mathcal{X})$, let
\[
L^2(\mu) = \left\{ \varphi \ \text{measurable} \ \middle|\ \int_\mathcal{X}|\varphi|^2 d\mu < +\infty \right\}. 
\]

\item For any  $\varphi \in C_c^\infty(\mathcal{X})$, we introduce the seminorm
\[
|\varphi|_{\mathscr{H}^1(\mu)}^2 = \int_\mathcal{X} \mathscr{C}(\varphi, \varphi) d\mu, 
\]
and the equivalence relation $\sim_1$ through: $\varphi \sim_1 \psi$ if and only if $|\varphi - \psi|_{\mathscr{H}^1(\mu)} = 0$. 

\item $\mathscr{H}^1(\mu)$ is the closure of $C_c^\infty(\mathcal{X})$ quotiented by $\sim_1$ for the norm $|\cdot|_{\mathscr{H}^1(\mu)}$. Note that $\mathscr{H}^1(\mu)$ and $L^2(\mu)$ are not subspaces of each other in general, but $\mathscr{H}^1(\mu) \subset L^2(\mu)$ for instance if $\mu$ satisfies a Poincaré inequality and $\mathcal{S}$ is positive definite. The difference between $L^2(\mu)$ and $\mathscr{H}^1(\mu)$ is however important for degenerate dynamics; see the application in Section~4.2 in \cite{LDP-GG}. We now construct a space dual to $\mathscr{H}^1(\mu)$ with the same density argument by introducing the seminorm: for $\varphi \in C_c^\infty(\mathcal{X})$,
\[
|\varphi|_{\mathscr{H}^{-1}(\mu)}^2 = \sup_{\psi \in C_c^\infty(\mathcal{X})} \left\{ 2 \int_\mathcal{X} \varphi \psi d\mu - |\psi|_{\mathscr{H}^1(\mu)}^2 \right\}.
\]

We define similarly the equivalence relation $\sim_{-1}$ on $C_c^\infty(\mathcal{X})$ by $\varphi \sim_{-1} \psi$ if and only if $|\varphi - \psi|_{\mathscr{H}^{-1}(\mu)} = 0$. The space $\mathscr{H}^{-1}(\mu)$ is then the closure of $C_c^\infty(\mathcal{X})$ quotiented by $\sim_{-1}$. This is actually the dual space of $\mathscr{H}^1(\mu)$. 
\end{itemize}

Next, we introduce the definitions of various topologies that are used in our paper.

\begin{itemize}

\item $\tau$-topology is the weak topology on $\mathcal{P}(\mathcal{X})$ associated with the convergence of measures tested against functions belonging to $B_{\infty}(\mathcal{X})$ (we may also use the notation $\sigma(\mathcal{P}(\mathcal{X}), B_{\infty})$). This means that for a sequence $(\nu_n)_{n \in \mathbb{N}}$ in $\mathcal{P}(\mathcal{X})$, $\nu_n \to \nu$ in the $\tau$-topology if $\nu_n(\varphi) \to \nu(\varphi)$ for any $\varphi \in B_{\infty}(\mathcal{X})$. Recall that the $\tau$-topology is stronger than the usual weak topology $\sigma(\mathcal{P}(\mathcal{X}), C_b(\mathcal{X}))$ on $\mathcal{P}(\mathcal{X}$), which corresponds to the convergence $\nu_n(\varphi) \to \nu(\varphi)$ for any $\varphi \in C_b(\mathcal{X})$. 

\item $\tau_W$ is the associated topology $\sigma(\mathcal{P}_W(\mathcal{X}), B_{W}^{\infty}(\mathcal{X}))$, accounting for the convergence of measures tested against the larger class of functions $\varphi \in B_{W}^{\infty}(\mathcal{X})$. 

\item The weighted total variation distance between two probability measures $\nu,\eta\in\mathcal{P}_{W}(\chi)$ is given by
\[
d_W(\nu, \eta) = \sup_{\|\varphi\|_{B^\infty_W} \leq 1} \left\{ \int_{\mathcal{X}} \varphi \, d\nu - \int_{\mathcal{X}} \varphi \, d\eta \right\} = \int_{\mathcal{X}} W(x) \, |\nu - \eta| (dx),\]
where $|\nu - \eta|$ is the total variation measure associated with 
$\nu - \eta$.
\end{itemize}

Finally, the notation \([ \cdot, \cdot ]\) refers to the Lie bracket of two vector fields. The Lie bracket is used to describe the symmetry or commutation properties between two vector fields and is defined as follows:
\[[A_i, A_j](f) = A_i(A_j(f)) - A_j(A_i(f)),\]
where \(f\) is a smooth function. Intuitively, the Lie bracket measures the non-commutativity of two vector fields, i.e., the difference in their actions when applied in different orders.

%%%%%%%%%%%%%%%%%%%%%%%%%%%%%%%%%%%%%%%%%%%%%%%%%%
\section{Technical Proofs}\label{Technical Proofs}		

\subsection{Proof of Theorem~\ref{thm:main}}\label{app:pf-thm:main}
\begin{proof} 
Under Assumption \ref{Hypoellipticity}-\ref{Lyapunov condition}, the empirical measure $\pi_t$ of the generalized Langevin dynamics defined in (\ref{APP}) satisfies a large deviation principle in the $\tau^\kappa$-topology, and the proof is similar as Theorem~2.11 in \cite{LDP-GG} and is omitted here. In the rest of the proof, we will obtain the rate function from the large deviation principle through the decomposition of the infinitesimal generator $\mathcal{L_\tau}$ of the dynamics (\ref{APP}) which can be easily computed as
\begin{equation}\label{FF}
\mathcal L_\tau=\nabla^{*}(\mathcal D+\mathcal Q)\nabla
=\underbrace{\nabla^{*}\mathcal D\nabla}_{\mathcal L_S}
+\underbrace{\nabla^{*}\mathcal Q\nabla}_{\mathcal L_A},
\end{equation}
where $\mathcal D(z)=\mathcal D(z)^{\mathsf T}\ge0,\mathcal Q(z)=- \mathcal Q(z)^{\mathsf T}$.
Adjoint of the gradient in $L^2(\mu): \tilde{\nabla}= -\operatorname{div} + \nabla U.$ By the definition of adjoint operator, we have for any smooth \(g,h\),\begin{equation}\label{KK}
\int_{\mathcal{X}} g(\mathbf{z})\mathcal{L}_{\tau}h(\mathbf{z})\ d\mu=\int_{\mathcal{X}} h(\mathbf{z})\mathcal{L}^*_\tau g(\mathbf{z})\ d\mu,
\end{equation}
where $\mathcal{L}_{\tau}$ is the infinitesimal generator of the generalized Langevin dynamics \eqref{APP} that is given by formula (\ref{FF}).
By plugging formula (\ref{FF}) into the left hand side of (\ref{KK}), and taking $\mu(d\mathbf{z})=e^{-H(\mathbf{z})}d\mathbf{z}$, where $H(\mathbf{z})$ is given in \eqref{eqn:Hamiltonian}, we get
\begin{equation}\label{MM}
\begin{aligned}
\int g\,\mathcal L_\tau h\,d\mu
&= \int g\,\nabla^{*}(\mathcal D+\mathcal Q)\nabla h\,d\mu
 = \int \nabla g\cdot (\mathcal D+\mathcal Q)\nabla h\,d\mu \\
&= \int \nabla h\cdot (\mathcal D+\mathcal Q)^{\mathsf T}\nabla g\,d\mu
 = \int \nabla h\cdot (\mathcal D-\mathcal Q)\nabla g\,d\mu \\
&= \int h\,\nabla^{*}(\mathcal D-\mathcal Q)\nabla g\,d\mu .
\end{aligned}
\end{equation}
Hence,
\[
\mathcal L_\tau^{*} \;=\; \nabla^{*}(\mathcal D-\mathcal Q)\,\nabla.
\]
Next, it follows from equation \eqref{GG} that
\begin{equation}\label{BH}
\mathcal{L}_S=\nabla^{*}\mathcal D\nabla,
\qquad
\mathcal{L}_A=\nabla^{*}\mathcal Q\nabla.
\end{equation}
By applying Lemma~\ref{RF} and Corollary~3.4 in~\cite{LDP-GG}, the $\mathscr{H}^{-1}(\nu)$–norm admits the energy
representation associated with the Dirichlet form induced by $\mathcal D$; hence, for any measure $\nu\in\mathcal{P}_\kappa({\mathcal{X}})$ such that $d\nu=e^\upsilon\ d\mu$ with $\upsilon\in\mathscr{H}^1(\nu)$ and  $\mathcal{L}_A\upsilon\in\mathscr{H}^{-1}(\nu)$
\[
\bigl\lvert \mathcal L_A(\upsilon)\bigr\rvert_{\mathscr{H}^{-1}(\nu)}^{2}
 \;=\; \int \nabla\psi_\upsilon \cdot \mathcal D\,\nabla\psi_\upsilon \, d\nu,
\]
where $\psi_\upsilon$ solves $-\operatorname{div}(\upsilon\,\mathcal D\,\nabla\psi_\upsilon)=\mathcal L_A(\upsilon)$.
This is consistent with the carré-du-champ identity in Eq.~(\ref{Gamma function}),
$\mathcal C(\varphi,\psi)=\nabla\varphi\cdot \mathcal S\,\nabla\psi$. Then we obtain
$$
I_\tau(\nu)=\frac{1}{4}\int_{\mathcal{X}}\nabla\upsilon\cdot \mathcal{D}\nabla\upsilon\ d\nu\ +\frac{1}{4}\int_{\mathcal{X}}\nabla\psi_\upsilon\cdot \mathcal{D}\nabla\psi_\upsilon\ d\nu,
$$
where $\psi_\upsilon$ is the unique solution in $\mathscr{H}^{1}(\nu)$ to the Poisson equation
$$
\tilde{\nabla}(\mathcal{D}\nabla\psi_\upsilon)=\mathcal{L}_A\upsilon.
$$
This completes the proof of Theorem~\ref{thm:main}.
\end{proof}
%%%%%%%%%%%%%%%%%%%%%%%%%%%%%%%%%%%%%%%%%%%%%%%%%%%%%%%%%%%%

%%%%%%%%%%%%%%%%%%%%%%%%%%%%%%%%%%%%%%%%%%%%%%%%%%%%%%%%%%%%
\subsection{Proof of Theorem~\ref{Mirror}\label{app:pf-Mirror}}

\begin{proof}
We know that the generator of the mirrored Langevin dynamics \eqref{mirror Langevin dynamics} is elliptic and therefore it automatically satisfies the hypoellipticity of the generator and the controllability (i.e., irreducibility of the dynamics) properties. By letting $\mathcal{D}=[\nabla^2\phi(\mathbf{z})]^{-1}$, the Lyapunov condition is satisfied by applying Lemma~\ref{Lyapunovkappa}. Finally, by applying Theorem~\ref{thm:main}, the conclusion follows.
\end{proof}

%%%%%%%%%%%%%%%%%%%%%%%%%%%%%%%%%%%%%%%%%%%%%%%%%%%%%%%%%%%%%%%%%%%%%%%%%%%%

\subsection{Proof of Lemma~\ref{Hypoellipticity-Controllability}\label{app:pf-Hypoellipticity-Controllability}}

\begin{proof}
The Hypoellipticity can be established as in Section~2 in \cite{2011OttobrePavliotis}, and hence we omit the proof here.
Next, we show the Controllability. Given $t > 0$ and two pair of points
$(\theta_0,p_0,r_0)$ and $(\theta_t,p_t,r_t)$, let $\phi(s)$ be any $\mathcal{C}^3$ path in $\mathbb{R}^d$ which satisfy $\phi(0) = \theta_0$, $\phi(t) = \theta_t$,
$\frac{d}{dt}\phi(0)=p_0$ and $\frac{d}{dt}\phi(t)=p_t$, $\frac{d^2}{dt^2}\phi(0)=-\nabla U(\theta_0)+ \gamma r_0$ and $\frac{d^2}{dt^2}\phi(t)=-\nabla U(\theta_t)+ \gamma r_t$. Consider the control $u$ given by
\begin{equation}\label{High-order control}
 u_t=\frac{1}{\sqrt{2\alpha}}\left(\frac{1}{\gamma}\frac{d^3}{dt^3}\phi(t)+\frac{\alpha}{\gamma}\frac{d^2}{dt^2}\phi(t)+\left(\frac{1}{\gamma}\nabla^2U(\phi(t))+\gamma\right)\frac{d}{dt}\phi(t)+\frac{\alpha}{\gamma}\nabla U(\phi(t))\right).
\end{equation}
By definition, $(\phi(t), \frac{d}{dt}\phi(t), \frac{1}{\gamma}(\frac{d^2}{dt^2}\phi(t)+\nabla U(\phi(t)))$ is a solution of the control system with control $u_t$. 
Therefore the Controllability (Assumption~\ref{irreducibility}) is satisfied.
The proof is complete.
\end{proof}

\begin{remark} The $u_t$ drives the system from $(\theta_0,p_0,r_0)$ to $(\theta_t,p_t,r_t)$. We denote by $G_t(x)$ the set of accessible points from $x$ in time $t$.
This implies that $G_t(x)= \mathbb{R}^{3d}$ for all $t > 0$ and all $x \in \mathbb{R}^{3d}$. From the support theorem (Theorem 6.1 in \cite{Bellet2006}) we conclude that $P_t(x,F) > 0$ for
all $t > 0$, all $x \in \mathbb{R}^{3d}$, and all open set $F$.
\end{remark}

%%%%%%%%%%%%%%%%%%%%%%%%%%%%%%%%%%%%%%%%%%%%%%%%%%%%%%%%%%%%%%%%%%%%%%%%%%%
\subsection{Proof of Proposition~\ref{Witten-Lyapunov condition-high:order}\label{app:pf-Witten-Lyapunov condition-high:order}}

\begin{proof}
Let us construct the Lyapunov functions $W_\delta : \mathbb{R}^{3d} \to [1, +\infty)$ satisfying  Lyapunov condition (\ref{LPC}). The construction is inspired from Proposition~2.13 in \cite{2024WuGullin}. We define the vector field $L$ as follows. Let $\chi : \mathbb{R}^d \to [0,1]$ be a smooth function such that $\chi(\theta) = 0$ if $|\theta| \leq 1$ and $\chi(\theta) = 1$ if $|\theta| \geq 2$. We define $J(\theta) := \theta |\theta|^{\beta - 1} \chi(\theta)$ with $\beta \in [0,1]$. Note that $J$ is $C^1$ and the first derivatives of $J$ are bounded over 
$\mathbb{R}^d$ (because $\beta \leq 1$), say by $C_J := \sup_{\theta\in\mathbb{R}^{d}} \|\text{Jac}(J(\theta))\|_2 > 0$ (where $\|M\|_2 := \sup\{|My|, |y| = 1\}, M \in \mathcal{M}_d(\mathbb{R})$). One then sets:
\begin{equation}\label{1}
L:= \kappa J, \quad \kappa := \frac{\gamma}{2C_J},
\end{equation}
so that
\begin{equation}\label{2}
C_L := \sup_{\theta\in\mathbb{R}^d} \|\text{Jac}(L(\theta))\|_2 \leq \gamma/2. 
\end{equation}

For all $(\theta,p,r) \in \mathbb{R}^{3d}, b \geq 0$, and $h,a > 0$, we define:
\begin{equation}\label{3}
\varphi_0(\theta,p,r) := hH(\theta,p,r) + a L(\theta)\cdot p+bp\cdot r,
\end{equation}
with $H(\theta,p,r)$ given in \eqref{H:function:high:order}.
The parameter $\beta > 0$ will be chosen such that 
\begin{equation}\label{4}
\inf_{(\theta,p,r)\in\mathbb{R}^{3d}} \varphi_0(\theta,p,r) \in \mathbb{R}.
\end{equation} 
We then set $\varphi_{\mathrm{HL}}(\theta,p,r) := \varphi_0(\theta,p,r) - \inf_{(\theta,p,r)\in\mathbb{R}^{3d}} \varphi_0(\theta,p,r) + 1$ and
\begin{equation}\label{11}
W_\delta(\theta,p,r) := \exp\left(\varphi_{\mathrm{HL}}^\delta(\theta,p,r)\right), \quad \text{where} \quad 1 - \frac{\beta}{k} < \delta < 1.
\end{equation}

In the following, for ease of notation we simply write $\varphi$ for $\varphi_{\mathrm{HL}}$. Since $\varphi^{1-\delta} \geq 1$, a straightforward computation implies that over $\mathbb{R}^{3d}$,
\begin{equation}\label{5}
\frac{\mathcal{L}_{H} W_\delta}{W_\delta} \le \frac{\delta}{ \varphi^{1-\delta}} \left[\mathcal{L}_{H} \varphi + \delta \alpha | \nabla_{r}\varphi|^2 \right], 
\end{equation}
where
\begin{equation}\label{High-order generator1}
\mathcal{L}_{H}=\alpha\Delta_r- (\gamma p+\alpha r)\cdot\nabla_r  + p\cdot\nabla_\theta -( \nabla U - \gamma r)\cdot\nabla_p.
\end{equation}
We also have
\begin{align*}
&\nabla_\theta \varphi(\theta,p,r) = h\nabla U(\theta) + a \text{Jac}(L(\theta)) p, \\
&\nabla_p \varphi(\theta,p,r) = hp +aL(\theta) +br, 
\\
&\nabla_r \varphi(\theta,p,r) = hr + bp.
\end{align*}
Consequently, one has for all $(\theta,p,r) \in \mathbb{R}^{3d}$,
\begin{align}\label{9}
\mathcal{L}_{H} \varphi(\theta,p,r) &= - \gamma b |p|^2 + ap \cdot \text{Jac}(L(\theta)) p - a \nabla U(\theta) \cdot L(\theta)+( b\gamma- \alpha h)|r|^2\nonumber\\
                              &\qquad\qquad+a\gamma r\cdot L(\theta)-b\nabla U(\theta) \cdot r-\alpha br\cdot p+hd\alpha, 
\end{align}
and 
\begin{align}\label{10}
\delta \alpha |\nabla_{\mathbf{r}} \varphi|^2(\theta,p,r)= \delta \alpha | h r + bp |^2,
\end{align}
with $\mathbf{z}=(\theta,p,r)$.
Set 
\begin{equation}\label{6}
b=a>0,~~~ k\in(1,2],~~~ \text{and}~~~ \beta=k-1.
\end{equation}  

Let us first check (\ref{4}). Let $p_1 = k / (k - 1) > 1$ and $q_1 = p_1 / (p_1 - 1) = k \leq 2$. Using Assumption \ref{HHJ}, we have for all $(\theta,p,r) \in \mathbb{R}^{3d}$, if $|\theta| \geq c_{U}$,
\begin{equation}\label{7}
\begin{aligned}
\varphi_0(\theta,p,r)&=h\left(U(\theta)+ \frac12|p|^2
+ \frac12|r|^2\right)+ a\kappa\theta |\theta|^{k-2} \chi(\theta)\cdot p+ap\cdot r\\
&\ge h m_U|\theta|^{k}+ \frac{h}{2}|p|^2+ \frac{h}{2}|r|^2-\frac{a\kappa}{k}|p|^k -\frac{a\kappa}{p_1}|\theta|^k- \frac{a}{2}|p|^2- \frac{a}{2}|r|^2\\
&\geq m_U h |\theta|^k + \frac{h - a}{2} |p|^2 + \frac{h - a}{2} |r|^2 - \frac{a\kappa}{p_1} |\theta|^k - \frac{a\kappa}{k} |p|^k.
\end{aligned}
\end{equation}
For the second inequality, we apply Cauchy-Schwarz and Young's inequalities to obtain:
\[
\begin{aligned}
|\theta|^{k-2}\theta\cdot p
&\ge -|\theta|^{k-2}\,|\theta|\,|p|
= -|\theta|^{k-1}\,|p|
\qquad(\text{Cauchy-Schwarz inequality})\\[2mm]
&\ge -\frac{\big(|\theta|^{k-1}\big)^{\,p_1}}{p_1}
   - \frac{|p|^{k}}{k}
\qquad\ (\text{Young's inequality with }p_1=\tfrac{k}{k-1},\,q=k)\\[2mm]
&=-\frac{|\theta|^{k}}{p_1}-\frac{|p|^{k}}{k}.
\end{aligned}
\]
Then, for $h > 0$, choose $a > 0$ small enough such that
\begin{equation}\label{8}
\frac{a\kappa}{p_1} < m_U h \quad \text{and} \quad \frac{a\kappa}{k}+\frac{a}{2} < \frac{h}{2}. 
\end{equation}
Then (\ref{4}) holds. Note also that when (\ref{6}) is satisfied,
\[
\varphi_0(\theta,p,r) \leq c'(|\theta|^k + |p|^2 + |r|^2) + C'.
\]

Next, let us verify that a Lyapunov function $W_\delta$ of the form (\ref{11}) which satisfies the Lyapunov condition (\ref{LPC}). Recall $b = a$. In the following, $(\theta,p,r) \in \mathbb{R}^{3d}$, $|\theta| \geq \max(2, c_U)$, and $\eta = \sqrt{a}$. Using Assumption \ref{HHJ}, (\ref{9}), and (\ref{10}), one has:
\begin{align*}
&\big(\mathcal{L}_{H} \varphi +  \delta \alpha |\nabla_r F|^2\big)(\theta,p,r) \\
&\leq \alpha h d - \gamma a |p|^2 + ap \cdot \text{Jac}(L(\theta)) p - a \nabla U(\theta) \cdot L(\theta) - \alpha h |r|^2 + \gamma a|r|^2 \\
&\quad + \gamma a|r||L(\theta)| + a |\nabla U(\theta)||r| + \alpha a|r||p| + \delta  \alpha |hr + ap|^2 \\
&\leq ah d - \gamma a |p|^2 + a C_L |p|^2 - a \kappa m_U |\theta|^{2(k-1)} - \alpha h |r|^2 + \gamma a|r|^2 \\
&\quad + \gamma \kappa a |r| |\theta|^{k-1} + a M_U |\theta|^{k-1} |r| + \alpha a |p| |r| + \delta \alpha |hr + ap|^2 \\
&\leq ah d - \gamma a |p|^2 + a C_L |p|^2 - a \kappa m_U |\theta|^{2(k-1)} - \alpha h |r|^2 + \gamma a|r|^2+ \frac{\gamma \kappa a |r|^2}{2\eta} + \gamma \kappa \eta a |\theta|^{2(k-1)} \\
&\quad +  \frac{a M_U |r|^2}{2\eta} + a\eta M_U |\theta|^{2(k-1)} + \frac{\alpha a|r|^2}{2\eta} + \frac{\alpha a\eta |p|^2}{2}+ 2\delta \alpha h^2 |r|^2 +2\delta \alpha a^2 |p|^2 \\
&\leq ah d + |p|^2 \left(-\gamma a + a C_L + 2 \delta \alpha a^2 +\alpha a^{\frac{3}{2}}/2\right)\\
&\qquad + |\theta|^{2(k-1)} \left(-\kappa m_U a + \gamma \kappa a^{3/2}/2 + M_U a^{3/2}/2\right) \\
&\qquad\qquad + |r|^2 \left(-\alpha h + 2 \delta\alpha h^2 + \gamma \kappa \sqrt{a}/2 + M_U \sqrt{a}/2 + \alpha \sqrt{a}/2 + \gamma a\right).
\end{align*}
With \eqref{11} and \eqref{6} we know that $\delta>\frac{1}{k}$. Let $h > 0$ such that $-\alpha h + 2 \delta h^2 < 0$. Using also (\ref{2}), it holds
$-\gamma + C_L \leq -\gamma / 2.$ We then choose $a > 0$ such that (\ref{8}) holds and
\begin{align*}
&-\gamma \alpha / 2 + 2 \delta a^2 \alpha + \alpha a^{3/2} / 2 < 0,
\\
&-\kappa m_U a + \gamma \kappa a^{3/2} / 2 + M_U a^{3/2} / 2 < 0,
\\
&-\alpha h + 2 \delta \alpha h^2 + \gamma \kappa \sqrt{a} / 2 + M_U \sqrt{a} / 2 + \alpha \sqrt{a} / 2 + \gamma a < 0.
\end{align*}
Then we proved inequality \eqref{LPC}. 

By the Witten identity, $\mathcal L_H W_\delta/W_\delta=\mathcal L_H\phi+\Gamma_H(\phi)$ with $\phi=\log W_\delta$ and $\Gamma_H(\phi)=\alpha|\nabla_r\phi|^2$.
Inequality \eqref{LPC} yields
\[
-(\mathcal L_H\phi+\Gamma_H(\phi))\ \ge\ A|\theta|^{2(k-1)}+B'|p|^2+C|r|^2-D.
\]
Since $W_\delta=e^{\varphi_{\mathrm{HL}}^\delta}$ with $\varphi_{\mathrm{HL}}=hH+a\,L(\theta)\!\cdot\!p+a\,p\!\cdot\!r$, one has on $K^c=\{\varphi_{\mathrm{HL}}\ge R\}$ that
$|\nabla_r\phi|^2\gtrsim |p|^2+|r|^2-C$ (choose $a>0$ small, $h>0$ large), and hence
\[
-(\mathcal L_H\phi+\Gamma_H(\phi))
\ \ge\ c_0\,\Gamma_H(\phi)-B\,\mathbf 1_K,
\]
for some $c_0>0$, $B<\infty$ and compact $K$. Equivalently,
\[
\mathcal L_H(\log W_\delta)+\Gamma_H(\log W_\delta)\ \le\ -\,c_0\,\Gamma_H(\log W_\delta)+B\,\mathbf 1_K,
\]
which is the Witten--Lyapunov inequality \eqref{eq:WL-high}.

It remains to show that the above estimate implies the Lyapunov condition
in Assumption~\ref{Lyapunov condition}. To avoid confusion of notation, set
\[
\Phi:=\varphi_{\mathrm{HL}},
\qquad
\psi:=\log W_\delta=\Phi^\delta .
\]
Recall that \(W_\delta=e^\psi\). By the Witten identity,
\[
\frac{\mathcal L_H W_\delta}{W_\delta}
=
\mathcal L_H\psi+\Gamma_H(\psi),
\qquad
\Gamma_H(\psi)=\alpha|\nabla_r\psi|^2 .
\]
From the estimate proved above, after increasing the constant if necessary,
there exist constants \(c,C>0\) such that
\[
\mathcal L_H\Phi+\delta\alpha|\nabla_r\Phi|^2
\le
C-c\Big(|\theta|^{2(k-1)}+|p|^2+|r|^2\Big).
\]
Using \eqref{5}, we obtain
\[
\frac{\mathcal L_H W_\delta}{W_\delta}
\le
\frac{\delta}{\Phi^{1-\delta}}
\left[
C-c\Big(|\theta|^{2(k-1)}+|p|^2+|r|^2\Big)
\right].
\]
Equivalently,
\[
-\frac{\mathcal L_H W_\delta}{W_\delta}
\ge
c_1
\frac{|\theta|^{2(k-1)}+|p|^2+|r|^2}{\Phi^{1-\delta}}
-C_1
\]
for some constants \(c_1,C_1>0\).

We now show that the right-hand side is coercive. Since \(\Phi\ge 1\) and,
by the construction of \(\Phi\),
\[
\Phi(\theta,p,r)
\le
C_2\bigl(1+|\theta|^k+|p|^2+|r|^2\bigr),
\]
we have
\[
\frac{|\theta|^{2(k-1)}+|p|^2+|r|^2}{\Phi^{1-\delta}}
\to +\infty,
\qquad |(\theta,p,r)|\to\infty .
\]
Indeed, in the \(\theta\)-direction the growth exponent is
\[
2(k-1)-k(1-\delta)
=
k\delta+k-2>0,
\]
because \(\delta>1/k\) and \(k>1\). In the \(p,r\)-directions, the quotient
grows like
\[
(|p|^2+|r|^2)^\delta\to+\infty .
\]
Therefore
\[
-\frac{\mathcal L_H W_\delta}{W_\delta}
\to+\infty,
\qquad |(\theta,p,r)|\to\infty .
\]
Hence there exist a compact set \(K\subset\mathbb R^{3d}\) and constants
\(\lambda,C_K>0\) such that
\[
\mathcal L_H W_\delta
\le
-\lambda W_\delta+C_K\mathbf 1_K .
\]
This is precisely the Foster--Lyapunov drift condition required in
Assumption~\ref{Lyapunov condition}.

Moreover, the same estimate also gives the Witten--Lyapunov form. Since
\[
\Gamma_H(\psi)
=
\alpha|\nabla_r\psi|^2
=
\alpha\delta^2\Phi^{2\delta-2}|\nabla_r\Phi|^2 ,
\]
and
\[
\nabla_r\Phi=hr+ap,
\]
we have
\[
\Gamma_H(\psi)
\le
C
\frac{|p|^2+|r|^2}{\Phi^{2(1-\delta)}}
\le
C
\frac{|\theta|^{2(k-1)}+|p|^2+|r|^2}{\Phi^{1-\delta}},
\]
where we used \(\Phi\ge1\). Combining this bound with the previous coercive
estimate, we obtain, for some \(c_0>0\), \(B<\infty\), and compact set \(K\),
\[
-\left(\mathcal L_H\psi+\Gamma_H(\psi)\right)
=
-\frac{\mathcal L_H W_\delta}{W_\delta}
\ge
c_0\Gamma_H(\psi)-B\mathbf 1_K .
\]
Equivalently,
\[
\mathcal L_H(\log W_\delta)+\Gamma_H(\log W_\delta)
\le
-c_0\Gamma_H(\log W_\delta)+B\mathbf 1_K .
\]This concludes the proof.
\end{proof}

%%%%%%%%%%%%%%%%%%%%%%%%%%%%%%%%%%%%%%%%%%%%%%%%%%%%%%%%%%%%
\subsection{Proof of Theorem~\ref{thm:high:order}\label{app:pf-thm:high:order}}

\begin{proof}
By Lemmas~\ref{Hypoellipticity-Controllability} and \ref{Witten-Lyapunov condition-high:order}, we know the infinite generator $\mathcal{L}_{H}$
satisfies the Hypoellipticity (Assumption~\ref{Hypoellipticity}), Controllability (Assumption~\ref{irreducibility}) and Witten-Lyapunov condition (Assumption~\ref{Lyapunov condition}). As a special case of the constructed model, by Theorem~\ref{thm:main}, the infinitesimal generator is decomposed 
such that the rate function is the sum of symmetric and anti-symmetric parts.
Through formula (\ref{High-order generator}) and (\ref{Gamma function}), we have
$$\mathscr{C}(\upsilon,\upsilon)=\frac{1}{2}\big(\mathcal{L}_H(\upsilon^2)-2\upsilon\mathcal{L}_H\upsilon\big)=\alpha|\nabla_r\upsilon|^2.$$
Hence by Theorem ~\ref{thm:main}, we obtain
$$
I_{H}(\nu)=\frac{\alpha}{4}\int_{\mathcal{X}}|\nabla_r\upsilon|^2\ d\nu + \frac{\alpha}{4}\int_{\mathcal{X}}|\nabla_r\psi|^2\ d\nu.
$$
Write $\psi_\upsilon=\alpha \psi$ where $\psi_\upsilon$ is the unique solution in $\mathscr{H}^{1}(\nu)$ to the Poisson equation (\ref{Poisson-equation-High-order})
and this completes of the proof of Theorem~\ref{thm:high:order}.
\end{proof}

%%%%%%%%%%%%%%%%%%%%%%%%%%%%%%%%%%%%%%%%%%%%%%%%%%%%%%%%%%%%%%%%%%%%%%%%%%%%

\subsection{Proof of Lemma~\ref{Witten-Lyapunov condition-Hessian-Free}\label{app:pf-Witten-Lyapunov condition-Hessian-Free}}
\begin{proof}
Note that $\mathcal{L}_R$ is an elliptic operator, so Controllability (Assumption~\ref{irreducibility}) is satisfied. Now we check the Witten-Lyapunov condition for the Hessian-free high-resolution dynamics \eqref{Hessian-free high-resolution}. 
Let
\[
\varphi(\theta,r)=aH(\theta,r)
=
a\left(U(\theta)+\frac12|r|^2\right),
\qquad
W_a=e^{\varphi}.
\]
we obtain
\[
\begin{aligned}
\frac{\mathcal L_R W_a}{W_a}
&=
(-\beta\nabla U(\theta)+r)\cdot a\nabla U(\theta)
+
(-\alpha r-\nabla U(\theta))\cdot ar
\\
&\qquad
+\beta\left(a\Delta U(\theta)+a^2|\nabla U(\theta)|^2\right)
+\alpha\left(ad+a^2|r|^2\right).
\end{aligned}
\]
Hence
\begin{equation}\label{eq:minusLW}
-\frac{\mathcal L_R W_a}{W_a}
=
a(1-a)\beta|\nabla U(\theta)|^2
+a(1-a)\alpha|r|^2
-a\beta\Delta U(\theta)
-a\alpha d.
\end{equation}

We now estimate the two terms involving \(U\). By Assumption~\ref{HFHR}(b),
\(\nabla U\) is globally \(m_U\)-Lipschitz. Since \(U\in\mathscr S\), this implies
\[
\|\nabla^2 U(\theta)\|_{\mathrm{op}}\le m_U,
\qquad \theta\in\mathbb R^d.
\]
Consequently,
\[
\Delta U(\theta)
=
\operatorname{tr}\nabla^2 U(\theta)
\le d m_U,
\qquad \theta\in\mathbb R^d,
\]
and hence
\begin{equation}\label{eq:Delta-bound}
-\Delta U(\theta)\ge -dm_U .
\end{equation}

Next, by Assumption~\ref{HFHR}(a),
\[
\theta\cdot\nabla U(\theta)
\ge c_1|\theta|^2-C_1,
\qquad \theta\in\mathbb R^d .
\]
Let
\[
C_1^+:=\max\{C_1,0\},
\qquad
R_U:=\max\left\{1,\sqrt{\frac{2C_1^+}{c_1}}\right\}.
\]
Then, for all \(|\theta|\ge R_U\),
\[
c_1|\theta|^2-C_1
\ge c_1|\theta|^2-C_1^+
\ge \frac{c_1}{2}|\theta|^2.
\]
Therefore,
\[
\theta\cdot\nabla U(\theta)
\ge \frac{c_1}{2}|\theta|^2,
\qquad |\theta|\ge R_U.
\]
On the other hand, by the Cauchy--Schwarz inequality,
\[
\theta\cdot\nabla U(\theta)
\le |\theta|\,|\nabla U(\theta)|.
\]
Since \(R_U\ge1\), we may divide by \(|\theta|\) for \(|\theta|\ge R_U\), and obtain
\[
|\nabla U(\theta)|
\ge \frac{c_1}{2}|\theta|,
\qquad |\theta|\ge R_U.
\]
Thus,
\[
|\nabla U(\theta)|^2
\ge \frac{c_1^2}{4}|\theta|^2,
\qquad |\theta|\ge R_U.
\]
It follows that the global lower bound
\begin{equation}\label{eq:gradU-lower-global}
|\nabla U(\theta)|^2
\ge
\frac{c_1^2}{4}|\theta|^2
-
\frac{c_1^2}{4}R_U^2,
\qquad \theta\in\mathbb R^d,
\end{equation}
holds. Indeed, when \(|\theta|\ge R_U\), this follows from the previous estimate; when
\(|\theta|<R_U\), the right-hand side of \eqref{eq:gradU-lower-global} is non-positive,
while the left-hand side is non-negative.

Combining \eqref{eq:minusLW}, \eqref{eq:Delta-bound}, and
\eqref{eq:gradU-lower-global}, we get
\[
\begin{aligned}
-\frac{\mathcal L_R W_a}{W_a}
&\ge
a(1-a)\beta
\left(
\frac{c_1^2}{4}|\theta|^2
-
\frac{c_1^2}{4}R_U^2
\right)
+
a(1-a)\alpha|r|^2
-a\beta d m_U
-a\alpha d
\\
&=
\frac{a(1-a)\beta c_1^2}{4}|\theta|^2
+
a(1-a)\alpha|r|^2
-
\left(
\frac{a(1-a)\beta c_1^2}{4}R_U^2
+a\beta d m_U
+a\alpha d
\right).
\end{aligned}
\]
Therefore \eqref{LPC1} holds with, for example,
\[
A=\frac{a(1-a)\beta c_1^2}{4},
\qquad
B=a(1-a)\alpha,
\]
and
\[
C=
\frac{a(1-a)\beta c_1^2}{4}R_U^2
+a\beta d m_U
+a\alpha d.
\]
Since \(a\in(0,1)\) and \(\alpha,\beta>0\), we have \(A,B>0\). Thus
\(W_a\) satisfies the desired Witten--Lyapunov condition.

Finally, since \(U\) has compact level sets and
$H(\theta,r)=U(\theta)+\frac12|r|^2,$ the function \(H\) has compact level sets in \(\mathbb R^d\times\mathbb R^d\).
Hence \(W_a=e^{aH}\to\infty\) at infinity. Therefore \(W_a\) is a Lyapunov
function. This completes the proof.
\end{proof}
%%%%%%%%%%%%%%%%%%%%%%%%%%%%%%%
\subsection{Proof of Theorem~\ref{thm: Hessian-Free}\label{app:pf-thm: Hessian-Free}}

\begin{proof}
From formula (\ref{Gamma function}), we have
\begin{align*}
&\mathscr{C}(\upsilon,\upsilon)=\frac{1}{2}\big(\mathcal{L}_R(\upsilon^2)-2\upsilon\mathcal{L}_R\upsilon\big)=\beta|\nabla_\theta\upsilon|^2 + \alpha|\nabla_r\upsilon|^2,
\\
&\mathscr{C}(\psi_\upsilon,\psi_\upsilon)=\beta|\nabla_\theta\psi_\upsilon|^2 + \alpha|\nabla_r\psi_\upsilon|^2,
\end{align*}
where $\mathcal{L}_R$ is defined in \eqref{BCB}.
For $\varphi\in C_c^\infty(\mathcal{X})$, we have the seminorm
$$
|\varphi|_{\mathscr{H}^1(\nu)}^2=\int_\mathcal{X}\mathscr{C}(\varphi,\varphi)\ d\nu=\int_\mathcal{X}\left(\beta|\nabla_\theta\varphi|^2 + \alpha|\nabla_r\varphi|^2\right)\ d\nu,
$$
and
$$
|\varphi|_{\mathscr{H}^{-1}(\nu)}^2=\sup_{\psi\in C_c^\infty}\left\{2\int_\mathcal{X}\varphi\psi\ d\nu-|\psi|_{\mathscr{H}^1(\nu)}^2\right\}.
$$
By Theorem ~\ref{thm:main}, we deduce that 
\begin{align*}
I_A(\nu)&=\frac{1}{4}\big|\mathcal{L}_{\mathrm{RA}}(\upsilon)\big|_{\mathscr{H}^{-1}(\nu)}^2\\
        &=\frac{1}{4}\sup_{\psi\in C_c^\infty}\left\{2\int_\mathcal{X}\mathcal{L}_{\mathrm{RA}}(\upsilon)\psi\ d\nu-|\psi|_{\mathscr{H}^1(\nu)}^2\right\}\\
        &=-\frac{1}{2}\inf_{\psi\in C_c^\infty}\left\{\int_\mathcal{X}\frac{1}{2}\mathscr{C}(\psi,\psi)-\mathcal{L}_{\mathrm{RA}}(\upsilon)\psi\ d\nu\right\},
\end{align*}
where $\mathcal{L}_{\mathrm{RA}}$ is the anti-symmetric part of $\mathcal{L}_{R}$.
Now it comes to a variational problem to find the infimum of $\int_\mathcal{X}\left(\frac{1}{2}\mathscr{C}(\psi,\psi)-\mathcal{L}_A(\upsilon)\psi\right)d\nu$.
Using Euler-Lagrange equation and by Lemma~\ref{RF} we can get that 
\begin{equation}\label{rate function Hessian Free1}
I_R(\nu)=\frac{\beta}{4}\left(\int_{\mathcal{X}}|\nabla_\theta\upsilon|^2\ d\nu + \int_{\mathcal{X}}|\nabla_\theta\psi_\upsilon|^2\ d\nu\right)
+ \frac{\alpha}{4}\left(\int_{\mathcal{X}}|\nabla_r\upsilon|^2\ d\nu + \int_{\mathcal{X}}|\nabla_r\psi_\upsilon|^2\ d\nu\right),
\end{equation}
where $\psi_{\upsilon}$ is the unique solution in $\mathscr{H}^{1}(\nu)$ to the Poisson equation
\begin{equation}\label{Poisson-equation-Hessian-Free1}
\tilde{\nabla}_{\nu}\!\left(
\begin{pmatrix}\beta\,\nabla_{\theta}\psi_{\upsilon}\\[2pt]
\alpha\,\nabla_{r}\psi_{\upsilon}\end{pmatrix}
\right)
=\mathcal{L}_{\mathrm{RA}}\upsilon,
\end{equation}
where $\tilde{\nabla}_{\nu}$ denotes the $L^{2}(\nu)$–adjoint of the gradient, i.e. the
weighted divergence:
\[
\tilde{\nabla}_{\nu}g=-\,\mathrm{div}_{\nu} g
:=-\,e^{\upsilon}\,\nabla\!\cdot\!\big(e^{-\upsilon} g\big)
= -\big(\nabla\!\cdot g+\langle\nabla\upsilon,\,g\rangle\big),
\qquad g\in L^{2}(\nu;\mathbb{R}^{2d}).
\]
Equivalently,
\[
-\,e^{\upsilon}\,\nabla\!\cdot\!\left(
e^{-\upsilon}
\begin{pmatrix}\beta\,\nabla_{\theta}\psi_{\upsilon}\\
\alpha\,\nabla_{r}\psi_{\upsilon}\end{pmatrix}
\right)
=\mathcal{L}_{\mathrm{RA}}\upsilon,
\]
i.e. the Poisson equation \eqref{Poisson-equation-Hessian-Free}
and this completes of the proof of Theorem~\ref{thm: Hessian-Free}.
\end{proof}

%%%%%%%%%%%%%%%%%%%%%%%%%%%%%%%%%%%%%
\subsection{Proof of Corollary~\ref{thm:accelerate}\label{app:pf-thm:accelerate}}

\begin{proof}
It follows from Theorem~\ref{thm:main} that
$$
I_\tau(\nu)=\frac{1}{4}\int_{\mathcal{X}}\nabla\upsilon\cdot \mathcal{D}\nabla\upsilon\ d\nu\ +\frac{1}{4}\int_{\mathcal{X}}\nabla\psi_\upsilon\cdot \mathcal{D}\nabla\psi_\upsilon\ d\nu.
$$
Since $\mathcal{D}$ is a positive semidefinite diffusion matrix, we get
$\nabla\psi_\upsilon\cdot \mathcal{D}\nabla\psi_\upsilon\ge 0$, which implies
$$
\frac{1}{4}\int_{\mathcal{X}}\nabla\psi_\upsilon\cdot \mathcal{D}\nabla\psi_\upsilon\ d\nu\ge 0.
$$
It can be seen from \eqref{eq:Io-variational} that
$$
I_o(\nu)=\frac{1}{4}\int_{\mathcal{X}}|\nabla\upsilon|^2\ d\nu.
$$
Next, we define
$$
J_\tau(\nu):=I_\tau(\nu)-I_o(\nu)=\frac{1}{4}\int_{\mathcal{X}}\nabla\upsilon\cdot (\mathcal{D}-\mathbf{I})\nabla\upsilon\ d\nu\ +\frac{1}{4}\int_{\mathcal{X}}\nabla\psi_\upsilon\cdot \mathcal{D}\nabla\psi_\upsilon\ d\nu,
$$
where $\mathbf{I}$ is the identity matrix. If $\mathcal{D}-\mathbf{I}$ is a positive semidefinite diffusion matrix, by the definition of $J_{\tau}(\nu)$, it is clear that $J_\tau(\nu)\ge 0$. This completes the proof.
\end{proof}

%%%%%%%%%%%%%%%%%%%%%%%%%%%%%%%%%%%%%%%%%%%%%%%%%%%%%%%%%%%%%%%%%%%
\subsection{Proof of Corollary~\ref{Mirror VS overdamped}\label{app:pf-Mirror VS overdamped}}

\begin{proof} In this example just let $\mathcal{D}=[\nabla^2\phi(\mathbf{z})]^{-1}$, then this corollary is a direct application of Corollary~\ref{thm:accelerate}.
\end{proof}

%%%%%%%%%%%%%%%%%%%%%%%%%%%%%%%%%%%%%%%%%%%%%%%%%%%%%%%%%%%%%%%%%%%
\subsection{Proof of Proposition~\ref{Hessian-free high-resolution VS overdamped}\label{app:pf-Hessian-free high-resolution VS overdamped}}

\begin{proof} Recall that the LDP rate function for Hessian-free high-resolution dynamics is in equation (\ref{rate function Hessian Free}), that is 
$$I_R(\nu)=\frac{\beta}{4}\left(\int_{\mathcal{X}}|\nabla_\theta\upsilon|^2\ d\nu + \int_{\mathcal{X}}|\nabla_\theta\psi_\upsilon|^2\ d\nu\right)
+ \frac{\alpha}{4}\left(\int_{\mathcal{X}}|\nabla_r\upsilon|^2\ d\nu + \int_{\mathcal{X}}|\nabla_r\psi_\upsilon|^2\ d\nu\right).$$ 
Since  $\nabla\upsilon=(\nabla_\theta\upsilon, \nabla_r\upsilon)$, we have 
$$
\int_{\mathcal{X}}|\nabla\upsilon|^2\ d\nu=\int_{\mathcal{X}}|\nabla_\theta\upsilon|^2\ d\nu+\int_{\mathcal{X}}|\nabla_r\upsilon|^2\ d\nu. 
$$
Hence we conclude that if $\min(\alpha,\beta)\geq 1$, then we have $I_R(\nu)\ge I_{eo}(\nu)$. This completes the proof.
\end{proof}

%%%%%%%%%%%%%%%%%%%%%%%%%%%%%%%%%%%%%%%%%%%%%%%%%%%%%%%%%%%%%%%%%%%
\subsection{Proof of Corollary~\ref{cor:HFHR:vs:overdamped}}\label{app:pf-cor:HFHR:vs:overdamped}

\begin{proof}
This corollary follows immediately from Proposition~\ref{Hessian-free high-resolution VS overdamped}
and \eqref{I:R:theta}, \eqref{eq:Io-variational} and \eqref{I:o:contraction}.    
\end{proof}

%%%%%%%%%%%%%%%%%%%%%%%%%%%%%%%%%%%%%%%%%%%%%%%%%%%%%%%%%%%%%%%%%%%
\subsection{Proof of Proposition~\ref{underdamped VS overdamped}}\label{app:pf-underdamped VS overdamped}

\begin{proof} Recall that the LDP rate function for underdamped Langevin dynamics is in equation (\ref{rate function underdamped}), that is 
$$I_{u}(\nu) = \frac{\gamma}{4} \int_{\mathcal{X}} |\nabla_r \upsilon|^2 d\nu + \frac{1}{4\gamma} \int_{\mathcal{X}} |\nabla_r \psi|^2 d\nu.$$ 
Similar as in the proof of Proposition ~\ref{Hessian-free high-resolution VS overdamped} 
$$I_{e2o}(\nu)=\frac{1}{4}\int_{\mathcal{X}}|\nabla\upsilon|^2\ d\nu=\frac{1}{4}\int_{\mathcal{X}}|\nabla_\theta\upsilon|^2\ d\nu+\frac{1}{4}\int_{\mathcal{X}}|\nabla_r\upsilon|^2\ d\nu.$$ 
When $\upsilon=\upsilon(r)$, we have $\nabla_\theta\upsilon=0$. If $\gamma\ge1$, we conclude that $I_u(\nu)\ge I_{e2o}(\nu)$ and this completes the proof.
\end{proof}

%%%%%%%%%%%%%%%%%%%%%%%%%%%%%%%%%%%%%%%%%%%%%%%%%%%%%%%%%%%%%%%%%%%
\subsection{Proof of Proposition~\ref{High-order VS overdamped}}\label{app:pf-High-order VS overdamped}

\begin{proof} Recall the LDP rate function for high-order Langevin dynamics from equation (\ref{eqn: rate function High order}): 
$$I_{H}(\nu)=\frac{\alpha}{4}\int_{\mathcal{X}}|\nabla_r\upsilon|^2\ d\nu + \frac{1}{4\alpha}\int_{\mathcal{X}}|\nabla_r\psi|^2\ d\nu.$$
Moreover,
$$I_{e3o}(\nu)=\frac{1}{4}\int_{\mathcal{X}}|\nabla\upsilon|^2\ d\nu=\frac{1}{4}\int_{\mathcal{X}}|\nabla_\theta\upsilon|^2\ d\nu+\frac{1}{4}\int_{\mathcal{X}}|\nabla_p\upsilon|^2\ d\nu+\frac{1}{4}\int_{\mathcal{X}}|\nabla_r\upsilon|^2\ d\nu.$$ 
Thus, when $\upsilon=\upsilon(r)$, we have $\nabla_\theta\upsilon=0$ and $\nabla_p\upsilon=0$. If $\alpha\ge1$, we conclude that $I_{H}(\nu)\ge I_{e3o}(\nu)$, and this completes the proof.
\end{proof}

%%%%%%%%%%%%%%%%%%%%%%%%%%%%

\section{Additional Technical Lemmas}\label{sec:technical:lemmas}

We first introduce the following technical lemma from \cite{LDP-GG}, which is used in the proof of Theorem~\ref{Mirror} and in particular showing the Lyapunov condition. In the sequel, we assume Assumption~\ref{Lyapunov condition}. 

\begin{lemma}[Proposition~2.9 in \cite{LDP-GG}]\label{Lyapunovkappa} 
Suppose that Assumption~\ref{Lyapunov condition} holds. Then Assumption~\ref{Witten-Lyapunov condition} is satisfied with
\[
  W_\eta(\mathbf{z})=e^{\eta W(\mathbf{z})},\qquad \mathscr W_\varepsilon(x)=e^{\varepsilon W(x)},
\]
for $\theta\in(0,1)$ and $\varepsilon<\theta/2$ small enough. In this case, it holds that
\[
\Psi_\tau:=-\frac{\mathcal{L_\tau}W_\eta}{W_\eta}=\eta\left(-\mathcal{L_\tau}W-\eta|\sqrt{\mathcal{D}}\nabla W|^2\right)\sim|\sqrt{\mathcal{D}}\nabla W|^{2}.
\]
\end{lemma}

%It is a Lyapunov function of the generalized Langevin dynamics, assume that the process satisfies the following conditions:
%\begin{aligned}
%&\frac{\sum\limits_i\sum\limits_j\Big(\frac{\partial}{\partial\mathbf{z}_j}\mathcal{Q}_{ij}(\bm{\mathbf{z}})-\mathcal{Q}(\bm{\mathbf{z}})\frac{\partial}
%{\partial\mathbf{z}_j}H(\mathbf{z})\Big)
%\frac{\partial}{\partial\mathbf{z}_i}H(\mathbf{z})}{\Psi}\xrightarrow[|\mathbf{z}|\rightarrow +\infty]\quad 0,\\
%&\frac{\sum\limits_i\sum\limits_j\Big(\frac{\partial}{\partial\mathbf{z}_j}\mathcal{D}_{ij}(\bm{\mathbf{z}})\Big)\frac{\partial}{\partial\mathbf{z}_i}
%H(\mathbf{z})}{\Psi}\xrightarrow[|\mathbf{z}|\rightarrow +\infty]\quad 0.
%\end{aligned}
%%Moreover
%&=\sum\limits_i\sum\limits_j\eta\bigg[(1-\eta)\mathcal{D}(\mathbf{z})\Big(\frac{\partial}{\partial\mathbf{z}_i}H(\mathbf{z})\frac{\partial}
%{\partial\mathbf{z}_j}H(\mathbf{z})\Big)-\mathcal{D}(\mathbf{z})
%\Big(\frac{\partial^2}{\partial\mathbf{z}_i\partial\mathbf{z}_j}H(\mathbf{z})\Big)\\
%&\qquad\qquad-\Big(\frac{\partial}{\partial\mathbf{z}_j}\mathcal{Q}_{ij}(\bm{\mathbf{z}})
%-\mathcal{Q}(\bm{\mathbf{z}})\frac{\partial}{\partial\mathbf{z}_j}H(\mathbf{z})\Big)\frac{\partial}{\partial\mathbf{z}_j}H(\mathbf{z})
%-\Big(\frac{\partial}{\partial\mathbf{z}_j}\mathcal{D}_{ij}(\bm{\mathbf{z}})
%\Big)\frac{\partial}{\partial\mathbf{z}_i}H(\mathbf{z})\bigg]

Next, let us recall the LDPs for overdamped Langevin dynamics \eqref{eq:overdamped-2}, underdamped Langevin dynamics \eqref{eqn:underdamped} and non-reversible Langevin dynamics \eqref{non:reversible} from the literature. Setting $\mathbf{z}=\theta\in\mathbb{R}^d$ and choosing $H(\mathbf{z})=U(\theta)$,
we get the overdamped Langevin dynamics \eqref{eq:overdamped-2}.

\begin{lemma}[Proposition~4.2 in \cite{LDP-GG}]\label{LDP-A}
Assume the potential $U\in\mathscr S$ has compact level sets, $e^{-U}\in L^1(\mathcal X)$, and that for any $\theta\in(0,1)$,
\begin{equation}\label{eq:A41-core}
(1-\theta)\,|\nabla U|^2-\Delta U \xrightarrow[|x|\to\infty]{} +\infty .
\end{equation}
Let $W(x):=e^{\theta U(x)}$ and define
\[
\Psi:=-\,\frac{\mathcal LW}{W}
=\theta\Big((1-\theta)|\nabla U|^2-\Delta U\Big),
\]
where $\mathcal L$ is the generator of the process \eqref{eq:overdamped-2}. Then the following hold. 

\noindent\textbf{(i) Lyapunov function and contraction.}
$W$ is a Lyapunov function in the sense of Assumption~2.7; in particular,
$\Psi$ has compact level sets. Moreover, for any fixed $\theta\in(0,1)$ there
exist constants $C,c>0$ such that for any initial law $\nu\in\mathcal P_W(\mathcal X)$,
\[
d_W\big(\nu P_t,\mu\big)\;\le\; C\,e^{-ct}\, d_W(\nu,\mu),
\]
where $\mu$ is the invariant measure, $P_t$ is the semigroup, and $d_W$ is the
corresponding (weighted) Wasserstein distance.

\noindent\textbf{(ii) Large deviations.}
If $\kappa:\mathcal X\to[1,\infty)$ belongs to $\mathscr S$, is either bounded or has
compact level sets, and satisfies
\[
\frac{\Psi(x)}{\kappa(x)} \xrightarrow[|x|\to\infty]{} +\infty,
\]
then the empirical measure
\[
L_t:=\frac1t\int_0^t \delta_{X_s}\,ds
\]
satisfies a large deviation principle in the $\tau^\kappa$–topology and the good rate function $I_o(\cdot)$ defined in \eqref{TR} for the large deviations of overdamped Langevin dynamics \eqref{eq:overdamped-2} reads
\begin{equation}\label{eq:Io-variational}
I_o(\nu)=\frac14\int_{X}\!\!|\nabla \upsilon|^{2}\, \mathrm d\nu,
\qquad
\text{for any }\ \nu\in \mathcal P_\kappa(X) \ \text{of the form }\ \mathrm d\nu=e^{\upsilon}\,\mathrm d\mu,
\end{equation}
and $I_o(\nu)=+\infty$ otherwise.
\end{lemma}

The LDP for the underdamped Langevin dynamics \eqref{eqn:underdamped} has been obtained in e.g. \cite{LDP-GG}, which is stated as follows. 

\begin{lemma}[\textbf{Theorem~4.6} in \cite{LDP-GG}]\label{lem:underdamped} 
Assume that $(X_t)_{t \geq 0} = (\theta_t, r_t)_{t \geq 0}$ is the underdamped Langevin dynamics in equation \eqref{eqn:underdamped} where $U$ satisfies Assumption \ref{ass:potential} and consider a smooth function $\kappa$ with $\kappa(\theta, r) = 1 + |\theta|^{\alpha} + |r|^{\beta}$ for $|\theta| + |r| \geq 1$ and $\alpha \in [0, 2)$, $\beta \in [0, 2)$. Then $(X_t)_{t \geq 0}$ is ergodic with respect to the measure $\mu$. Moreover, the empirical measure
\[\pi_t := \frac{1}{t} \int_0^t \delta_{(\theta_s, r_s)} \, ds\]
satisfies a LDP in the $\tau^\kappa$-topology. Finally, for any $\nu \in \mathcal{P}_{\kappa}(\mathcal{X})$ such that $d\nu = e^\upsilon d\mu$ with $\upsilon \in \mathscr{H}^1(\nu)$ and $\mathcal{L}_{ham} \upsilon \in \mathscr{H}^{-1}(\nu)$,
where the generator of the dynamics is
\[
\mathcal{L}_\gamma = \mathcal{L}_{\text{ham}} + \gamma \mathcal{L}_{\text{FD}},
\]
with
\[
\mathcal{L}_{\text{ham}} = r \cdot \nabla_\theta - \nabla U \cdot \nabla_r, \quad \mathcal{L}_{\text{FD}} = -r \cdot \nabla_r + \Delta_r.
\]
and the rate function reads
\begin{equation}\label{rate function underdamped}
I_{u}(\nu) = \frac{\gamma}{4} \int_{\mathcal{X}} |\nabla_r \upsilon|^2 d\nu + \frac{1}{4\gamma} \int_{\mathcal{X}} |\nabla_r \psi|^2 d\nu, \end{equation}
and the formula for the rate function $I_u$ holds for probability measures $\nu$ of the form
$d\nu = e^{\upsilon}\,d\mu$ with $\upsilon \in \mathscr{H}^{1}(\nu)$, and $-J \nabla U \cdot \nabla \upsilon \in \mathscr{H}^{-1}(\nu)$, where $\psi$ is the unique solution in $\mathscr{H}^1(\nu)$ to the Poisson problem:
\begin{equation}-\Delta_r \psi + (r - \nabla_r \upsilon) \cdot \nabla_r \psi = \mathcal{L}_{ham}\upsilon. 
\end{equation}
\end{lemma}

The LDP for the non-reversible Langevin dynamics \eqref{non:reversible} has been obtained in e.g. \cite{LDP-GG}, which 
can be derived from Proposition~4.3 in \cite{LDP-GG}

\begin{lemma}[\textbf{Proposition~4.3} in \cite{LDP-GG}]\label{lemma:non:reversible}
Let $\theta_{t}$ follow the non-reversible Langevin dynamics \eqref{non:reversible}.
Assume that the potential \( U \in \mathscr{S} \) has compact level sets, satisfies \( e^{-U} \in L^1(\mathcal{X}) \) and, for any \( \eta \in (0,1) \), it holds
\[
(1 - \eta) |\nabla U|^2 - \Delta U \xrightarrow[|\theta| \to +\infty]{} +\infty.
\]
Then, with the notation in \eqref{GG}, it holds \( \mathcal{L}_{{S}} = -\nabla U\cdot \nabla + \Delta \) and \( \mathcal{L}_{{A}} = -J\nabla U \cdot \nabla \) and  \(\mathcal{L}=\mathcal{L}_S+\mathcal{L}_A.\)  Moreover, $\pi_{t}=\frac{1}{t}\int_{0}^{t}\delta_{\theta_{s}}ds$ satisfies a large deviation principle (LDP) in the \( \tau^\kappa \)-topology for any function $\kappa\in \mathscr{S}$ with compact level sets and such that \( \kappa(x) \to +\infty \) as \( |x| \to \infty \). The associated rate function \( I_J \) is given by:
for any probability measure $\nu$
such that $d\nu = e^{\upsilon}\,d\mu$ with $\upsilon \in \mathscr{H}^{1}(\nu)$ and
$(-J\nabla U) \cdot \nabla \upsilon \in \mathscr{H}^{-1}(\nu)$,
\[
I_J(\nu) = \frac{1}{2} \int_{\mathbb{R}^n} \left( |\nabla \psi_\upsilon|^2 +  |\nabla \upsilon|^2\right) \, d\nu,
\]
where \( \psi_\upsilon\) is the unique \( \mathscr{H}^1(\upsilon) \)-solution to the Poisson equation
\[
-\Delta \psi + \nabla (U -\upsilon) \cdot \nabla \psi_\upsilon= (-J\nabla U) \cdot \nabla \upsilon.
\]
\end{lemma}%\proofat{app:pf-lemma:non:reversible}

\begin{proof}
This is a direct consequence of Proposition~4.3 in \cite{LDP-GG}.
From that $\pi_{t}=\frac{1}{t}\int_{0}^{t}\delta_{\theta_{s}}ds$ satisfies a large deviation principle (LDP) in the \( \tau_\kappa \)-topology under the assumption that 
\( e^{-U} \in L^1(\mathcal{X}) \) and, for any \( \eta \in (0,1) \), it holds
\[
(1 - \eta) |\nabla U|^2 - \Delta U \xrightarrow[|\theta| \to +\infty]{} +\infty, 
\]
and \( \nabla \cdot \left(J\Delta U e^{-U}\right) = 0 \) and
\[
\frac{J\nabla U \cdot \nabla U}{\Psi} \xrightarrow[|\theta| \to +\infty]{} 0, 
\]
where \[
\Psi = \eta \left((1 - \eta)|\nabla U|^2 - \Delta U \right).
\]
Then, with the notation of \eqref{GG}, it holds \( \mathcal{L}_{\mathcal{S}} = -\nabla U\cdot \nabla + \Delta \) and \( \mathcal{L}_{\mathcal{A}} = -J\nabla U \cdot \nabla \) and  \(\mathcal{L}=\mathcal{L}_S+\mathcal{L}_A.\)  Moreover, we have the following relation:
\begin{align*}
- \frac{(\mathcal{L} -J \nabla U \cdot \nabla) e^{\eta U}}{e^{\eta U}} &= \eta \left((1 - \eta) |\nabla U|^2 - \Delta U +J \nabla U \cdot \nabla U\right)
\\
&=\eta \left((1 - \eta) |\nabla U|^2 - \Delta U \right)=\Psi.
\end{align*}
Finally, the conclusion follows from \( \nabla \cdot \left(J\Delta U e^{-U}\right) = 0 \) and $J\nabla U\cdot\nabla U=0$
since $J$ is an anti-symmetric matrix.
\end{proof}

\end{document}